\documentclass[11pt,a4paper]{article}
\usepackage{amsfonts}
\usepackage{amsmath,amssymb}
\usepackage{mathrsfs}
\usepackage{hyperref}
\usepackage{graphicx}
\usepackage{float}

\newtheorem{thm}{Theorem}[section]
\newtheorem{cor}[thm]{Corollary}
\newtheorem{lem}[thm]{Lemma}
\newtheorem{prop}[thm]{Proposition}

\numberwithin{equation}{section}\allowdisplaybreaks

\def\leq{\leqslant}

\def\leq{\leqslant}
\def\geq{\geqslant}

\def\Real{{\mathbb{R}}}

\begin{document}

\title{\large\bf  Local  Well-Posedness for the Derivative  Nonlinear \\ Schr\"{o}dinger Equations with $L^2$ Subcritical Data}

\author{\normalsize \bf Shaoming Guo$^{\dag}$,  \ Xianfeng Ren$^{\ddag}$  \  and Baoxiang Wang$^{\ddag,}$\footnote{Corresponding author. The project was supported in part by NSFC, grant 10271023} \\
\footnotesize
\it $^\dag$Department of Mathematics,   Indiana University, Bloomington, Indiana 47405, USA
\\
\footnotesize
\it $^\ddag$LMAM, School of Mathematical Sciences, Peking University, Beijing 100871, PR of China \\
\footnotesize
\texttt{ Emails: shaoguo@iu.edu, \ xianfengren@pku.edu.cn, \ wbx@math.pku.edu.cn}
} \maketitle

\thispagestyle{empty}
\begin{abstract}
Considering the Cauchy problem of the derivative
 nonlinear Schr\"{o}dinger equation (DNLS)
\begin{align}
{\rm i} u_{t} + \partial^2_{xx} u = {\rm i} \mu
 \partial_x (|u|^2u),\quad
u(0,x)=u_0(x),
 \label{gNLS}
\end{align}
we will show its local well-posedness in modulation spaces $M^{1/2}_{2,q}({\Real})$ $(2\leq q<\infty) $. It is well-known that $H^{1/2}$ is a critical Sobolev space of DNLS so that it is locally well-posedness in $H^s$ for $s\geq 1/2$ and ill-posed in $H^{s'}$ with $s'<1/2.$  Noticing that that $M^{1/2}_{2,q} \subset B^{1/q}_{2,q}$ is a sharp embedding  and $L^2 \subset B^0_{2,\infty}$,  our result contains all of the subcritical data in $M^{1/2}_{2,q}$, which contains a class of functions in $L^2\setminus H^{1/2}$. \\


\end{abstract}

\section{Introduction}

In this paper, we consider the Cauchy problem for the derivative
nonlinear Schr\"{o}dinger equations with a derivative nonlinearity
(DNLS)
\begin{align}
{\rm i} u_{t} + \partial^2_{xx} u = {\rm i} \mu
 \partial_x (|u|^2u),\quad
u(x,0)=u_0(x),
 \label{gNLS}
\end{align}
where   $u$ is a complex valued function of
$(x,t)\in\Real \times [0,T]$ for some $T>0$, $\mu\in \mathbb{R}$.

Using the gauge transform
\begin{align}
\mathcal{G} u (x)= \exp\left(-i \int^x_{-\infty} |u(y)|^2 dy \right) u(x) \label{gauge}
\end{align}
Hayashi \cite{Ha93} obtained the global well-posedness of DNLS in $H^1$ (see also \cite{HaOz94a,HaOz94b,OzTs98}). Takaoka \cite{Ta99} considered the rougher data and he established the local well-posedness in $H^s$ with $s \geq 1/2$ by considering the equivalent equation of $v=\mathcal{G}u$:
\begin{align}
{\rm i} v_{t} + \partial^2_{xx} v = - {\rm i} \mu v^2
 \partial_x  \overline{v} - \frac{|\mu|^2}{2} |v|^4 v,\quad
 v(x,0)=v_0(x).
 \label{egNLS}
\end{align}
 The ill-posedness in the sense that the solution map $u_0\to u(t)$ in $H^s$ with $s<1/2$ is not uniformly continuous has been obtained by Biagioni and Linares \cite{BiLi01}. So, $H^{1/2} $ is the critical Sobolev space in all $H^s$ so that DNLS is well posed.

However, the critical space for DNLS in the scaling sense  is $L^2$, i.e., for any solution $u$ of DNLS,  the scaling solution $u_\sigma (t,x) : = \sigma^{1/2} u(\sigma^2 t, \ \sigma x)$  has an invariant norm in $L^2$ for any $\sigma>0$. This fact implies that there is a gap between $L^2$ and $H^{1/2}$ for the well-posedness of DNLS. One can naturally ask what is the reasonable well-posed space with the regularity  at the same level with $L^2$. In order to answer this question, Gr\"unrock applied the  $\widehat{H^s_p}$ spaces for which the norm is defined by
$$
\|u\|_{\widehat{H^s_p}} : = \|\langle \xi\rangle^s \widehat{u}\|_{p'}, \ \ 1/p+1/p'=1
$$
and he obtained that DNLS is local well-posed in $\widehat{H^{1/2}_p}$ for any $1<p\leq 2$. Using the scaling argument we see that $\widehat{H^{1/2}_p}$ $(1<p\leq 2)$ can be regarded as subcritical spaces. In this paper we consider the initial data in more general modulation spaces.

We  write $\Box_k = \mathscr{F}^{-1} \chi_{[k-1/2, k+1/2]} \mathscr{F},$ where $\mathscr{F}$ ($\mathscr{F}^{-1}$) denotes the (inverse) Fourier transform on $ \mathbb{R}$ ,  $\chi_{E}$ denotes the characteristic function on $E$.   Modulation spaces $M^s_{p,q}$ were introduced by Feichtinger \cite{Fei2} and one can refer to \cite{Groch} for their basic properties. The modulation space $ M^s_{2,q}$
can be equivalently defined in the following way (cf.
\cite{WaHud07,WaHu07,WaHaHu09,WaHaHuGu11}):
\begin{align}
\|f\|_{M^s_{2,q}(\mathbb{R})}= \left(\sum_{k\in \mathbb{Z} } \langle k \rangle^{sq}
\|\Box_k f\|^q_{L^2(\mathbb{R})} \right)^{1/q},  \label{mod-space1}
\end{align}
where $\langle k\rangle=(1+|k|^2)^{1/2}$. Let $2\leq q\leq \infty$, $1/q+1/q'=1$.   By Plancherel's inequality and H\"older's inequality, we see that $\widehat{H^s_{q'}} \subset M^s_{2,q}$.  Combining the inclusions between Besov and modulation spaces, we have (cf. \cite{SuTo07,WaHaHuGu11})
$$
\widehat{H^{1/2}_{q'}} \subset M^{1/2}_{2,q} \subset B^{1/q}_{2,q}
$$
and these inclusions are optimal. On the other hand, from the scaling argument we see that for the scaling solution $u_\sigma$  (cf. Sugimoto and Tomita \cite{SuTo07}, and Han and Wang \cite{HaWa14}),
$$
 \|u_\sigma (t,\ \cdot) \|_{M^{1/2}_{2,q}}  \lesssim  \sigma^{1/q} \|u\|_{M^{1/2}_{2,q}}, \ \ 0< \sigma <1, \ \ q\geq 2.
$$
Hence, $M^{1/2}_{2,q}$ ($2\leq q <\infty$) can be regarded as subcritical spaces and $M^{1/2}_{2, \infty}$ is a critical space for the DNLS.

\begin{thm}
Let $2\leq q <\infty$,  $u_0 \in  M^{1/2}_{2,q}$. Then there exists a $T>0$ such that DNLS \eqref{gNLS} is locally well posed in $C([0,T];  M^{1/2}_{2,q}) \cap X^s_q([0,T])$, where $X^s_q$ is defined in \eqref{mod-space2}.
\end{thm}

The regularity index $1/2$ in $M^{1/2}_{2,q}$ is optimal. In fact, there is an ill-posedness for the DNLS in $M^s_{2,q}$ if $s<1/2$; cf. \cite{Wa13}. We conject that DNLS is also ill-posed in $M^{1/2}_{2,\infty}$. Let us observe the following examples:
$$
\widehat{u}_0 (\xi)= \langle \xi\rangle^{-1/2-\eta},  \ \ \widehat{v}_0 (\xi)=  \xi^{-1/2+ \theta} \chi_{(0,1)}(\xi).
$$
We see that $u_0\in L^2$ if and only if $u_0\in M^{1/2}_{2,q}$ for some $q>2$, if and only if $\eta>0$.   $v_0\in L^2$ if and only if  $v_0\in M^{1/2}_{2,q}$,  if and only if $\theta>0$. The above two examples indicate that the gaps between $L^2$ and $M^{1/2}_{2,q}$ $(q\gg 2)$ are very narrow.

There are some recent papers which have been devoted to the study of nonlinear PDE with initial data in modulation spaces $M^0_{p,1}$, see \cite{BO09,BGOR07,CFS12,CN08,CN08b,CN09,Iw10,KaKoIt12,KaKoIt14,Wa13,Wang1}. An interesting feature is that modulation spaces $M^0_{p,1}$ contains a class of initial data out of the critical Sobolev spaces $H^{s_c}$, for which the nonlinear PDE is well-posed for $s>s_c$ and ill-posed for $s<s_c$.  Guo \cite{Gu16} considered a class of initial data in $M_{2,q}$ for the cubic NLS, where the case $q>2$ was first taken into account by using $U^p$ and $V^p$ spaces.

Let $c < 1$, $C>1$  denote positive universal constants, which can
be different at different places, $a\lesssim b$ stands for $a\leq C
b$, $a\sim b$ means that $a\lesssim b$ and
$b\lesssim a$. If $a=b+\ell$, $|\ell| \leq C, $ then we write $a\approx b$. $a\gg b$ means that $a>b+C$.  We write $a\wedge b =\min(a,b)$, $a\vee b
=\max(a,b)$.   We write $p'$ as the dual number of $p \in
[1,\infty]$, i.e., $1/p+1/p'=1$.  We will use
Lebesgue spaces $L^p:=L^p(\mathbb{R} )$, $\|\cdot\|_p
:=\|\cdot\|_{L^p}$, Sobolev spaces $H^{s}=(\mathrm{I}-\Delta)^{-s/2}L^2$.
Some properties of these function spaces can be found in
\cite{Bergh,Triebel}.

\section{$U^p$ and $V^p$ spaces} \label{Upspaces}

\subsection{Definitions}

$U^p$ and $V^p$, as a development of Bourgain's spaces \cite{Bo93I,Bo93II} were first applied by Koch and Tataru in the study of NLS, cf. \cite{KoTa05,KoTa07,KoTa12}.  Using $U^p$ and $V^p$, Hadac, Herr and Koch \cite{HaHeKo09} obtained the well-posedness and scattering results for the critical KP-II equation.   Let $\mathcal{Z}$ be the set of finite partitions $-\infty= t_0 <t_1<...< t_{K-1} < t_K =\infty$. Let $1\leq p <\infty$. For any $\{t_k\}^K_{k=0} \subset \mathcal{Z}$ and $\{\phi_k\}^{K-1}_{k=0} \in L^2$ with $\sum^{K-1}_{k=0} \|\phi_k\|^p_2=1$, $\phi_0=0$. A step function $a: \mathbb{R}\to L^2$ given by
$$
a= \sum^{K}_{k=1} \chi_{[t_{k-1}, t_k)} \phi_{k-1}
$$
is said to be a $U^p$-atom. All of the $U^p$ atoms is denoted by $\mathcal{A}(U^p)$.   The $U^p$ space is
$$
U^p=\left\{u= \sum^\infty_{j=1} c_j a_j : \ a_j \in \mathcal{A}(U^p), \ \ c_j \in \mathbb{C}, \ \ \sum^\infty_{j=1} |c_j|<\infty  \right\}
$$
for which the norm is given by
$$
\|u\|_{U^p}= \inf \left\{\sum^\infty_{j=1} |c_j| : \ \ u= \sum^\infty_{j=1} c_j a_j , \ \  \ a_j \in \mathcal{A}(U^p), \ \ c_j \in \mathbb{C} \right\}.
$$

We define $V^p$ as the normed space of all functions $v: \mathbb{R} \to L^2$ such that $\lim_{t\to \pm \infty} v(t)$ exist and for which the norm
$$
\|v\|_{V^p} := \sup_{\{t_k\}^K_{k=0} \in \mathcal{Z}} \left( \sum^K_{k=1} \|v(t_k)-v(t_{k-1})\|^p_{L^2}\right)^{1/p}
$$
is finite, where we use the convention that $v(-\infty) = \lim_{t\to \infty} v(t)$ and $v(\infty)=0$ (here $v(\infty)$ and $\lim_{t\to  \infty} v(t)$ are different notations). Likewise, we denote by $V^p_-$ the subspace of all $v\in V^p$ so that $v(-\infty) =0$. Moreover, we define the closed subspace $V^p_{rc}$ $(V^p_{-,rc})$ as all of the right continuous functions in $V^p$ $(V^p_-)$.

We define
$$
U^p_{\Delta} := e^{\cdot \ {\rm i}\Delta} U^p, \ \ \|u\|_{U^p_{\Delta}} = \|e^{-{\rm i} t \Delta} u \|_{U^p}.
$$
$$
V^p_{\Delta} := e^{\cdot \ {\rm i}\Delta} V^p, \ \ \|u\|_{V^p_{\Delta}} = \|e^{-{\rm i} t \Delta} u \|_{V^p}.
$$
Similarly for the definition of $V^p_{rc, \Delta}$, $V^p_{-, \Delta}$, $V^p_{-, rc, \Delta}$.

We introduce the frequency-uniform localized $U^2_\Delta$-spaces $X^s_{q}(I)$ and $V^2_{\Delta}$-spaces $Y^s_{q} (I)$ for which the norms are defined by
\begin{align}
 \|u\|_{X^s_q(I)} & = \left(\sum_{\lambda \in I \cap \mathbb{Z}} \langle \lambda\rangle^{sq}\|\Box_\lambda u\|^q_{U^2}\right)^{1/q}, \quad  X^s_q:= X^s_q(\mathbb{R}),  \label{mod-space2}\\
 \|v\|_{Y^s_q(I)} & = \left(\sum_{\lambda \in I \cap \mathbb{Z}} \langle \lambda\rangle^{sq}\|\Box_\lambda v\|^q_{V^2}\right)^{1/q}, \quad  Y^s_q:= X^s_q(\mathbb{R}),  \label{mod-space2a} \\
\|u\|_{X^s_{q, \Delta}}: & = \| e^{-{\rm i} t \Delta}  u\|_{X^s_q}, \ \ \ \|v\|_{Y^s_{q, \Delta}}: = \| e^{-{\rm i} t \Delta}  v\|_{Y^s_q}
\end{align}
 Besov type Bourgain's spaces $\dot X^{s, b, q}$ are defined by
$$
\|u\|_{\dot X^{s,b,q}} := \left\| \|\chi_{|\tau-\xi^2|\in [2^{j-1}, 2^j)} |\xi|^{s} |\tau-\xi^2|^{b} \widehat{u}(\tau,\xi) \|_{L^2_{\xi,\tau}}  \right\|_{\ell^q_{j\in \mathbb{Z}}}.
$$

\subsection{Known results on $U^p$ and $V^p$}

We list some known results in $U^p$ and $V^p$ (cf. \cite{KoTa05,KoTa12,HaHeKo09}).

\begin{prop} \label{UVprop1}
{\rm (Embedding)} Let $1\leq p <q <\infty$. We have the following results.
\begin{itemize}
 \item[\rm (i)]  $U^p$ and $V^p$, $V^p_{rc}$, $V^p_{-}$, $V^p_{rc, -}$ are Banach spaces.

\item[\rm (ii)] $U^p\subset V^p_{rc, -} \subset U^q \subset L^\infty (\mathbb{R}, L^2)$. Every $u\in U^p$ is right continuous on $t\in \mathbb{R}$

\item[\rm (iii)] $V^p \subset V^q$,   $V^p_{-} \subset V^q_{-} $,   $V^p_{rc} \subset V^q_{rc} $,  $V^p_{rc, -} \subset V^q_{rc, -} $.

\item[\rm (iv)] $\dot X^{0, 1/2, 1} \subset U^2_{\Delta} \subset V^2_{\Delta} \subset \dot X^{0, 1/2, \infty}$.
\end{itemize}
\end{prop}
It is known that, for the free solution of the Schr\"odinger equation $u(x,t)= e^{{\rm i}t\Delta} u_0$, $\widehat{u}(\xi,\tau)$ is supported on a curve $\tau+\xi^2=0$, which is said to be the {\it dispersion relation}. For the solution $u$ of DNLS, $\widehat{u}(\xi,\tau)$ can be supported in $(\xi,\tau) \in \mathbb{R}^2$,  we need to consider the size of $|\xi^2+\tau|$, which is said to be the {\it dispersion modulation}. By the last inclusion of (iv) in Proposition \ref{UVprop1}, we see that have

\begin{lem}[\rm Dispersion Modulation Decay]
Suppose that the dispersion modulation $|\tau+\xi^2| \gtrsim \mu$ for a function $u\in L^2_{x,t}$, then we
\begin{align}
  \|u \|_{L^2_{x,t} }  \lesssim \mu^{-1/2} \|u\|_{V^2_\Delta}. \label{dispersiondecay}
\end{align}
\end{lem}

\begin{prop} \label{UVprop2}
{\rm (Interpolation)} Let $1\leq p <q <\infty$.  There exists a positive constant $\epsilon(p,q)>0$, such that for any $u\in V^p $ and $M>1$,  there exists a decomposition $u=u_1+u_2$ satisfying
\begin{align}
 \frac{1}{M} \|u_1\|_{U^p} + e^{\epsilon M}  \|u_2\|_{U^q} \lesssim \|u\|_{V^q}. \label{interp}
\end{align}
\end{prop}

\begin{prop} \label{UVprop3}
{\rm (Duality)} Let $1\leq p   <\infty$, $1/p+1/p'=1$.  Then $(U^p)^* = V^{p'}$ in the sense that
\begin{align}
T: V^{p'}  \to (U^p)^*; \ \ T(v)=B(\cdot,v), \label{dual}
\end{align}
is an isometric mapping.  The bilinear form $B: U^p\times V^{p'}$ is defined in the following way: For a partition $\mathrm{t}:= \{t_k\}^K_{k=0} \in \mathcal{Z}$, we define
 \begin{align}
B_{\mathrm{t}} (u,v) = \sum^K_{k=1} \langle u(t_{k-1}), \ v(t_k)-v(t_{k-1})\rangle. \label{dual2}
\end{align}
Here $\langle \cdot, \cdot \rangle$ denotes the inner product on $L^2$. For  any $u\in U^p$, $v\in V^{p'}$, there exists a unique number $B(u,v)$ satisfying the following property. For any $\varepsilon>0$, there exists a partition $\mathrm{t}$ such that
$$
|B(u,v)- B_{\mathrm{t}'} (u,v)| <\varepsilon, \ \ \forall \  \mathrm{t}'    \supset \mathrm{t}.
$$
Moreover,
$$
|B(u,v)| \leq \|u\|_{U^p} \|v\|_{V^{p'}}.
$$
In particular, let $u\in V^1_{-}$ be absolutely continuous on compact interval, then for any $v\in V^{p'}$,
$$
 B(u,v) =\int \langle u'(t), v(t)\rangle dt.
$$
\end{prop}

\subsection{Duality of $X^s_q$ and $Y^{-s}_{q'}$}

\begin{prop} \label{UVprop4}
{\rm (Duality)} Let $1\leq q   <\infty$.  Then $(X^s_q)^* = Y^{-s}_{q'} $ in the sense that
\begin{align}
T: Y^{-s}_{q'}   \to (X^s_q)^* ; \ \ T(v)=B(\cdot,v), \label{dualprop4}
\end{align}
is an isometric mapping, where the bilinear form $B(\cdot,\cdot)$ is defined in  Proposition \ref{UVprop3}. Moreover, we have
$$
|B(u,v)| \leq   \|u\|_{X^s_q} \|v\|_{Y^{-s}_{q'}}.
$$
\end{prop}
{\bf Proof.} By the orthogonality, we see that
$$
B_{\mathrm{t}}(\Box_k u, \Box_l v) =0, \ \ k\neq l,
$$
which implies that
$$
B (\Box_k u, \Box_l v) =0, \ \ k\neq l.
$$
For any $v\in  Y^{-s}_{q'}$, by Proposition \ref{UVprop3} and H\"older's inequality,  we have
\begin{align}
|B(u,v)| & = \left|\sum_{k\in \mathbb{Z}} B (\Box_k u, \Box_k v)\right| \nonumber\\
 & \leq  \sum_{k\in \mathbb{Z}} \|\Box_k u\|_{U^2} \|\Box_k v\|_{V^2}  \leq \|u\|_{X^s_q} \|v\|_{Y^{-s}_{q'}} \nonumber
\end{align}
It follows that $Y^{-s}_{q'} \subset (X^s_q)^*$ and $\|v\|_{(X^s_q)^*} \leq \|v\|_{Y^{-s}_{q'}}$.

Conversely,  considering the map
$$
X^s_q \ni f \to \{\Box_k f\} \in \ell^s_q (\mathbb{Z}; U^2),
$$
where
$$
\ell^s_q (\mathbb{Z}; U^2):= \left\{\{f_k\}_{k\in \mathbb{Z}}: \ \  \|\{f_k\}\|_{\ell^s_q (\mathbb{Z}; U^2)} := \left\| \{\langle k\rangle^s\|f_k\|_{U^2} \} \right\|_{\ell^q} <\infty \right\},
$$
we see that it is an isometric mapping from $X^s_q$ into a subspace of $\ell^s_q (\mathbb{Z}; U^2)$.  So,  $v \in (X^s_q)^*$ can be regarded as a continuous functional in a subspace of $\ell^s_q (\mathbb{Z}; U^2)$. In view of Hahn-Banach Theorem, it can be extended onto $\ell^s_q (\mathbb{Z}; U^2)$ (the extension is written as $\tilde{v}$) and its norm will be preserved.  In view of the well-known duality $(\ell^s_q (\mathbb{Z}; X))^* = \ell^{-s}_{q'} (\mathbb{Z}; X^*)$, we have
$$
(\ell^s_q (\mathbb{Z}; U^2))^* = \ell^{-s}_{q'} (\mathbb{Z}; V^2),
$$
and there exists $\{v_k\}_{k\in \mathbb{Z}} \in \ell^{-s}_{q'} (\mathbb{Z}; V^2)$ such that
$$
\langle \tilde{v}, \ \{f_k\}\rangle = \sum_{k\in \mathbb{Z}} B(f_k, v_k), \ \ \forall \ \{f_k\} \in \ell^s_q (\mathbb{Z}; U^2).
$$
Moreover, $\|v\|_{(X^s_q)^*}= \|\{v_k\}\|_{\ell^{-s}_{q'} (\mathbb{Z}; V^2)}.$  Hence,  for any $u\in X^s_q$,
$$
\langle v, u\rangle = \langle \tilde{v}, \{\Box_k u\} \rangle =  \sum_{k\in \mathbb{Z}} B(\Box_k u, v_k).
$$
From $B_{\mathrm{t}}(\Box_k u, \ v) = B_{\mathrm{t}}( u, \ \Box_k v)$ we see that $B (\Box_k u, \ v) = B ( u, \ \Box_k v)$. It follows that
$$
v= \sum_{k\in \mathbb{Z}} \Box_k v_k.
$$
Obviously, we have
$$
\|v\|_{Y^{-s}_{q'}} \leq \|\{v_k\}\|_{\ell^{-s}_{q'} (\mathbb{Z}; V^2)} = \|v\|_{(X^s_q)^*}.
$$
This proves $(X^s_q)^* \subset Y^{-s}_{q'}$.  $\hfill\Box$

Now we apply the duality to the norm calculation to the inhomogeneous part of the solution of DNLS in $X^s_{q,\Delta}$. It is known that \eqref{gNLS} is equivalent to the following integral equation:
\begin{align}
u(t)= e^{{\rm i}t\Delta} u_0 - \mathcal{A} \left( {\rm i} \mu u^2 \partial_x \overline{u} + \frac{|\mu|^2}{2}|u|^4 u \right),   \label{IDNLS}
\end{align}
where
$$
\mathcal{A} (f) = \int^t_0 e^{{\rm i}(t-s)\Delta} f(s) ds.
$$
By Propositions \ref{UVprop3} and \ref{UVprop4}, we see that, for ${\rm supp} \ v \subset \mathbb{R}\times [0,T]$,
\begin{align}
\|\mathcal{A} (f) \|_{X^{1/2}_{p,\Delta}} &  = \sup \left\{ \left|B\left(\int^t_0 e^{-{\rm i}s\Delta} f(s) ds, v \right) \right| : \  \|v\|_{Y^{-1/2}_{p'}} \leq 1 \right\} \nonumber\\
& \leq \sup_{\|v\|_{Y^{-1/2}_{p'}} \leq 1} \left|\int_{ [0,T]} \langle f (s), \  e^{{\rm i}s \Delta} v(s) \rangle ds  \right| \nonumber\\
& \leq \sup_{\|v\|_{Y^{-1/2}_{p', \Delta}} \leq 1} \left|\int_{ [0,T]} \langle f (s), \  v(s) \rangle ds  \right|.  \label{normest}
\end{align}

\section{Frequency localized estimates in  $L^4$}

Let $I\subset \mathbb{R}$ be an interval with finite length. For simply, we denote
$$
 u_\lambda = \Box_\lambda u, \ \  u_I = \sum_{\lambda\in I \cap \mathbb{Z}} u_\lambda,
$$

\begin{lem}\label{L4} {\rm (\cite{Gu16})}
Let $I\subset \mathbb{R}$ with $|I|<\infty$.  For  any  $\theta\in (0,1)$, $\beta>0$, we have
\begin{align}
 \|u_I\|_{L^4_{x,t\in [0,T]}} \lesssim (T^{1/8}+ T^{(1-\theta)/8}|I|^{\beta+(1-\theta)/4 }) \|u\|_{X^{0}_{4,\Delta}(I)}. \label{lebesgue4}
\end{align}
In particular, if $1 \lesssim  |I| <\infty$, $0<T<1$, then for any $0< \varepsilon \ll 1,$
\begin{align}
 \|u_I\|_{L^4_{x,t\in [0,T]}} \lesssim   T^{\varepsilon/3}|I|^{\varepsilon}  \|u\|_{X^{0}_{4,\Delta} (I)}. \label{lebesgue4a}
\end{align}
\end{lem}

\begin{lem}\label{L4toX}
Let $I\subset \mathbb{R}$ with $1\lesssim  |I|<\infty$,   $0<T<1$, then for any $4\leq q <\infty$,  $0< \varepsilon \ll 1,$ we have
\begin{align}
 \|u_I\|_{L^4_{x,t\in [0,T]}} & \lesssim   T^{\varepsilon/3}|I|^{1/4-1/q + \varepsilon} \max_{\lambda \in I} \langle \lambda\rangle^{-1/2} \|u\|_{X^{1/2}_{q,\Delta}(I)},  \label{lebesgue4b}\\
\|u_I\|_{L^\infty_t L^2_{x} \cap V^2_\Delta} & \lesssim    |I|^{1/2-1/q } \max_{\lambda \in I} \langle \lambda\rangle^{-1/2} \|u\|_{X^{1/2}_{q,\Delta}(I)},  \label{lebesgue4c}\\
\|u_I\|_{L^\infty_{x,t}} & \lesssim    |I|^{1-1/q } \max_{\lambda \in I} \langle \lambda\rangle^{-1/2} \|u\|_{X^{1/2}_{q,\Delta}(I)}.   \label{lebesgue4d}
\end{align}
\end{lem}
{\bf Proof.} By \eqref{lebesgue4a} and H\"older's inequality, we have \eqref{lebesgue4b}. Using $M^0_{2,2}=L^2$,   $V^2_\Delta \subset L^\infty_t L^2_x$  and  H\"older's inequality, we have \eqref{lebesgue4c}. In view of Bernstein's inequality,
$$
\|u_I\|_{L^\infty_{x,t}}   \lesssim    |I|^{1/2 } \|u_I\|_{L^\infty_{t}L^2_x}.
$$
Combining it with  \eqref{lebesgue4c}, we have \eqref{lebesgue4d}. $\hfill\Box$

\begin{cor}\label{L4toXj}
Let  $I_j =[a2^j, b2^j]$, $0<a<b$,    $0<T<1$, then for any $4\leq q <\infty$,  $0< \varepsilon \ll 1,$ we have
\begin{align}
 \|u_{I_j}\|_{L^4_{x,t\in [0,T]}} & \lesssim   T^{\varepsilon/3}2^{j(-1/4-1/q + \varepsilon)}   \|u\|_{X^{1/2}_{q,\Delta}(I_j)},  \label{lebesgue4bb}\\
\|u_{I_j}\|_{L^\infty_t L^2_{x} \cap V^2_\Delta} & \lesssim    2^{-j/q}     \|u\|_{X^{1/2}_{q,\Delta}(I_j)},  \label{lebesgue4cc}\\
\|u_{I_j}\|_{L^\infty_{x,t}} & \lesssim    2^{j(1/2-1/q) }      \|u\|_{X^{1/2}_{q,\Delta}(I_j)}.   \label{lebesgue4dd}
\end{align}
\end{cor}

\begin{cor}\label{L4toXjlambda}
Let $\lambda_0\gg 1$, $\lambda\in \mathbb{N}$,  $I_j =\lambda_0 -[a2^j, b2^j] \subset [c\lambda_0, \lambda_0]$ for some $0<a<b$,    $0<T<1$, then for any $4\leq q <\infty$,  $0< \varepsilon \ll 1,$ we have
\begin{align}
 \|u_{I_j}\|_{L^4_{x,t\in [0,T]}} & \lesssim   T^{\varepsilon/3} \langle\lambda_0\rangle^{-1/2} 2^{j(1/4-1/q + \varepsilon)}   \|u\|_{X^{1/2}_{q,\Delta}(I_j)},  \label{lebesgue4bbb}\\
\|u_{I_j}\|_{L^\infty_t L^2_{x}} & \lesssim  \langle\lambda_0\rangle^{-1/2}  2^{j(1/2-1/q)}     \|u\|_{X^{1/2}_{q,\Delta}(I_j)},  \label{lebesgue4ccc}\\
\|u_{I_j}\|_{L^\infty_{x,t}} & \lesssim    \langle\lambda_0\rangle^{-1/2}  2^{j(1-1/q)}     \|u\|_{X^{1/2}_{q,\Delta}(I_j)}.   \label{lebesgue4ddd}
\end{align}
\end{cor}

\begin{cor}\label{L4toXlambda}
Let $\lambda_0\in  \mathbb{R}$ with $\lambda_0\gg 1$, $I\subset [c\lambda_0, \lambda_0]$ or $I\subset [-\lambda_0, -c\lambda_0]$ for some $c\in (0,1)$,   $0<T<1$, then for any $4\leq q <\infty$,  $0< \varepsilon \ll 1,$ we have
\begin{align}
 \|u_I\|_{L^4_{x,t\in [0,T]}} & \lesssim   T^{\varepsilon/3} \langle \lambda_0\rangle^{-1/4-1/q + \varepsilon}  \|u\|_{X^{1/2}_{q,\Delta}(I)},  \label{lebesgue4bbbb}\\
\|u_I\|_{L^\infty_t L^2_{x} \cap V^2_\Delta} & \lesssim     \langle \lambda_0\rangle^{-1/q }  \|u\|_{X^{1/2}_{q,\Delta}(I)},  \label{lebesgue4cccc}\\
\|u_I\|_{L^\infty_{x,t}} & \lesssim    \langle \lambda_0\rangle^{1/2-1/q }  \|u\|_{X^{1/2}_{q,\Delta}(I)}.   \label{lebesgue4dddd}
\end{align}
\end{cor}

\section{Bilinear Estimates}

The following bilinear estimate can be found in \cite{Gr05}, \cite{Gu16}.

\begin{lem} [\rm Bilinear Estimate 1]  \label{Bilinear}
Let $0<T\leq 1$. Suppose that $\widehat{u}, \ \widehat{v}$ are localized in some compact intervals $I_1,I_2$ with $dist(I_1, I_2)\geq \lambda$. Then for any $0<\varepsilon \ll 1$, we have
\begin{align}
 \|u\overline{v} \|_{L^2_{x, t\in[0,T]}} \lesssim T^{\varepsilon/4} \lambda^{-1/2 + \varepsilon} \|u\|_{V^2_\Delta} \|v\|_{V^2_\Delta}. \label{bilinear}
\end{align}
\end{lem}
Similar to Gr\"unrock's bilinear estimate, we can consider the following bilinear estimate, which is useful for our late purpose.
\begin{lem}[\rm Bilinear Estimate 2] \label{V2decay2}
Let $0<T\leq 1$. Suppose that $\widehat{u}, \widehat{v}$ are localized in some compact intervals $I_1,I_2$ with $dist(I_1, I_2)\geq \lambda$. Then for any $0<\varepsilon \ll 1$, we have
\begin{align}
 \|u v \|_{L^2_{x, t\in[0,T]}} \lesssim T^{\varepsilon/4} \lambda^{-1/2 + \varepsilon} \|u\|_{V^2_\Delta} \|v\|_{V^2_\Delta}. \label{v2decay}
\end{align}
\end{lem}
{\bf Proof.} First, we show that if $\widehat{u}_0, \widehat{v}_0$ are localized in $I_1, \ I_2$ with $dist (I_1, I_2)\geq \lambda$, then
\begin{align}
\left\| e^{it\Delta}u_0 e^{it\Delta}v_0 \right\|_{L^2_{x,t\in \mathbb{R}}}  \lesssim \lambda^{-1/2}  \|u_0\|_2 \|v_0\|_2. \label{bilinearest}
\end{align}
We have
\begin{align}
& \mathscr{F}_x  \left( e^{it\Delta}u_0  e^{it\Delta}v_0  \right)  = \int e^{-it(\xi^2-2\xi\xi_1 +2\xi^2_1)} \widehat{u}_0(\xi-\xi_1) \widehat{v}_0(\xi_1) d\xi_1.
 \label{bilinearest1}
\end{align}
It follows that
\begin{align}
& \mathscr{F}_{x,t}  \left( e^{it\Delta}u_0 e^{it\Delta}v_0  \right)   = \int  \delta(2\xi^2_1-2\xi\xi_1+ \xi^2+\tau)  \widehat{u}_0(\xi-\xi_1) \widehat{v}_0(\xi_1) d\xi_1.
 \label{bilinearest2}
\end{align}
Denote
$$
g(\xi_1) = 2\xi^2_1-2\xi\xi_1+ \xi^2+\tau,
$$
we see that
$$
g'(\xi_1) = 4\xi_1 -2\xi, \ \  g(\xi^\pm_1)=0, \ \ \xi^\pm = \frac{\xi}{2} \pm \sqrt{\frac{\xi^2}{4}- \frac{\xi^2+\tau}{2}}:=\frac{\xi}{2} \pm y .
$$
Recall that $\delta(g(\xi_1)) = \delta (\xi_1-\xi^+)/ |g'(\xi^+)| - \delta (\xi_1-\xi^-)/ |g'(\xi^-)|= \delta (\xi_1-\xi^+)/ 4y -\delta (\xi_1-\xi^-)/4y $, we have
\begin{align}
& \mathscr{F}_{x,t}  \left( e^{it\Delta}u_0 e^{it\Delta}v_0  \right)  = \frac{1}{4y} \widehat{u}_0 \left(\frac{\xi}{2} - y\right) \widehat{v}_0\left(\frac{\xi}{2} +y\right) - \frac{1}{4y} \widehat{u}_0\left(\frac{\xi}{2} + y\right) \widehat{v}_0 \left(\frac{\xi}{2} -y\right). \label{bilinear4.6}
\end{align}
By symmetry, it suffices to estimate the first term in \eqref{bilinear4.6}.  Changing of variables $y= \sqrt{\frac{\xi^2}{4}- \frac{\xi^2+\tau}{2}}$, we see that
\begin{align}
\left\| e^{it\Delta}u_0 e^{it\Delta}v_0  \right\|^2_{L^2_{x,t}} &  \leq \int_{\mathbb{R}^2}\frac{c}{|y|}\left |\widehat{u}_0 \left(\frac{\xi}{2} - y \right)\right|^2 \left | \widehat{v}_0\left(\frac{\xi}{2} +y\right)\right|^2 dy d\xi \  \nonumber\\
&  \lesssim \int_{\mathbb{R}^2}\frac{1}{|\xi_1-\xi_2|} |\widehat{u}_0(\xi_1)|^2 | \widehat{v}_0(\xi_2)|^2 d\xi_1d\xi_2  \notag\\
&  \lesssim \lambda^{-1} \int_{\mathbb{R}^2} |\widehat{u}_0(\xi_1)|^2 | \widehat{v}_0(\xi_2)|^2 d\xi_1d\xi_2=C \lambda^{-1}\|u_0\|^2_2\|v_0\|^2_2,
\label{bilinearestmm}
\end{align}
where in the last inequality, we have applied $dist(I_1, I_2) >\lambda$. By testing atoms in $U^2$, then applying the interpolation in Proposition \ref{UVprop2}, we have the Bilinear Estimate 2. $\hfill\Box$

\section{Trilinear estimates}

We need to have a bound of the second term of the integral equation \eqref{IDNLS} in $X^{1/2}_{p,\Delta}$. More precisely, we want to show that
\begin{align}
& \left\|\int_{0}^{t} e^{(t-s)\Delta} (u^2\partial_{x}\bar{u})(s) \, ds \right\|_{X^{1/2}_{p,\Delta}}\lesssim T^{\varepsilon}\|u\|^{3}_{X^{1/2}_{p,\Delta}},  \label{trilinear1}\\
& \left \|\int_{0}^{t} e^{(t-s)\Delta} (uv \partial_x \bar{w})(s) \, ds \right\|_{X^{1/2}_{p,\Delta}}\lesssim T^{\varepsilon} \|u\|_{X^{1/2}_{p,\Delta}}\|v\|_{X^{1/2}_{p,\Delta}}\|w\|_{X^{1/2}_{p,\Delta}}.  \label{trilinear2}
\end{align}
{\bf Proof of \eqref{trilinear1}.} In view of \eqref{normest},  it suffices to show that
\begin{align}
&  \left|\int_{\mathbb{R}\times [0,T]} v u^2 \partial_x \overline{u} dxdt \right|\lesssim T^{\varepsilon}\|u\|^{3}_{X^{\frac{1}{2}}_{2,p}}  \|v\|_{Y^{-1/2}_{p'}} . \label{trilinear3}
\end{align}
We perform a uniform decomposition with $u,v$ in the left hand side of \eqref{trilinear3}, it suffice to prove that
\begin{align}
\left\lvert \sum_{\lambda_0,...,\lambda_3} \langle \lambda_0 \rangle^{1/2}\int_{[0,T]\times\mathbb{R}} \overline{v}_{\lambda_0} u_{\lambda_1}u_{\lambda_2}\partial_{x}\overline{u}_{\lambda_3} \, dxdt \right \rvert \lesssim T^{\varepsilon} \| u\| ^{3}_{X^{1/2}_{p,\Delta}}  \|v\|_{Y^{0}_{p',\Delta}}. \label{trilinear4}
\end{align}
In order to keep the left hand side of \eqref{trilinear4} nonzero, we have the frequency constraint condition (FCC)
\begin{equation}\label{E:frequencyrelation}
\lambda_1+\lambda_2-\lambda_3-\lambda_0  \thickapprox  0
\end{equation}
and dispersion modulation constraint condition (DMCC)
\begin{equation}\label{E:modulationrelation}
\max_{0\leq k \leq 3} |\xi^2_k + \tau_k|  \gtrsim |(\xi_0- \xi_1)(\xi_0-\xi_2)|,
\end{equation}
where we assume that $\widehat{v}_{\lambda_0}, \widehat{u}_{\lambda_1},..., \widehat{u}_{\lambda_3}$ in \eqref{trilinear4} are functions of $(\xi_0, \tau_0),..., (\xi_3, \tau_3)$, respectively. It suffices to consider the cases that $\lambda_0$ is minimal or secondly minimal number in $\lambda_0,...,\lambda_3$ (In the opposite case, one can instead $\lambda_0,...,\lambda_3$ by $-\lambda_0,...,-\lambda_3$).

{\bf Step 1.} We assume that $\lambda_0= \min_{0\leq k \leq 3} \lambda_k$. We separate the proof into two subcase $\lambda_0 \geq 0$ and $\lambda_0<0$.

{\bf Step 1.1.} We consider the case $\lambda_0 \geq 0$.  Let us denote $I_0=[0,1)$,  $I_{j}=[ 2^{j-1}, 2^{j}), \ j\geq 1.$ We decompose $\lambda_0 +[0,\infty)$  by dyadic decomposition, i.e.,
$$
\lambda_k \in \lambda_0 +[0,\infty) =\bigcup_{j_k\geq 0}(\lambda_0 +I_{j_k}), \ \ k=1,2,3.
$$
Recall that
 \[
u_{\lambda_0+I_{j_k}}=\sum_{\lambda_k\in \lambda_0+I_{j_k}} u_{\lambda_k}, \, \ \ k=1,2,3.
\]
By the symmetry we can assume that $j_1 \leq j_2$. Moreover, in view of  FCC \eqref{E:frequencyrelation}, we see that  $j_2 \thickapprox j_3$.
It follows that
\begin{align}
  \mathscr{L}^+ (u,v) &  :=  \sum_{0\leq \lambda_0 \leq \lambda_k, \ k=1,2,3 } \left\langle\lambda_0\right\rangle^{1/2}\int_{\left[0,T\right] \times \mathbb{R}} \overline{v}_{\lambda_0} {u_{\lambda_1}} {u_{\lambda_2}} \partial_x \overline{u}_{\lambda_3} \,dxdt   \notag\\
&   = \sum_{\lambda_0\geq 0, \ j_1\leq j_2, \ j_3\thickapprox j_2 } \left\langle\lambda_0\right\rangle^{1/2}\int_{\left[0,T\right] \times \mathbb{R}} \overline{v}_{\lambda_0} {u_{\lambda_0+I_{j_1}}} {u_{\lambda_0+ I_{j_2}}} \partial_x \overline{u}_{\lambda_0+I_{j_3}} \,dxdt   \notag\\
&  \leq  \left (\sum_{\lambda_0\geq 0, \ j_1\leq j_2 \thickapprox j_3 \lesssim 1 } + \sum_{\lambda_0\geq 0, \ j_1 \leq 10, \  j_2 \thickapprox j_3 \gg  1 }  + \sum_{\lambda_0\geq 0, \ 10< j_1 \leq  j_2 \thickapprox j_3 } \right)  \langle\lambda_0 \rangle^{1/2} \nonumber\\
& \ \ \ \ \times \int_{\left[0,T\right] \times \mathbb{R}} |\overline{v}_{\lambda_0} {u_{\lambda_0+I_{j_1}}} {u_{\lambda_0+ I_{j_2}}} \partial_x \overline{u}_{\lambda_0+I_{j_3}}| \,dxdt   \notag\\
& : = \mathscr{L}^+_{l} (u,v) + \mathscr{L}^+_{m} (u,v) + \mathscr{L}^+_{h} (u,v). \label{decompL}
 \end{align}
It is easy to see that in $\mathscr{L}^+_{l}(u,v)$, the frequency of $u_{\lambda_0+I_{j_k}}$ ($k=1,2,3$) are all localized in a neighbourhood of $\lambda_0$. So, by H\"older's inequality,  $\|\Box_k u\|_4 \leq \|\Box_k u\|_2$ and $U^2_\Delta \subset V^2_\Delta \subset L^\infty_tL^2_x$,  we have
 \begin{align}
|\mathscr{L}^+_{l}(u,v)|
 & \lesssim \sum_{\lambda_0 \thickapprox \lambda_1 \thickapprox \lambda_2  \thickapprox \lambda_3} \left\langle\lambda_0\right\rangle^{3/2} \|v_{\lambda_0}\|_{L^{4}_{t}L^{4}_{x}} \|u_{\lambda_1}\|_{L^{4}_{t}L^{4}_{x}} \|u_{\lambda_2}\|_{L^{4}_{t}L^{4}_{x}} \|u_{\lambda_3}\|_{L^{4}_{t}L^{4}_{x}} \nonumber \\
 & \lesssim T  \|u\|^{3}_{X^{1/2}_{p, \Delta}}\|v\|_{Y^0_{p', \Delta}}.    \label{trilinear5}
\end{align}
In $\mathscr{L}^+_{m} (u,v)$,  we easily see that the frequency of $v_{\lambda_0}$ and $u_{\lambda_0+I_{j_1}}$ are localized near $\lambda_0$, which are much less than those of $u_{\lambda_0+I_{j_2}}$ and $u_{\lambda_0+I_{j_3}}$. So, we can use bilinear estimate \eqref{bilinear} and Lemma \ref{L4toX} to obtain that
\begin{align}
|\mathscr{L}^+_{m} (u,v)|
& \lesssim  \sum_{\lambda_0\geq 0, \ j_1 \leq 10, \  j_2 \thickapprox j_3 \gg  1 }  \langle\lambda_0 \rangle^{1/2}
\|{u_{\lambda_0+I_{j_1}}} {\partial_x \overline{u}_{\lambda_0+I_{j_3}}} \|_{L^2_{x,t}}  \|{\overline{v}_{\lambda_0}} {u_{\lambda_0+I_{j_2}}}\|_{L^2_{x,t}} \notag\\
& \lesssim T^{\varepsilon/2} \sum_{\lambda_0\geq 0, \ j_1 \leq 10, \  j_2 \thickapprox j_3 \gg  1 }  \langle\lambda_0 \rangle^{1/2}\langle \lambda_0+ 2^{j_3}\rangle  2^{(-1/2+\varepsilon)j_3 }
\|{u_{\lambda_0+I_{j_1}}} \|_{V^2_\Delta}\|{u_{\lambda_0+I_{j_3}}}\|_{V^2_\Delta} \nonumber\\
 & \ \ \ \ \ \times 2^{(-1/2+\varepsilon)j_2} \|{u_{\lambda_0+I_{j_2}}}\|_{V^2_\Delta}\|{v_{\lambda_0}}\|_{V^2_\Delta} \notag \\
 & \lesssim T^{\varepsilon/2} \!\!\!\!\!\!\!\! \sum_{\lambda_0\geq 0,    j \gg  1, |\ell_1|\vee |\ell_2| \lesssim 1 }  \!\!\!\!\!\!\!\! \langle\lambda_0 \rangle^{1/2}  \langle \lambda_0+ 2^{j}\rangle   2^{ (-1+2 \varepsilon) j  }\|{v_{\lambda_0}}\|_{V^2_\Delta}
\| u_{\lambda_0+ \ell_1 }  \|_{V^2_\Delta}\|{u_{\lambda_0+I_{j+\ell_2}}}\|_{V^2_\Delta}
  \| u_{\lambda_0+I_{j}} \|_{V^2_\Delta}\notag  \\
& \lesssim T^{\varepsilon/2} \!\!\!\!\!\!\!\! \sum_{\lambda_0\geq 0,  \  j \gg  1, |\ell |  \lesssim 1 } \!\!\!\!\!\!\!\!  2^{(-2/p + 2 \varepsilon) j  }\|{v_{\lambda_0}}\|_{V^2_\Delta}
\langle\lambda_0 \rangle^{1/2}  \| u_{\lambda_0+ \ell  }  \|_{V^2_\Delta}
  \| u  \|^2_{X^{1/2}_{p,\Delta}}. \label{trilinear9}
 \end{align}
Noticing that $0<\varepsilon <1/p$,  we have from H\"older's inequality that
\begin{align}
|\mathscr{L}^+_{m} (u,v)|  &  \lesssim T^{\varepsilon/2} \| u  \|^2_{X^{1/2}_{p,\Delta}}  \sum_{\lambda_0\geq 0,    \ |\ell |  \lesssim 1 }     \|{v_{\lambda_0}}\|_{V^2_\Delta}
\langle\lambda_0 \rangle^{1/2}  \| u_{\lambda_0+ \ell  }  \|_{V^2_\Delta}    \notag \\
&  \lesssim T^{\varepsilon/2} \| u  \|^3_{X^{1/2}_{p,\Delta}} \|v\|_{X^0_{p', \Delta}}. \label{trilinear10}
\end{align}
Now we estimate $\mathscr{L}^+_{h} (u,v)$.

{\it Case 1.}   $v_{\lambda_0}$ has the highest dispersion modulation in the right hand side of $\mathscr{L}^+_{h} (u,v)$ in \eqref{decompL}. In view of DMCC \eqref{E:modulationrelation}, we have
$$
|\xi^2_0 + \tau_0| \gtrsim 2^{j_1+ j_2}.
$$
It follows from the dispersion modulation decay \eqref{dispersiondecay} that
\begin{align} \label{highestmod}
\|v_{\lambda_0}\|_{L^2_{x,t}} \lesssim 2^{-(j_1+ j_2)/2} \|v_{\lambda_0}\|_{V^2_\Delta}.
\end{align}
By H\"older's inequality, we have
$$
|\mathscr{L}^+_{h} (u,v)| \leq  \!\!\! \!\!\! \sum_{\lambda_0\geq 0, \ 10< j_1 \leq  j_2 \thickapprox j_3 } \langle \lambda_0\rangle^{1/2}  \|v_{\lambda_0}\|_{L^2_tL^\infty_x}  \|u_{\lambda_0+I_{j_1}} \|_{L^\infty_t L^2_{x}}  \|u_{\lambda_0+I_{j_2}} \|_{L^4_{x,t\in [0,T]}} \|u_{\lambda_0+I_{j_3}} \|_{L^4_{x,t\in [0,T]}}.
$$
Using $\|v_\lambda\|_\infty \lesssim  \|v_\lambda\|_2$, the bound of the highest modulation in \eqref{highestmod} and Lemma \ref{L4toX}, we have
\begin{align}
|\mathscr{L}^+_{h} (u,v)| & \lesssim  T^{\varepsilon/2} \sum_{\lambda_0\geq 0, \ 10< j_1 \leq  j_2 \thickapprox j_3 } \langle \lambda_0\rangle^{1/2}  2^{-(j_1+ j_2)/2}  \|v_{\lambda_0}\|_{V^2_\Delta} \langle \lambda_0 +2^{j_3}\rangle    \nonumber\\
& \ \ \ \ \times  2^{ j_1(1/2-1/p)+ j_2(1/4-1/p + \varepsilon)+ j_3(1/4-1/p+ \varepsilon)}  \prod^3_{k=1} \langle \lambda_0 +2^{j_k}\rangle^{-1/2}  \|u  \|_{X^{1/2}_{p,\Delta}(\lambda_0+I_{j_k})}. \nonumber
\end{align}
Making the summation on $j_2, j_3$, applying H\"older's inequality on $\lambda_0$ and finally summing over all $j_1$,  we have for $0< \varepsilon <1/2p$,
\begin{align}
|\mathscr{L}^+_{h} (u,v)| & \lesssim  T^{\varepsilon/2} \sum_{\lambda_0\geq 0,  \  j_1>10 }   2^{ j_1(-3/p +2 \varepsilon) }  \|v_{\lambda_0}\|_{V^2_\Delta}     \|u  \|_{X^{1/2}_{p,\Delta}(\lambda_0+I_{j_1})}     \|u  \|^2_{X^{1/2}_{p,\Delta}}. \nonumber\\
&  \lesssim T^{\varepsilon/2} \| u  \|^3_{X^{1/2}_{p,\Delta}} \|v\|_{X^0_{p', \Delta}}. \label{trilinear11}
\end{align}

{\it Case 2.}   $u_{\lambda_0+I_{j_1}}$ has the highest dispersion modulation in the right hand side of $\mathscr{L}^+_{h} (u,v)$ in \eqref{decompL}.  In order to use DMCC \eqref{E:modulationrelation}, $u_{\lambda_0+I_{j_1}}$  takes $L^2_{x,t}$ norm,  $u_{\lambda_0+I_{j_2}}$ and $u_{\lambda_0+I_{j_3}}$ take $L^4_{x,t}$ norms. Indeed, we have
$$
|\mathscr{L}^+_{h} (u,v)| \leq  \!\!\! \!\!\! \sum_{\lambda_0\geq 0, \ 10< j_1 \leq  j_2 \thickapprox j_3 } \langle \lambda_0\rangle^{1/2}  \|v_{\lambda_0}\|_{L^\infty_{x,t}}  \|u_{\lambda_0+I_{j_1}} \|_{L^2_{x,t}}  \|u_{\lambda_0+I_{j_2}} \|_{L^4_{x,t\in [0,T]}} \|u_{\lambda_0+I_{j_3}} \|_{L^4_{x,t\in [0,T]}}.
$$
Using the fact $\|v_{\lambda_0}\|_{L^\infty_{x,t}}  \lesssim  \|v_{\lambda_0}\|_{L^\infty_{ t} L^2_x}  \lesssim \|v_{\lambda_0}\|_{V^2_\Delta}$,  the dispersion modulation decay estimate \eqref{highestmod} and Lemma \ref{L4toX}, this case reduces to the same estimate as in Case 1, see \eqref{trilinear11}.

{\it Case 3.}   $u_{\lambda_0+I_{j_2}}$ has the highest dispersion modulation in the right hand side of $\mathscr{L}^+_{h} (u,v)$ in \eqref{decompL}. If $j_1 \thickapprox j_3$ in the right hand side of $\mathscr{L}^+_{h} (u,v)$ in \eqref{decompL},  then we can repeat the proof  as in Case 2 to obtain the desired estimates. So, it suffuces to consider the case $j_1\ll j_3$.  $j_1\ll j_3$ implies that the frequency of $u_{\lambda_0+I_{j_3}}$ is much higher than that of $u_{\lambda_0+I_{j_1}}$ and we can use the bilinear estimate \eqref{v2decay}.
By H\"older's inequality,
$$
|\mathscr{L}^+_{h} (u,v)| \leq  \!\!\! \!\!\! \sum_{\lambda_0\geq 0, \ 10< j_1 \ll  j_2 \thickapprox j_3 } \langle \lambda_0\rangle^{1/2}  \|v_{\lambda_0}\|_{L^\infty_{x,t}}  \|u_{\lambda_0+I_{j_1}}\partial_x \overline{u}_{\lambda_0+I_{j_3}} \|_{L^2_{x,t}}  \|u_{\lambda_0+I_{j_2}} \|_{L^2_{x,t\in [0,T]}}.
$$
Applying the bilinear estimate \eqref{bilinear}, DMCC \eqref{E:modulationrelation} and Lemma \ref{L4toX}, we have
\begin{align}
|\mathscr{L}^+_{h} (u,v)| & \leq  \!\!\! \!\!\! \sum_{\lambda_0\geq 0, \ 10< j_1 \ll  j_2 \thickapprox j_3 } \langle \lambda_0\rangle^{1/2}  \|v_{\lambda_0}\|_{V^2_\Delta} T^{\varepsilon/4} 2^{j_3(-1/2+\varepsilon)} \|u_{\lambda_0+I_{j_1}}\|_{V^2_\Delta} \| u_{\lambda_0+I_{j_3}} \|_{V^2_\Delta} \nonumber\\
&  \ \ \ \ \ \times \langle\lambda_0 + 2^{j_3}\rangle  2^{-(j_1+j_2)/2}\|u_{\lambda_0+I_{j_2}} \|_{V^2_\Delta} \nonumber\\
& \leq  T^{\varepsilon/4}  \!\!\!   \sum_{\lambda_0\geq 0, \ 10< j_1 \ll  j_2 \thickapprox j_3 }  2^{-(j_1+j_2 +j_3)/p + \varepsilon j_3} \|v_{\lambda_0}\|_{V^2_\Delta} \|u \|_{X^{1/2}_{p,\Delta}(\lambda_0+I_{j_1}) } \| u  \|^2_{X^{1/2 }_{p,\Delta}} \nonumber
  \end{align}
Making the summation on $j_2, j_3$, we have
\begin{align}
|\mathscr{L}^+_{h} (u,v)|
  \leq  T^{\varepsilon/4}  \!\!\!   \sum_{\lambda_0\geq 0, \   j_1 >10}  2^{j_1(-3 /p + \varepsilon)} \|v_{\lambda_0}\|_{V^2_\Delta} \|u \|_{X^{1/2}_{p,\Delta}(\lambda_0+I_{j_1}) } \| u  \|^2_{X^{1/2 }_{p,\Delta}}. \nonumber
  \end{align}
This reduces to the same estimate as in the first inequality of \eqref{trilinear11}.

{\it Case 4.}   $u_{\lambda_0+I_{j_3}}$ has the highest modulation in the right hand side of $\mathscr{L}^+_{h} (u,v)$ in \eqref{decompL}. Noticing that $j_2 \thickapprox j_3$ in $\mathscr{L}^+_{h} (u,v)$, using the bilinear estimate\eqref{v2decay} instead of \eqref{bilinear} we see that this case is similar to Case 3.

{\bf Step 1.2.}  We consider the case $\lambda_0 <0$.  We can assume that $\lambda_0 \ll 0$.  For short, we denote by $\lambda  \in h_-$ and $\lambda\in l_-$ the fact $\lambda \in [\lambda_0, 3\lambda_0/4)$ and  $\lambda \in [3\lambda_0/4, 0)$, respectively. We need to consider the following three cases $\lambda_1 \in h_-$, $\lambda_1\in l_-$ and $\lambda_1\in [0, \infty)$ separately.

{\it Case A.} $\lambda_1\in h_-$.  We consider the following four subcases as in Table 1.
\begin{table}[h]
\begin{center}

\begin{tabular}{|c|c|c|c|}
\hline
 ${\rm Case}$   &  $\lambda_1\in $ & $\lambda_2\in $ & $\lambda_3\in $    \\
\hline
$h_-h_-a$  & $ [\lambda_0, 3\lambda_0/4)$  & $ [\lambda_0, 3\lambda_0/4)$  &  $[\lambda_0, \infty)$  \\
\hline
$h_-l_-l_-$ & $ [\lambda_0, 3\lambda_0/4)$  & $ [3\lambda_0/4, \ 0)$ & $ [\lambda_0, 0)$  \\
\hline
$h_-l_- a_+$ & $ [\lambda_0, 3\lambda_0/4)$  & $ [3\lambda_0/4,0)$  & $ [0,  \infty)$\\
\hline
$h_-a_+ a$ & $ [\lambda_0, 3\lambda_0/4)$  & $ [0, \infty)$  & $ [\lambda_0,  \infty)$\\
\hline
\end{tabular}
\end{center}
\caption{In $\lambda_1,...,\lambda_3$, there is at least one frequency near $\lambda_0$}
\end{table}

{\it Case $h_-h_-a$.}   By FCC \eqref{E:frequencyrelation}, we see that $\lambda_3 \in [\lambda_0, \lambda_0/2+ C)$. Hence, $\lambda_3$ is also near $\lambda_0$.  Using the same way as in the proof of Step 1, we can get the result and the details are omitted.

{\it Case $h_-l_-l_-$.} By FCC \eqref{E:frequencyrelation}, we see that $\lambda_3 \geq 3\lambda_0/4 -C.$ We decompose $\lambda_k$ in the following way:
$$
\lambda_1\in \lambda_0 +[0, - \frac{\lambda_0}{4} ) = \bigcup_{j_1\geq 0} (\lambda_0 + I_{j_1});\footnote{In order to realize the identity, we can assume that $I_{j_1} = [2^{j-1}, 2^j) \cap [0, -\lambda_0/4)$, it is convenient to use such a kind of notations below. } \ \lambda_2 \in  [ \frac{3\lambda_0}{4}, 0) = \bigcup_{j_2\geq 0} - I_{j_2}; \ \  \lambda_3 \in  [\frac{3\lambda_0}{4} -C, 0) = \bigcup_{j_3\geq 0} - I_{j_3}.
$$
Again, in view of \eqref{E:frequencyrelation} we see that
$$
(j_1\vee j_3 ) \   \thickapprox j_2.
$$
We have the following three subcases:
$$
j_1\ll j_3,  \ j_2 \thickapprox j_3; \ \ or \ \  j_3 \ll j_1, \  j_1 \thickapprox j_2; \ \ or \ \ j_1 \thickapprox j_2 \thickapprox j_3
$$
It follows that
\begin{align}
  & \sum_{  \lambda_0 \ll 0 } \left\langle\lambda_0\right\rangle^{1/2}\int_{\left[0,T\right] \times \mathbb{R}} |\overline{v}_{\lambda_0} {u_{[\lambda_0, \ 3\lambda_0/4)}} {u_{[3\lambda_0/4, \ 0)}} \partial_x \overline{u}_{[3\lambda_0/4-C, \ 0)}| \,dxdt   \notag\\
&  \leq  \sum_{\lambda_0\ll 0, \ (j_1\vee  j_3) \thickapprox j_2 \lesssim \ln \langle \lambda_0\rangle } \left\langle\lambda_0\right\rangle^{1/2}\int_{\left[0,T\right] \times \mathbb{R}} |\overline{v}_{\lambda_0} {u_{\lambda_0+I_{j_1}}} {u_{- I_{j_2}}} \partial_x \overline{u}_{-I_{j_3}}| \,dxdt   \notag
 \end{align}
We consider the case $j_1\ll j_3,  \ j_2 \thickapprox j_3.$  We have
\begin{align}
  \mathscr{L}^-_{h_-l_-l_-} (u,v) & :=   \sum_{\lambda_0\ll 0, \ j_1\ll  j_2 \thickapprox j_3 \lesssim \ln \langle \lambda_0\rangle } \left\langle\lambda_0\right\rangle^{1/2}\int_{\left[0,T\right] \times \mathbb{R}} |\overline{v}_{\lambda_0} {u_{\lambda_0+I_{j_1}}} {u_{- I_{j_2}}} \partial_x \overline{u}_{-I_{j_3}}| \,dxdt   \notag\\
&  \leq \left ( \sum_{\lambda_0\ll 0, \ j_1 \leq 10 \ll  j_2 \thickapprox j_3 \lesssim \ln \langle \lambda_0\rangle }  + \sum_{\lambda_0\ll 0, \ 10< j_1 \ll  j_2 \thickapprox j_3 \lesssim \ln \langle \lambda_0\rangle } \right)  \langle\lambda_0 \rangle^{1/2} \nonumber\\
& \ \ \ \ \times \int_{\left[0,T\right] \times \mathbb{R}}| \overline{v}_{\lambda_0} {u_{\lambda_0+I_{j_1}}} {u_{-I_{j_2}}} \partial_x \overline{u}_{-I_{j_3}} | \,dxdt   \notag\\
& : =  \mathscr{L}^{-,m}_{h_-l_-l_-} (u,v) + \mathscr{L}^{-,h}_{h_-l_-l_-} (u,v). \label{decompL<0}
 \end{align}
In order to estimate $\mathscr{L}^{-,m}_{h_-l_-l_-} (u,v)$, we follow the same ideas as in \eqref{trilinear9}. We may assume that $ 2^{10}  \leq -\lambda_0/16$. It follows that $\lambda_0 + 2^{j_1} + 2^{j_2} \leq \lambda_0/16$. Now we can use the bilinear estimate \eqref{bilinear}, Corollary \ref{L4toXj}.  For $0<\varepsilon <1/2p$, we have
\begin{align}
   \mathscr{L}^{-,m}_{h_-l_-l_-} (u,v)  &
    \lesssim    \sum_{\lambda_0\ll 0, \ j_1 \leq 10 \ll j_2 \thickapprox j_3 \lesssim \ln \langle \lambda_0\rangle }    \langle\lambda_0 \rangle^{1/2}  \|\overline{v}_{\lambda_0}\partial_x \overline{u}_{-I_{j_3}} \|_{L^2_{x,t\in [0,T]}}  \| {u_{\lambda_0+I_{j_1}}} {u_{-I_{j_2}}} \|_{L^2_{x,t\in [0,T]}}   \notag\\
&  \lesssim  T^{\varepsilon/2} \!\!\!\!\!\! \sum_{\lambda_0\ll 0, \ |\ell|, |\ell_1|  \lesssim 1, \     1 \ll j  \lesssim \ln \langle \lambda_0\rangle } \!\!\!\!\!\!    \langle\lambda_0 \rangle^{-1/2+ 2\varepsilon} 2^{j}  \|{v}_{\lambda_0}\|_{V^2_\Delta} \|{u}_{-I_{j}} \|_{V^2_\Delta}  \| {u_{\lambda_0+ \ell}}\|_{V^2_\Delta} \|{u_{-I_{j+\ell_1}}} \|_{V^2_\Delta}   \notag\\
&  \lesssim  T^{\varepsilon/2} \!\!\!\!\!\! \sum_{\lambda_0\ll 0, \ |\ell| \lesssim 1, \     1 \ll j  \lesssim \ln \langle \lambda_0\rangle }   \!\!\!\!\!\!  \langle\lambda_0 \rangle^{-1/2+ 2\varepsilon}  2^{j-2j/p} \|{v}_{\lambda_0}\|_{V^2_\Delta}  \| {u_{\lambda_0+ \ell}}\|_{V^2_\Delta}\|{u}  \|^2_{X^{1/2}_{p,\Delta}} \notag\\
& \lesssim    T^{\varepsilon/2}   \|v\|_{Y^0_{p',\Delta}} \|u  \|^3_{X^{1/2}_{p,\Delta}} .    \label{}
 \end{align}
For the estimate of $\mathscr{L}^{-,h}_{h_-l_-l_-} (u,v)$, we need to us DMCC \eqref{E:modulationrelation}, we have
\begin{equation}\label{highestmod2}
\max_{0\leq k \leq 3} |\xi^2_k +\tau_k|  \gtrsim |(\xi_0- \xi_1)(\xi_0-\xi_2)| \gtrsim 2^{j_1} \langle \lambda_0\rangle.
\end{equation}
If $v_{\lambda_0}$ has the highest dispersion modulation,  we have
\begin{align}
\mathscr{L}^{-,h}_{h_-l_-l_-} (u,v)
& \lesssim \!\!\!\!\!\!\!\!\!\!\!\!  \sum_{\lambda_0\ll 0, \ 10< j_1 \ll  j_2 \thickapprox j_3 \lesssim \ln \langle \lambda_0\rangle  } \!\!\!\!\!\!\!\!\!\!\!\!   \langle \lambda_0 \rangle ^{1/2} 2^{j_3}  \|v_{\lambda_0}\|_{L^{2}_{t}L^{\infty}_{x}}\|u_{\lambda_0+I_{j_1}}\|_{L^{\infty}_{t}L^{2}_{x}}\|u_{-I_{j_2}}\|_{L^{4}_{x,t}}\|u_{-I_{j_3}}\|_{L^{4}_{x, t}}. \label{trilinear16}
\end{align}
In view of \eqref{dispersiondecay}, Corrollary \ref{L4toXj}, we have
\begin{align}
\mathscr{L}^{-,h}_{h_-l_-l_-} (u,v)
& \lesssim  T^{2\varepsilon/3} \sum_{\lambda_0\ll 0, \ 10< j_1 \ll  j_2 \thickapprox j_3 \lesssim \ln \langle \lambda_0\rangle }  \langle \lambda_0 \rangle ^{1/2} 2^{j_3} 2^{- j_1/2 }\langle\lambda_0\rangle^{-1/2} \|v_{\lambda_0}\|_{V^{2}_{\Delta}} \notag \\
& \  \ \times \langle\lambda_0\rangle^{-1/2} 2^{j_1(1/2-1/p)} \|u\|_{X^{1/2}_{p,\Delta}(\lambda_0+I_{j_1})}  2^{(j_2+j_3)(-1/4-1/p +\varepsilon)} \|u\|^2_{X^{1/2}_{p,\Delta}}  \notag\\
& \lesssim  T^{2\varepsilon/3}  \sum_{\lambda_0\ll 0, \ 10< j_1 \ll  j_3 \lesssim \ln \langle \lambda_0\rangle } 2^{-j_1/p}  2^{j_3/2 + (-2/p+ 2\varepsilon)j_3} \langle\lambda_0\rangle^{-1/2}  \|v_{\lambda_0}\|_{V^{2}_{\Delta}}  \notag\\
  & \ \ \times \|u\|_{X^{1/2}_{p,\Delta}(\lambda_0+I_{j_1})} \|u\|^2_{X^{1/2}_{p,\Delta}}. \label{trilinear17}
\end{align}
Making the summation on $j_3$, then using the same way as in \eqref{trilinear11}, one has that for $0< \varepsilon < 1/2p$,
\begin{align}
\mathscr{L}^{-,h}_{h_-l_-l_-} (u,v)
& \lesssim  T^{2\varepsilon/3}  \sum_{\lambda_0\ll 0, \ 10< j_1   }    2^{ (-3/p+ 2\varepsilon)j_1}  \|v_{\lambda_0}\|_{V^{2}_{\Delta}}
   \|u\|_{X^{1/2}_{p,\Delta}(\lambda_0+I_{j_1})} \|u\|^2_{X^{1/2}_{p,\Delta}}. \notag \\
  & \lesssim T^{ \varepsilon/2}\|u\|^{3}_{X^{1/2}_{p, \Delta}}\|v\|_{Y_{p', \Delta}}. \label{trilinear18}
\end{align}
If $u_{\lambda_0+I_{j_1}}$ has the highest dispersion modulation,  we have
\begin{align}
\mathscr{L}^{-,h}_{h_-l_-l_-} (u,v)
& \lesssim \!\!\!\!\!\!  \sum_{\lambda_0\ll 0, \ 10<j_1 \ll j_2 \thickapprox j_3 \lesssim \ln \langle \lambda_0\rangle  } \!\!\!\!\!\!   \langle \lambda_0 \rangle ^{1/2} 2^{j_3}  \|v_{\lambda_0}\|_{ L^{\infty}_{x,t}}\|u_{\lambda_0+I_{j_1}}\|_{ L^{2}_{x,t}}\|u_{-I_{j_2}}\|_{L^{4}_{x,t}}\|u_{-I_{j_3}}\|_{L^{4}_{x, t}}.   \label{trilinear19}
\end{align}
Applying the dispersion modulation decay estimate \eqref{dispersiondecay} to $u_{\lambda_0+I_{j_1}}$, and  $\|v_{\lambda_0}\|_\infty \lesssim  \|v_{\lambda_0}\|_2$, we can reduce the estimate of   \eqref{trilinear19} to the case as in \eqref{trilinear17} and \eqref{trilinear18}, the details are omitted.

Now we consider the case that $u_{-I_{j_2}}$ has the highest modulation.  In $\mathscr{L}^{-,h}_{h_-l_-l_-} (u,v)$, it is easy to see that $j_1\ll \ln \langle \lambda_0\rangle$, which implies that $\lambda_0 + 2^{j_1} + 2^{j_3} \leq \lambda_0 /16.$  Using the bilinear estimate \eqref{bilinear} to $\partial_x \overline{u}_{-I_{j_3}} u_{\lambda_0+I_{j_1}}$, and the dispersion modulation decay estimate \eqref{dispersiondecay}, one has that
\begin{align}
\mathscr{L}^{-,h}_{h_-l_-l_-} (u,v)
& \lesssim \!\!\!\!\!\!  \sum_{\lambda_0\ll 0, \ 10< j_1 \ll j_2 \thickapprox j_3 \lesssim \ln \langle \lambda_0\rangle  } \!\!\!\!\!\!   \langle \lambda_0 \rangle ^{1/2}    \|v_{\lambda_0}\|_{ L^{\infty}_{x,t}}\|u_{\lambda_0+I_{j_1}} \partial_x \overline{u}_{-I_{j_3}} \|_{ L^{2}_{x,t}}\|u_{-I_{j_2}}\|_{L^{2}_{x,t}} \nonumber\\
& \lesssim \!\!\!\!\!\!  \sum_{\lambda_0\ll 0, \ 10< j_1 \ll j_2 \thickapprox j_3 \lesssim \ln \langle \lambda_0\rangle  } \!\!\!\!\!\!   \langle \lambda_0 \rangle ^{1/2}    \|v_{\lambda_0}\|_{V^2_\Delta} T^{\varepsilon/4} \langle \lambda_0\rangle^{-1/2+ \varepsilon}\|u_{\lambda_0+I_{j_1}}\|_{V^2_\Delta} \| \partial_x u_{-I_{j_3}} \|_{V^2_\Delta}  \nonumber\\
& \ \ \ \ \times 2^{-j_1/2} \langle \lambda_0\rangle^{-1/2} \|u_{-I_{j_2}}\|_{_{V^2_\Delta}}
   \label{trilinear20}
\end{align}
It follows from Lemma \ref{L4toX} and \eqref{trilinear20} that
\begin{align}
\mathscr{L}^{-,h}_{h_-l_-l_-} (u,v)
  \lesssim T^{\varepsilon/4} \!\!\!\!  \sum_{\lambda_0\ll 0, \ j_1 \leq j_3 \lesssim \ln \langle \lambda_0\rangle }  \!\!\!\!\!   \|v_{\lambda_0}\|_{V^2_\Delta}  \langle \lambda_0\rangle^{-1+ \varepsilon}  2^{-j_1/p}  2^{j_3(1-2/p)}  \|u \|_{X^{1/2}_{p, \Delta}(\lambda_0+I_{j_1})}
   \|u\|^2_{X^{1/2}_{p, \Delta} }.
   \label{trilinear21}
\end{align}
Making the summation on $j_3$ and using the same way as in \eqref{trilinear11}, we have for $0<\varepsilon <1/2p$,
\begin{align}
\mathscr{L}^{-,h}_{h_-l_-l_-} (u,v)
  & \lesssim T^{\varepsilon/4} \!\!  \sum_{\lambda_0\ll 0, \ j_1  }    \|v_{\lambda_0}\|_{V^2_\Delta}     2^{-2 j_1/p}    \|u \|_{X^{1/2}_{p, \Delta}(\lambda_0+I_{j_1})} \|u\|^2_{X^{1/2}_{p, \Delta} } \nonumber\\
 & \lesssim T^{ \varepsilon/2}\|u\|^{3}_{X^{1/2}_{p, \Delta}}\|v\|_{Y_{p', \Delta}}.    \label{trilinear22}
\end{align}

We consider the case $j_3 \ll j_1 \thickapprox j_2.$  We denote
\begin{align}
  \mathscr{H}^-_{h_-l_-l_-} (u,v) &  :=    \sum_{\lambda_0\ll 0, \ j_3\ll  j_1 \thickapprox j_2 \lesssim \ln \langle \lambda_0\rangle } \left\langle\lambda_0\right\rangle^{1/2}\int_{\left[0,T\right] \times \mathbb{R}} |\overline{v}_{\lambda_0} {u_{\lambda_0+I_{j_1}}} {u_{- I_{j_2}}} \partial_x \overline{u}_{-I_{j_3}}| \,dxdt
  \label{decompL<0a}
 \end{align}
In $\mathscr{H}^{-,h}_{h_-l_-l_-} (u,v)$, observing DMCC \eqref{E:modulationrelation}, we have
\begin{equation}\label{highestmod3}
\max_{0\leq k \leq 3} |\xi^2_k +\tau_k|  \gtrsim |(\xi_0- \xi_1)(\xi_0-\xi_2)| \gtrsim 2^{j_1} \langle \lambda_0\rangle.
\end{equation}
In view of \eqref{highestmod3} we see that the lower bound of the highest dispersion modulation is the same as that of the case $j_1 \ll j_3 \thickapprox j_2$. Moreover,  $\partial_x u_{-I_{j_3}}$ has the lower frequency, which leads to that the derivative in front of $u_{-I_{j_3}}$ becomes easier to handle.  So, we omit the details of the proof in this case.

For the case $j_1 \thickapprox j_2  \thickapprox j_3$, we divide the proof into two subcases: $ j_3 \lesssim 1$ and $j_3 \gg 1.$ The first case is very easy, since both $u_{-I_{j_2}}$ and $u_{-I_{j_3}}$ have the low frequency in a neighbourhood of $0$. One can directly use H\"older's inequality to get the desired estimate. For the case $j_3\gg 1$,  we consider that the highest dispersion modulation (larger than $   2^{j_1} \langle \lambda_0\rangle$) is due to $v_{\lambda_0}, $ $u_{\lambda_0 +I_{j_1}}$, $u_{-I_{j_2}}$ and $u_{-I_{j_3}}$, separately.  Denote
\begin{align*}
\mathscr{J}^{-,h}_{h_-l_-l_-} (u,v)
& \lesssim \!\!\!\!\!\!  \sum_{\lambda_0\ll 0, \ 1 \ll j_1 \thickapprox  j_2 \thickapprox j_3 \lesssim \ln \langle \lambda_0\rangle  } \!\!\!\!\!\!   \langle \lambda_0 \rangle ^{1/2}  \int_{\left[0,T\right] \times \mathbb{R}} |\overline{v}_{\lambda_0} {u_{\lambda_0+I_{j_1}}} {u_{- I_{j_2}}} \partial_x \overline{u}_{-I_{j_3}}| \,dxdt.
\end{align*}
If $v_{\lambda_0}$ has the highest dispersion modulation,  we have
\begin{align*}
\mathscr{J}^{-,h}_{h_-l_-l_-} (u,v)
& \lesssim \!\!\!\!\!\!  \sum_{\lambda_0\ll 0, \ 1 \ll j_1 \thickapprox  j_2 \thickapprox j_3 \lesssim \ln \langle \lambda_0\rangle  } \!\!\!\!\!\!   \langle \lambda_0 \rangle ^{1/2} 2^{j_3}  \|v_{\lambda_0}\|_{L^{2}_{t}L^{\infty}_{x}}\|u_{\lambda_0+I_{j_1}}\|_{L^{\infty}_{t}L^{2}_{x}}\|u_{-I_{j_2}}\|_{L^{4}_{x,t}}\|u_{-I_{j_3}}\|_{L^{4}_{x, t}}. \notag
\end{align*}
We can follow the same way as in dealing with \eqref{trilinear16} to have the desired estimate. If $u_{\lambda_0 +I_{j_1}}$ has the highest dispersion modulation, we have
\begin{align*}
\mathscr{J}^{-,h}_{h_-l_-l_-} (u,v)
& \lesssim \!\!\!\!\!\!  \sum_{\lambda_0\ll 0, \ 1 \ll j_1 \thickapprox  j_2 \thickapprox j_3 \lesssim \ln \langle \lambda_0\rangle  } \!\!\!\!\!\!   \langle \lambda_0 \rangle ^{1/2} 2^{j_3}  \|v_{\lambda_0}\|_{ L^{\infty}_{x,t}}\|u_{\lambda_0+I_{j_1}}\|_{L^{2}_{x,t}}\|u_{-I_{j_2}}\|_{L^{4}_{x,t}}\|u_{-I_{j_3}}\|_{L^{4}_{x, t}}. \notag
\end{align*}
If $u_{-I_{j_2}}$ has the highest dispersion modulation, we have
\begin{align*}
\mathscr{J}^{-,h}_{h_-l_-l_-} (u,v)
& \lesssim \!\!\!\!\!\!  \sum_{\lambda_0\ll 0, \ 1 \ll j_1 \thickapprox  j_2 \thickapprox j_3 \lesssim \ln \langle \lambda_0\rangle  } \!\!\!\!\!\!   \langle \lambda_0 \rangle ^{1/2} 2^{j_3}  \|v_{\lambda_0}\|_{ L^{\infty}_{x,t}}\|u_{\lambda_0+I_{j_1}}\|_{L^{4}_{x,t}}\|u_{-I_{j_2}}\|_{L^{2}_{x,t}}\|u_{-I_{j_3}}\|_{L^{4}_{x, t}}. \notag
\end{align*}
Then, using \eqref{dispersiondecay}, Lemma \ref{L4toX} and similar estimates as in the above, we can get the the result, as desired.

{\it Case $h_-l_-a_+$.} By FCC \eqref{E:frequencyrelation}, we see that $\lambda_3 <  -\lambda_0/4 + C.$ We decompose $\lambda_k$:
$$
\lambda_1\in \lambda_0 +[0, - \frac{\lambda_0}{4} ) = \bigcup_{j_1\geq 0} (\lambda_0 + I_{j_1}); \ \lambda_2 \in  [ \frac{3\lambda_0}{4}, 0) = \bigcup_{j_2\geq 0} - I_{j_2}; \ \  \lambda_3 \in  [0, \ -\frac{ \lambda_0}{4}+C) = \bigcup_{j_3\geq 0}   I_{j_3}.
$$
Again, in view of  FCC \eqref{E:frequencyrelation} we see that
$$
(j_2\vee j_3 ) \   \thickapprox j_1, \ \ j_1 \leq \ln \langle \lambda_0\rangle +C.
$$
It follows that
$$
j_2 \ll j_3 \thickapprox j_1; \ \ or \ \  j_3 \ll j_2 \thickapprox j_1; \ \ or \ \ j_1 \thickapprox j_2 \thickapprox j_3.
$$
We have
\begin{align}
  & \sum_{  \lambda_0 \ll 0 } \left\langle\lambda_0\right\rangle^{1/2}\int_{\left[0,T\right] \times \mathbb{R}} |\overline{v}_{\lambda_0} {u_{[\lambda_0, \ 3\lambda_0/4)}} {u_{[3\lambda_0/4, \ 0)}} \partial_x \overline{u}_{[0, -\lambda_0/4+C)}| \,dxdt   \notag\\
&  \leq  \sum_{\lambda_0\ll 0, \ (j_2\vee  j_3) \thickapprox j_1 \lesssim \ln \langle \lambda_0\rangle } \left\langle\lambda_0\right\rangle^{1/2}\int_{\left[0,T\right] \times \mathbb{R}} |\overline{v}_{\lambda_0} {u_{\lambda_0+I_{j_1}}} {u_{- I_{j_2}}} \partial_x \overline{u}_{I_{j_3}}| \,dxdt   \notag
 \end{align}
From the proof of the Case $h_-l_-l_-$, it suffices to consider the case $j_2\ll j_3 \thickapprox j_1.$  We estimate
\begin{align}
  \mathscr{L}^-_{h_-l_-l} (u,v) & :=   \sum_{\lambda_0\ll 0, \ j_2\ll  j_3 \thickapprox j_1 \lesssim \ln \langle \lambda_0\rangle } \left\langle\lambda_0\right\rangle^{1/2}\int_{\left[0,T\right] \times \mathbb{R}} |\overline{v}_{\lambda_0} {u_{\lambda_0+I_{j_1}}} {u_{- I_{j_2}}} \partial_x \overline{u}_{I_{j_3}}| \,dxdt. \label{decompL<0hll}
 \end{align}
In the right hand side of $ \mathscr{L}^-_{h_-l_-l} (u,v)$, by DMCC \eqref{E:modulationrelation},
\begin{equation}\label{highestmod5}
\max_{0\leq k \leq 3} |\xi^2_k +\tau_k|   \gtrsim 2^{j_1} \langle \lambda_0\rangle.
\end{equation}
If $v_{\lambda_0}$ has the highest dispersion modulation, in a similar way as in \eqref{trilinear16}--\eqref{trilinear18},  we have for $0<\varepsilon<1/2p$,
\begin{align}
\mathscr{L}^{-}_{h_-l_-l} (u,v)
& \lesssim \!\!\!\!\!\!  \sum_{\lambda_0\ll 0, \  j_2 \ll  j_3 \thickapprox j_1 \lesssim \ln \langle \lambda_0\rangle  } \!\!\!\!\!\!   \langle \lambda_0 \rangle ^{1/2} 2^{j_3}  \|v_{\lambda_0}\|_{L^{2}_{t}L^{\infty}_{x}}\|u_{\lambda_0+I_{j_1}}\|_{L^{4}_{x,t}}\|u_{-I_{j_2}}\|_{L^\infty_t L^{2}_{x}}\|u_{I_{j_3}}\|_{L^{4}_{x, t}} \nonumber\\
& \lesssim  T^{2\varepsilon/3}  \sum_{\lambda_0\ll 0, \    j_1 \lesssim \ln \langle \lambda_0\rangle }   2^{(-2/p+ 2\varepsilon)j_1}   \|v_{\lambda_0}\|_{V^{2}_{\Delta}}  \|u\|_{X^{1/2}_{p,\Delta}(\lambda_0+I_{j_1})} \|u\|^2_{X^{1/2}_{p,\Delta}} \notag\\
  & \lesssim T^{ \varepsilon/2}\|u\|^{3}_{X^{1/2}_{p, \Delta}}\|v\|_{Y_{p', \Delta}}. \label{trilinear18+}
\end{align}
Now we consider the case that $u_{-I_{j_2}}$ has the highest dispersion modulation.
\begin{align}
\mathscr{L}^{-}_{h_-l_-l} (u,v)
& \lesssim \!\!\!\!\!\!  \sum_{\lambda_0\ll 0, \ j_2 \ll j_3 \thickapprox j_1 \lesssim \ln \langle \lambda_0\rangle  } \!\!\!\!\!\!   \langle \lambda_0 \rangle ^{1/2} 2^{j_3}  \|v_{\lambda_0}\|_{ L^{\infty}_{x,t}}\|u_{\lambda_0+I_{j_1}}\|_{ L^{4}_{x,t}}\|u_{-I_{j_2}}\|_{L^{2}_{x,t}}\|u_{I_{j_3}}\|_{L^{4}_{x, t}}.   \label{trilinear19+}
\end{align}
Applying the dispersion modulation decay estimate \eqref{dispersiondecay} to $u_{-I_{j_2}}$,   we can reduce the estimate of   \eqref{trilinear19+} to the case as in  \eqref{trilinear18+}, the details are omitted.

Assume that $u_{\lambda_0+I_{j_1}}$ has the highest dispersion modulation.
   Using the bilinear estimate \eqref{bilinear} to $u_{-I_{j_2}} \partial_x \overline{u}_{I_{j_3}}$, and the dispersion modulation decay estimate \eqref{dispersiondecay}, then applying Lemma \ref{L4toX},  one has that for $0<\varepsilon <1/2p$,
\begin{align}
\mathscr{L}^{-}_{h_-l_-l} (u,v)
& \lesssim \!\!\!\!\!\!  \sum_{\lambda_0\ll 0, \ j_2 \ll j_3 \thickapprox j_1 \lesssim \ln \langle \lambda_0\rangle  } \!\!\!\!\!\!   \langle \lambda_0 \rangle ^{1/2}    \|v_{\lambda_0}\|_{ L^{\infty}_{x,t}}\|u_{\lambda_0+I_{j_1}} \|_{ L^{2}_{x,t}}\|u_{-I_{j_2}} \partial_x \overline{u}_{I_{j_3}}\|_{L^{2}_{x,t}} \nonumber\\
& \lesssim T^{\varepsilon/4} \!\!\!\!  \sum_{\lambda_0\ll 0, \ j_1   \lesssim \ln \langle \lambda_0\rangle }  \!\!\!\!\!   \|v_{\lambda_0}\|_{V^2_\Delta}  2^{j_1(-2/p+ \varepsilon)}  \|u \|_{X^{1/2}_{p, \Delta}(\lambda_0+I_{j_1})}
   \|u\|^2_{X^{1/2}_{p, \Delta} } \nonumber\\
 & \lesssim T^{ \varepsilon/4}\|u\|^{3}_{X^{1/2}_{p, \Delta}}\|v\|_{Y_{p', \Delta}}.    \label{trilinear22+}
\end{align}

If $u_{I_{j_3}}$ has the highest dispersion modulation, we can use the same way as in the case that $u_{\lambda_0+I_{j_1}}$ has the highest dispersion modulation to obtain the result. The details are omitted.

{\it Case $h_-a_+a$. } First, using FCC \eqref{E:frequencyrelation}, we have $\lambda_3\geq -C$.  Moreover, if $\lambda_3\leq C$, we see that $\lambda_2 \leq C$. We have
\begin{align}
  & \sum_{  \lambda_0 \ll 0 } \left\langle\lambda_0\right\rangle^{1/2}\left|\int_{\left[0,T\right] \times \mathbb{R}}\overline{v}_{\lambda_0} {u_{[\lambda_0, \ 3\lambda_0/4)}} {u_{[0,\ \infty)}} \partial_x \overline{u}_{[ \lambda_0, \infty)} \,dxdt \right|   \notag\\
 & \leq   \sum_{  \lambda_0 \ll 0 } \left\langle\lambda_0\right\rangle^{1/2}\left| \int_{\left[0,T\right] \times \mathbb{R}}\overline{v}_{\lambda_0} {u_{[\lambda_0, \ 3\lambda_0/4)}} {u_{[0,\ C)}} \partial_x \overline{u}_{[-C, C)} \,dxdt \right| \notag\\
& \quad \ \ +  \sum_{  \lambda_0 \ll 0 } \left\langle\lambda_0\right\rangle^{1/2} \left | \int_{\left[0,T\right] \times \mathbb{R}}\overline{v}_{\lambda_0} {u_{[\lambda_0, \ 3\lambda_0/4)}} {u_{[0,\ \infty)}} \partial_x \overline{u}_{[C, \ \infty)} \,dxdt \right|  \notag\\
& := L^{-}_{h_-a_+l} (u,v) + L^{-}_{h_-a_+h} (u,v).
 \end{align}
In  $L^{-}_{h_-a_+l} (u,v)$, there are two lower frequency in a neighbourhood of the origin and  two higher frequency near $\lambda_0$.  Hence, we can use the bilinear estimate \eqref{bilinear} to obtain that
\begin{align}
\mathscr{L}^{-}_{h_-a_+l} (u,v)
 & \lesssim T^{ \varepsilon/2}\|u\|^{3}_{X^{1/2}_{p, \Delta}}\|v\|_{Y_{p', \Delta}}.    \label{trilinear31}
\end{align}
We estimate  $L^{-}_{h_-a_+h} (u,v)$. Considering the dyadic version of $\lambda_k$:
$$
\lambda_1\in \lambda_0 +[0, - \frac{\lambda_0}{4} ) = \bigcup_{j_1\geq 0} (\lambda_0 + I_{j_1}); \ \lambda_2 \in  [0,\infty) = \bigcup_{j_2\geq 0}  I_{j_2}; \ \  \lambda_3 \in  [C, \ \infty) = \bigcup_{j_3 \gtrsim 1}   I_{j_3}.
$$
Again, in view of  FCC \eqref{E:frequencyrelation} we see that
$$
(j_1 \vee j_2 ) \   \thickapprox j_3, \ \ j_1 \leq \ln \langle \lambda_0\rangle +C.
$$
  We have the following three subcases:
$$
j_1 \ll j_2  \thickapprox j_3 ; \ \ or \ \  j_2 \ll j_1 \thickapprox j_3  ; \ \ or \ \ j_1 \thickapprox j_2 \thickapprox j_3.
$$
Using the dyadic version, we have
\begin{align}
  \mathscr{L}^-_{h_-a_+h} (u,v) & :=   \left(\sum_{\lambda_0\ll 0, \ j_1\ll  j_2 \thickapprox j_3   }  + \sum_{\lambda_0\ll 0, \ j_2\ll  j_1 \thickapprox j_3 } + \sum_{\lambda_0\ll 0, \ j_1 \thickapprox j_2 \thickapprox j_3  } \right)   \langle\lambda_0 \rangle^{1/2} \nonumber \\
& \ \ \ \ \ \times \int_{\left[0,T\right] \times \mathbb{R}} |\overline{v}_{\lambda_0} u_{\lambda_0+I_{j_1}} u_{I_{j_2}}  \partial_x \overline{u}_{I_{j_3}}| \,dxdt := I+II+III. \label{trilinear32}
 \end{align}
We mainly estimate $I$. Using the bilinear estimate \eqref{v2decay}, Corollary \ref{L4toXj}, H\"older's inequality,  we have
\begin{align}
I_l   & :=    \sum_{\lambda_0\ll 0, \ j_1 \leq 10 \ll  j_2 \thickapprox j_3  }    \langle\lambda_0 \rangle^{1/2}
\int_{\left[0,T\right] \times \mathbb{R}} |\overline{v}_{\lambda_0} u_{\lambda_0+I_{j_1}} u_{I_{j_2}}  \partial_x \overline{u}_{I_{j_3}}| \,dxdt  \nonumber\\
&  \leq    \sum_{\lambda_0\ll 0, \ j_1 \leq 10 \ll  j_2 \thickapprox j_3  }    \langle\lambda_0 \rangle^{1/2}
 \|\overline{v}_{\lambda_0}\partial_x \overline{u}_{I_{j_3}}\|_{L^2_{x,t}}  \|u_{\lambda_0+I_{j_1}} u_{I_{j_2}}\|_{L^2_{x,t}}   \nonumber\\
&  \leq  T^{\varepsilon/2}  \sum_{\lambda_0\ll 0, \ j_1 \leq 10 \ll  j_2 \thickapprox j_3  }    \langle\lambda_0 \rangle^{1/2} (\langle\lambda_0 \rangle + 2^{j_2})^{-1+2\varepsilon} 2^{j_2}
 \| {v}_{\lambda_0}\|_{V^2_\Delta}\| {u}_{I_{j_3}}\|_{V^2_\Delta}  \|u_{\lambda_0+I_{j_1}}\|_{V^2_\Delta}\| u_{I_{j_2}}\|_{V^2_\Delta}   \nonumber\\
&  \leq  T^{\varepsilon/2}  \sum_{\lambda_0\ll 0,  |\ell|\lesssim 1,    j_2 \gg 1  }    \langle\lambda_0 \rangle^{1/2}  2^{j_2(-2/p+2\varepsilon)}
 \| {v}_{\lambda_0}\|_{V^2_\Delta}  \|u_{\lambda_0+\ell}\|_{V^2_\Delta} \| u \|^2_{X^{1/2}_{p,\Delta}}   \nonumber\\
 & \lesssim T^{ \varepsilon/2}\|u\|^{3}_{X^{1/2}_{p, \Delta}}\|v\|_{Y_{p', \Delta}}.  \label{trilinear33}
 \end{align}
We estimate
\begin{align}
I_h   & :=    \sum_{\lambda_0\ll 0, \ 10 \leq j_1   \ll  j_2 \thickapprox j_3 }    \langle\lambda_0 \rangle^{1/2}
\int_{\left[0,T\right] \times \mathbb{R}} |\overline{v}_{\lambda_0} u_{\lambda_0+I_{j_1}} u_{I_{j_2}}  \partial_x \overline{u}_{I_{j_3}}| \,dxdt. \label{trilinear34}
 \end{align}
In the right hand side of $I_h$, in view of DMCC \eqref{E:modulationrelation}, the highest dispersion modulation
\begin{equation}\label{highestmod35}
\max_{0\leq k \leq 3} |\xi^2_k +\tau_k|   \gtrsim 2^{j_1} (\langle \lambda_0\rangle +2^{j_2}).
\end{equation}
If $v_{\lambda_0}$ has the highest dispersion modulation, in a similar way as in \eqref{trilinear16}--\eqref{trilinear18},  we have
\begin{align}
I_h
& \lesssim  \!\!\!\!  \sum_{\lambda_0\ll 0, \ 10<  j_1 \ll  j_2 \thickapprox j_3  }  \!\!\!\!   \langle \lambda_0 \rangle ^{1/2} 2^{j_3}  \|v_{\lambda_0}\|_{L^{2}_{t}L^{\infty}_{x}}\|u_{\lambda_0+I_{j_1}}\|_{L^{\infty}_{t}L^{2}_{x}}\|u_{I_{j_2}}\|_{L^{4}_{x,t}}\|u_{I_{j_3}}\|_{L^{4}_{x, t}} \nonumber\\
& \lesssim  \!\!\!\!  \sum_{\lambda_0\ll 0, \ 10<  j_1 \ll  j_2 \thickapprox j_3  }  \!\!\!\!   \langle \lambda_0 \rangle ^{1/2} 2^{j_3} 2^{-j_1/2} (\langle \lambda_0 \rangle+ 2^{j_2}) ^{-1/2} \|v_{\lambda_0}\|_{V^2_\Delta}\|u_{\lambda_0+I_{j_1}}\|_{V^2_\Delta}\|u_{I_{j_2}}\|_{L^{4}_{x,t}}\|u_{I_{j_3}}\|_{L^{4}_{x, t}} \nonumber\\
& \lesssim T^{2\varepsilon/3} \!\!\!\! \sum_{\lambda_0\ll 0, \ 10<  j_1 \ll  j_2  } \!\!\!\!    \langle \lambda_0 \rangle ^{1/2} 2^{-j_1/2} (\langle \lambda_0\rangle +2^{j_2})^{-1/2} \|v_{\lambda_0}\|_{V^2_\Delta} \notag\\
& \ \ \ \ \ \ \ \times \langle \lambda_0 \rangle ^{-1/2} 2^{j_1(1/2-1/p)} \|u \|_{X^{1/2}_{p,\Delta}(\lambda_0+I_{j_1})} 2^{j_2(-1/2-2/p+2\varepsilon)}2^{j_2}\|u \|^2_{X^{1/2}_{p,\Delta}}.  \label{trilinear36}
\end{align}
Summing  over all $j_2 $, then using H\"older's on $\lambda_0$,  we obtain that for $0< \varepsilon<1/2 p$,
\begin{align}
I_h
  & \lesssim T^{ \varepsilon/2}\|u\|^{3}_{X^{1/2}_{p, \Delta}}\|v\|_{Y_{p', \Delta}}. \label{trilinear37}
\end{align}
If $u_{\lambda_0+I_{j_1}}$ has the highest dispersion modulation, in a similar way as in \eqref{trilinear36} and \eqref{trilinear37},  we have for $0<\varepsilon<1/2p$,
\begin{align}
I_h
& \lesssim  \!\!\!\!  \sum_{\lambda_0\ll 0, \ 10<  j_1 \ll  j_2 \thickapprox j_3  }  \!\!\!\!   \langle \lambda_0 \rangle ^{1/2} 2^{j_3}  \|v_{\lambda_0}\|_{ L^{\infty}_{x,t}}\|u_{\lambda_0+I_{j_1}}\|_{L^{2}_{x, t}}\|u_{I_{j_2}}\|_{L^{4}_{x,t}}\|u_{I_{j_3}}\|_{L^{4}_{x, t}} \nonumber\\
  & \lesssim T^{ \varepsilon/2}\|u\|^{3}_{X^{1/2}_{p, \Delta}}\|v\|_{Y_{p', \Delta}}. \label{trilinear38}
\end{align}
If $u_{ I_{j_3}}$ has the highest dispersion modulation,    we have
\begin{align}
I_h
& \lesssim  \!\!\!\!  \sum_{\lambda_0\ll 0, \ 10<  j_1 \ll  j_2 \thickapprox j_3  }  \!\!\!\!   \langle \lambda_0 \rangle ^{1/2} 2^{j_3}  \|v_{\lambda_0}\|_{ L^{\infty}_{x,t}}\|u_{\lambda_0+I_{j_1}} u_{I_{j_2}}\|_{L^{2}_{x, t}} \|u_{I_{j_3}}\|_{L^{2}_{x, t}}.  \label{trilinear39}
\end{align}
Using the bilinear estimate \eqref{v2decay}, and the dispersion modulation decay \eqref{highestmod35} and \eqref{dispersiondecay}, then applying Lemma \ref{L4toX},
\begin{align}
I_h & \lesssim  T^{ \varepsilon/4}  \sum_{\lambda_0\ll 0, \    j_1 >10 }   2^{(-3/p+  \varepsilon)j_1}   \|v_{\lambda_0}\|_{V^{2}_{\Delta}}  \|u\|_{X^{1/2}_{p,\Delta}(\lambda_0+I_{j_1})} \|u\|^2_{X^{1/2}_{p,\Delta}} \notag\\
  & \lesssim T^{ \varepsilon/4}\|u\|^{3}_{X^{1/2}_{p, \Delta}}\|v\|_{Y_{p', \Delta}}. \label{trilinear40}
\end{align}
If $u_{ I_{j_2}}$ has the highest dispersion modulation,    we have
\begin{align}
I_h
& \lesssim  \!\!\!\!  \sum_{\lambda_0\ll 0, \ 10<  j_1 \ll  j_2 \thickapprox j_3  }  \!\!\!\!   \langle \lambda_0 \rangle ^{1/2} \|v_{\lambda_0}\|_{ L^{\infty}_{x,t}}\| u_{I_{j_2}}\|_{L^{2}_{x, t}} \|u_{\lambda_0+I_{j_1}} \partial_x \overline{u}_{I_{j_3}}\|_{L^{2}_{x, t}}.  \label{trilinear41}
\end{align}
One sees that  $dist (\lambda_0+I_{j_1}, \ I_{j_3} ) \gtrsim
\langle \lambda_0\rangle +
2^{j_3}.
$
Using the bilinear estimate \eqref{bilinear} and the dispersion modulation decay \eqref{dispersiondecay} and \eqref{highestmod35},
\begin{align}
I_h
& \lesssim T^{\varepsilon/4} \!\!\!\! \sum_{\lambda_0\ll 0, \ 10<  j_1 \ll  j_2 \thickapprox j_3 }  \!\!\!\!   \langle \lambda_0 \rangle ^{1/2} \|v_{\lambda_0}\|_{V^2_\Delta} 2^{-j_1/2} (\langle \lambda_0\rangle + 2^{j_2})^{-1+ \varepsilon}  \| u_{I_{j_2}}\|_{V^2_\Delta} \nonumber\\
 & \ \ \ \  \times 2^{j_3} \|u_{\lambda_0+I_{j_1}}\|_{V^2_\Delta} \| {u}_{I_{j_3}}\|_{V^2_\Delta}.   \label{trilinear42}
\end{align}
Applying Corollaries \ref{L4toXj}, \ref{L4toXjlambda} and making the summations in turn to $j_3$ and $j_2$, $\lambda_0,j_1$, we have
\begin{align}
I_h & \lesssim  T^{ \varepsilon/4}  \sum_{\lambda_0\ll 0, \   j_2\geq  j_1 >10 } 2^{j_2(-2/p+ \varepsilon)}  2^{-j_1/p}   \|v_{\lambda_0}\|_{V^{2}_{\Delta}}  \|u\|_{X^{1/2}_{p,\Delta}(\lambda_0+I_{j_1})} \|u\|^2_{X^{1/2}_{p,\Delta}} \notag\\
  & \lesssim T^{ \varepsilon/4}\|u\|^{3}_{X^{1/2}_{p, \Delta}}\|v\|_{Y_{p', \Delta}}. \label{trilinear43}
\end{align}

{\it Case B.} $\lambda_1\in l_-$. Since $\lambda_1$ and $\lambda_2$ are symmetry, it suffices to consider the following four subcases as in Table 2.
\begin{table}[h]
\begin{center}

\begin{tabular}{|c|c|c|c|}
\hline
 ${\rm Case}$   &  $\lambda_1\in $ & $\lambda_2\in $ & $\lambda_3\in $    \\
\hline
$l_-l_-l$  & $ [3\lambda_0/4, 0)$  & $ [3\lambda_0/4,0)$  &  $[\lambda_0, -\lambda_0/16)$  \\
\hline
$l_-l_-h_+$ & $ [ 3\lambda_0/4, 0)$  & $ [3\lambda_0/4, \ 0)$ & $ [-\lambda_0/16, \infty)$  \\
\hline
$l_-a_+a_m$ & $ [3\lambda_0/4, 0)$  & $ [0, \infty)$  & $ [\lambda_0,  -\lambda_0)$\\
\hline
$l_-a_+h_+$ & $ [3\lambda_0/4, 0)$  & $ [0, \infty)$  & $ [-\lambda_0,  \infty)$\\
\hline
\end{tabular}
\end{center}
\caption{ $\lambda_1,\lambda_2,\lambda_3$  far away from $\lambda_0$.}
\end{table}

{\it Case $l_-l_-l$.} From FCC \eqref{E:frequencyrelation} it follows that $\lambda_3 \geq \lambda_0/2-C$. Moreover, we see that $\lambda_1,  \lambda_2< \lambda_0/8$.  So, it suffices to estimate
\begin{align}
\mathscr{L}^{-}_{l_-l_-l} (u,v) :=   \sum_{  \lambda_0 \ll 0 } \left\langle\lambda_0\right\rangle^{1/2}\left|\int_{\left[0,T\right] \times \mathbb{R}}\overline{v}_{\lambda_0}  u_{[3\lambda_0/4,\lambda_0/8)}   u_{[3\lambda_0/4,\lambda_0/8)}   \partial_x \overline{u}_{[ \lambda_0/2-C, -\lambda_0/16)} \,dxdt \right|.   \label{trilinear44}
 \end{align}
Using the dyadic decomposition and Corollary \ref{L4toXj}, we see that
\begin{align}
&  \|\partial_x \overline{u}_{[ \lambda_0/2-C, 0)} \|_{L^\infty_tL^2_x \ \cap \ V^2_\Delta} \leqslant \sum_{j\leq \ln \langle \lambda_0\rangle} \|\partial_x u_{I_j}\|_{ V^2_\Delta} \leq \langle \lambda_0\rangle^{1-1/p} \|u\|_{X^{1/2}_{p,\Delta}};  \label{trilinear45}\\
&  \|\partial_x \overline{u}_{[ \lambda_0/2-C, 0)} \|_{L^4_{x,t}} \leqslant T^{\varepsilon/3} \sum_{j\leq \ln \langle \lambda_0\rangle } \|\partial_x u_{I_j}\|_{ L^4_{x,t} } \leq \langle \lambda_0\rangle^{3/4-1/p+ \varepsilon} \|u\|_{X^{1/2}_{p,\Delta}}.  \label{trilinear46}
 \end{align}
Similarly, $\partial_x \overline{u}_{[0, -\lambda_0/16)}$ has the same estimates as in \eqref{trilinear45} and \eqref{trilinear46}. By DMCC \eqref{E:modulationrelation}, the highest dispersion modulation satisfies
\begin{align}
\bigvee^3_{k=0} |\xi^2_k + \tau_k| \gtrsim  \langle  \lambda_0\rangle^2.  \label{trilinear46a}
\end{align}
If $v_{\lambda_0}$ has the highest dispersion modulation, we have
\begin{align}
\mathscr{L}^{-}_{l_-l_-l} (u,v)  \leq    \sum_{  \lambda_0 \ll 0 }  \langle\lambda_0 \rangle^{1/2} \|\overline{v}_{\lambda_0}\|_{L^2_tL^\infty_x}  \|u_{[3\lambda_0/4,\lambda_0/8)}\|^2_{L^4_{x,t}}  \| \partial_x \overline{u}_{[ \lambda_0/2-C, -\lambda_0/16)}\|_{L^\infty_tL^2_x}.       \label{trilinear47}
 \end{align}
Applying the highest dispersion modulation decay estimate \eqref{dispersiondecay}, Corollary \ref{L4toXlambda}, \eqref{trilinear45} and \eqref{trilinear46}, we have
\begin{align}
\mathscr{L}^{-}_{l_-l_-l} (u,v)  \leq   T^{\varepsilon/2}  \sum_{  \lambda_0 \ll 0 }    \langle\lambda_0 \rangle^{-3/p+2\varepsilon} \|v_{\lambda_0}\|_{V^2_\Delta}  \|u\|^3_{X^{1/2}_{p,\Delta}} \lesssim   T^{\varepsilon/2} \|v\|_{Y^{0}_{p',\Delta}} \|u\|^3_{X^{1/2}_{p,\Delta}}.      \label{trilinear48}
 \end{align}
If $u_{[ \lambda_0/2-C,-\lambda_0/16)}$ gains the highest dispersion modulation, we can use an analogous way to get the result. In fact, taking $L^\infty_{x,t}$,  $L^4_{x,t}$ and $L^2_{x,t}$ norms for $v_{\lambda_0}$, $u_{[3\lambda_0/4,\lambda_0/8)}$ and $\partial_x u_{[ \lambda_0/2-C,-\lambda_0/16)}$, respectively, then using  the dispersion modulation decay estimate \eqref{dispersiondecay}, Corollary \ref{L4toXlambda}, \eqref{trilinear45}, we also have \eqref{trilinear48}. Finally,  in the case $u_{[3\lambda_0/4,\lambda_0/8)}$ having the highest dispersion modulation, taking $L^\infty_{x,t}$,  $L^2_{x,t}$, $L^4_{x,t}$ and $L^4_{x,t}$ norms for $v_{\lambda_0}$, $u_{[3\lambda_0/4,\lambda_0/8)}$, $u_{[3\lambda_0/4,\lambda_0/8)}$  and $\partial_x u_{[ \lambda_0/2-C,-\lambda_0/16)}$, respectively, then using  the dispersion modulation decay estimate \eqref{dispersiondecay}, Corollary \ref{L4toXlambda}, \eqref{trilinear45} and \eqref{trilinear47}, we also have \eqref{trilinear48}.

{\it Case $l_-l_-h_+$. } We easily see that $\lambda_3\leq -\lambda_0+C$.  We collect $\lambda_k$ in the following dyadic version:
$$
\lambda_k \in [3\lambda_0/4, 0) = \bigcup_{j_k} -I_{j_k}, \ \ k=1,2; \ \ \lambda_3 \in [-\lambda_0/16, -\lambda_0 +C) = \bigcup_{ j_3\geq -1}(-\lambda_0  -I_{j_3}), \ \ j_{-1}=[0, C).
$$
By FCC \eqref{E:frequencyrelation},  we see that $j_3 \thickapprox (j_1\vee j_2) \lesssim \ln \langle \lambda_0\rangle$.  We need to estimate
\begin{align}
\mathscr{L}^{-}_{l_-l_-h_+} (u,v) &  :=   \sum_{  \lambda_0 \ll 0 } \left\langle\lambda_0\right\rangle^{1/2}\left|\int_{\left[0,T\right] \times \mathbb{R}}\overline{v}_{\lambda_0}  u_{[3\lambda_0/4,0)}   u_{[3\lambda_0/4,0)}   \partial_x \overline{u}_{[ -\lambda_0/16 , -\lambda_0 +C)} \,dxdt \right| \notag\\
 &  \lesssim   \left(\sum_{\lambda_0\ll 0, \ j_1\ll  j_2 \thickapprox j_3 }  + \sum_{\lambda_0\ll 0, \ j_2\ll  j_1 \thickapprox j_3  } + \sum_{\lambda_0\ll 0, \ j_1 \thickapprox j_2 \thickapprox j_3  } \right)   \langle\lambda_0 \rangle^{1/2} \nonumber \\
& \ \ \ \ \  \ \ \ \   \ \ \ \  \ \ \ \  \ \ \ \  \ \ \ \  \ \ \ \  \times \int_{\left[0,T\right] \times \mathbb{R}} |\overline{v}_{\lambda_0} u_{-I_{j_1}} u_{-I_{j_2}}  \partial_x \overline{u}_{-\lambda_0-I_{j_3}}| \,dxdt \notag\\
&  := \sum_{i=1,2,3} \mathscr{L}^{-,i}_{l_-l_-h_+} (u,v). \label{trilinear49}
 \end{align}
In the right hand side of \eqref{trilinear49}, from DMCC \eqref{E:modulationrelation} it follows that
\begin{equation}\label{highestmod50}
\max_{0\leq k \leq 3} |\xi^2_k +\tau_k|   \gtrsim   \langle \lambda_0\rangle^2  .
\end{equation}
If $v_{\lambda_0}$ attains the highest dispersion modulation, using a similar way as in \eqref{trilinear16}--\eqref{trilinear18},  we have for $0<\varepsilon<1/2p$,
\begin{align}
\mathscr{L}^{-,1}_{l_-l_-h_+} (u,v)
& \lesssim \!\!\!\!\!\!  \sum_{\lambda_0\ll 0, \  j_1 \ll  j_2 \thickapprox j_3 } \!\!\!\!\!\!   \langle \lambda_0 \rangle ^{1/2}   \|v_{\lambda_0}\|_{L^{2}_{t}L^{\infty}_{x}}\|u_{-I_{j_1}}\|_{L^{\infty}_{t}L^{2}_{x}}\|u_{-I_{j_2}}\|_{L^{4}_{x,t}}\|\partial_x u_{-\lambda_0-I_{j_3}}\|_{L^{4}_{x, t}} \notag\\
& \lesssim \!\!\!\!\!\!  \sum_{\lambda_0\ll 0, \  j_1 \ll  j_2 \thickapprox j_3 } \!\!\!\!\!\!   \langle \lambda_0 \rangle ^{-1/2}   \|v_{\lambda_0}\|_{V^2_\Delta}\|u_{-I_{j_1}}\|_{V^2_\Delta}\|u_{-I_{j_2}}\|_{L^{4}_{x,t}}\|\partial_x u_{-\lambda_0-I_{j_3}}\|_{L^{4}_{x, t}} \notag\\
& \lesssim  T^{2\varepsilon/3}  \sum_{  j_1 \ll  j_3, \ \lambda_0\ll 0} 2^{-j_1/p}  2^{(-2/p+ 2\varepsilon)j_3}   \|v_{\lambda_0}\|_{V^{2}_{\Delta}} \|u\|_{X^{1/2}_{p,\Delta}(-\lambda_0-I_{j_3})} \|u\|^2_{X^{1/2}_{p,\Delta}} \notag\\
& \lesssim  T^{2\varepsilon/3}  \sum_{ j_3, \ \lambda_0\ll 0}  2^{(-2/p+ 2\varepsilon)j_3}   \|v_{\lambda_0}\|_{V^{2}_{\Delta}} \|u\|_{X^{1/2}_{p,\Delta}(-\lambda_0-I_{j_3})} \|u\|^2_{X^{1/2}_{p,\Delta}} \notag\\
 & \lesssim T^{ \varepsilon/2}\|u\|^{3}_{X^{1/2}_{p, \Delta}}\|v\|_{Y_{p', \Delta}}.  \label{trilinear51}
\end{align}
In the case $u_{-I_{j_1}}$ attaining the highest dispersion modulation, we have from \eqref{dispersiondecay} and \eqref{highestmod50} that
\begin{align}
\mathscr{L}^{-,1}_{l_-l_-h_+} (u,v)
& \lesssim \!\!\!\!\!\!  \sum_{\lambda_0\ll 0, \  j_1 \ll  j_2 \thickapprox j_3 } \!\!\!\!\!\!   \langle \lambda_0 \rangle ^{1/2}   \|v_{\lambda_0}\|_{ L^{\infty}_{x,t}}\|u_{-I_{j_1}}\|_{L^{2}_{x,t}}\|u_{-I_{j_2}}\|_{L^{4}_{x,t}}\|\partial_x u_{-\lambda_0-I_{j_3}}\|_{L^{4}_{x, t}} \notag\\
& \lesssim \!\!\!\!\!\!  \sum_{\lambda_0\ll 0, \  j_1 \ll  j_2 \thickapprox j_3 } \!\!\!\!\!\!   \langle \lambda_0 \rangle ^{-1/2}   \|v_{\lambda_0}\|_{ V^2_\Delta}\|u_{-I_{j_1}}\|_{V^2_\Delta}\|u_{-I_{j_2}}\|_{L^{4}_{x,t}}\|\partial_x u_{-\lambda_0-I_{j_3}}\|_{L^{4}_{x, t}},
 \label{trilinear52}
\end{align}
which is the same as in the right hand side of the second inequality as in \eqref{trilinear51}.

If $u_{-I_{j_2}}$ has the highest modulation, we have
\begin{align}
\mathscr{L}^{-,1}_{l_-l_-h_+} (u,v)
& \lesssim \!\!\!\!\!\!  \sum_{\lambda_0\ll 0, \  j_1 \ll  j_2 \thickapprox j_3 } \!\!\!\!\!\!   \langle \lambda_0 \rangle ^{1/2}   \|v_{\lambda_0}\|_{ L^{\infty}_{x,t}}\|u_{-I_{j_2}}\|_{L^{2}_{x,t}}\|u_{-I_{j_1}} \partial_x \overline{u}_{-\lambda_0-I_{j_3}}\|_{L^{2}_{x, t}}. \label{trilinear53}
\end{align}
Noticing that for any $\lambda_1\in -I_{j_1}, \  \lambda_3 \in -\lambda_0-I_{j_3}$, if $j_1\ll j_3$,   we have $|\lambda_3|- |\lambda_1| \gtrsim |\lambda_0| $.
Using \eqref{highestmod2}, the bilinear estimate \eqref{bilinear}, Corollaries \ref{L4toXj} and \ref{L4toXjlambda},  and noticing that $j_3 \lesssim \ln \langle \lambda_0\rangle$,  we have
\begin{align}
\mathscr{L}^{-,1}_{l_-l_-h_+} (u,v)
& \lesssim T^{\varepsilon/4} \!\!\!\!\!\! \sum_{\lambda_0\ll 0, \  j_1 \ll  j_3 } \!\!\!\!\!\!      2^{-j_1/p} 2^{ j_3(-2/p+ \varepsilon)}  \|v_{\lambda_0}\|_{V^2_\Delta} \|u \|_{X^{1/2}_{p,\Delta}(-\lambda_0-I_{j_3})} \|u \|^2_{X^{1/2}_{p,\Delta}}. \label{trilinear54}
\end{align}
$\mathscr{L}^{-,2}_{l_-l_-h_+} (u,v)$ and $\mathscr{L}^{-,3}_{l_-l_-h_+} (u,v)$ can be estimated in an analogous way as above, we omit the details.

{\it Case $l_-a_+a_m$.} In view of the frequency constraint condition we see that $\lambda_2 \in [0, -3\lambda_0/4+C)$ and $\lambda_3\in [-\lambda_0/4-C, -\lambda_0)$.  We consider the dyadic collections of $\lambda_k$:
$$
\lambda_1 \in [\frac{3\lambda_0}{4}, 0) = \bigcup_{j_1} -I_{j_1},  \ \ \lambda_2 \in [0, -\frac{3\lambda_0}{4} +C) = \bigcup_{j_2}  I_{j_2},  \ \ \lambda_3 \in [-\frac{\lambda_0}{4}-C, -\lambda_0) = \bigcup_{ j_3\geq 0}(-\lambda_0  -I_{j_3}).
$$
It is easy to see that $j_1 \thickapprox (j_2\vee j_3) \lesssim \ln \langle \lambda_0\rangle$.  We need to estimate
\begin{align}
\mathscr{L}^{-}_{l_-a_+h_m} (u,v) &  :=   \sum_{  \lambda_0 \ll 0 } \left\langle\lambda_0\right\rangle^{1/2}\left|\int_{\left[0,T\right] \times \mathbb{R}}\overline{v}_{\lambda_0}  u_{[3\lambda_0/4,0)}   u_{[0, -3\lambda_0/4+C)}   \partial_x \overline{u}_{[ -\lambda_0/4 -C , -\lambda_0)} \,dxdt \right| \notag\\
 &  \lesssim   \left(\sum_{\lambda_0\ll 0, \ j_3\ll  j_2 \thickapprox j_1 }  + \sum_{\lambda_0\ll 0, \ j_2\ll  j_1 \thickapprox j_3  } + \sum_{\lambda_0\ll 0, \ j_1 \thickapprox j_2 \thickapprox j_3  } \right)   \langle\lambda_0 \rangle^{1/2} \nonumber \\
& \ \ \ \ \ \ \ \ \ \ \ \ \ \ \ \ \ \ \ \ \ \ \ \times \int_{\left[0,T\right] \times \mathbb{R}} |\overline{v}_{\lambda_0} u_{-I_{j_1}} u_{I_{j_2}}  \partial_x \overline{u}_{-\lambda_0-I_{j_3}}| \,dxdt. \label{trilinear55}
 \end{align}
Using a similar way as in the estimate of \eqref{trilinear49}, we can get the desired result and we omit the details.

{\it Case $l_-a_+h_+$.}   We consider the dyadic collections of $\lambda_k$:
$$
\lambda_1 \in [\frac{3\lambda_0}{4}, 0) = \bigcup_{j_1} -I_{j_1},  \ \ \lambda_2 \in [0, \infty) = \bigcup_{j_2\geq 0}  I_{j_2},  \ \ \lambda_3 \in [ -\lambda_0, \infty) = \bigcup_{ j_3\geq 0}(-\lambda_0  +I_{j_3}).
$$
From FCC \eqref{E:frequencyrelation} it follows that $j_2 \thickapprox (j_1\vee j_3) $.  We need to estimate
\begin{align}
\mathscr{L}^{-}_{l_-a_+h_+} (u,v) &  :=   \sum_{  \lambda_0 \ll 0 } \left\langle\lambda_0\right\rangle^{1/2}\left|\int_{\left[0,T\right] \times \mathbb{R}}\overline{v}_{\lambda_0}  u_{[3\lambda_0/4,0)}   u_{[0, \infty)}   \partial_x \overline{u}_{[ -\lambda_0, \infty)} \,dxdt \right| \notag\\
 &  \lesssim   \left(\sum_{\lambda_0\ll 0, \ j_1\ll  j_2 \thickapprox j_3 }  + \sum_{\lambda_0\ll 0, \ j_3 \ll  j_1 \thickapprox j_2 } + \sum_{\lambda_0\ll 0, \ j_1 \thickapprox j_2 \thickapprox j_3  } \right)   \langle\lambda_0 \rangle^{1/2} \nonumber \\
& \ \ \ \ \ \times \int_{\left[0,T\right] \times \mathbb{R}} |\overline{v}_{\lambda_0} u_{-I_{j_1}} u_{I_{j_2}}  \partial_x \overline{u}_{-\lambda_0+ I_{j_3}}| \,dxdt:= \Gamma_1+ \Gamma_2+\Gamma_3. \label{trilinear56}
 \end{align}
 In the right hand side of \eqref{trilinear56}, by DMCC \eqref{E:modulationrelation} the highest dispersion modulation satisfies
\begin{equation}\label{highestmod57}
\max_{0\leq k \leq 3} |\xi^2_k +\tau_k|   \gtrsim   \langle \lambda_0\rangle (\langle \lambda_0\rangle  +2^{j_2}) .
\end{equation}
If $v_{\lambda_0}$ has the highest dispersion modulation, using a similar way as in \eqref{trilinear16}--\eqref{trilinear18},  we have
\begin{align}
\Gamma_1
& \lesssim \!\!\!\!\!\!  \sum_{\lambda_0\ll 0, \  j_1 \ll  j_2 \thickapprox j_3 } \!\!\!\!\!\!   \langle \lambda_0 \rangle ^{1/2}   \|v_{\lambda_0}\|_{L^{2}_{t}L^{\infty}_{x}}\|u_{-I_{j_1}}\|_{L^{\infty}_{t}L^{2}_{x}}\|u_{ I_{j_2}}\|_{L^{4}_{x,t}}\|\partial_x u_{-\lambda_0+I_{j_3}}\|_{L^{4}_{x, t}} \notag\\
& \lesssim \!\!\!\!\!\!  \sum_{\lambda_0\ll 0, \  j_1 \ll  j_2 \thickapprox j_3 } \!\!\!\!\!\!   (\langle \lambda_0 \rangle + 2^{j_2} )^{-1/2}   \|v_{\lambda_0}\|_{V^2_\Delta}\|u_{-I_{j_1}}\|_{V^2_\Delta}\|u_{ I_{j_2}}\|_{L^{4}_{x,t}}\|\partial_x u_{-\lambda_0+I_{j_3}}\|_{L^{4}_{x, t}} \notag\\
& \lesssim  T^{2\varepsilon/3}  \sum_{  j_1 \ll  j_3, \ \lambda_0\ll 0} 2^{-j_1/p}  2^{(-2/p+ 2\varepsilon)j_3}   \|v_{\lambda_0}\|_{V^{2}_{\Delta}} \|u\|_{X^{1/2}_{p,\Delta}(-\lambda_0+I_{j_3})} \|u\|^2_{X^{1/2}_{p,\Delta}} \notag\\
  & \lesssim T^{ \varepsilon/2}\|u\|^{3}_{X^{1/2}_{p, \Delta}}\|v\|_{Y_{p', \Delta}}.  \label{trilinear58}
\end{align}
If $u_{-I_{j_1}}$ has the highest dispersion modulation, the estimate of $\Gamma_1$ is similar to \eqref{trilinear58} by considering
\begin{align}
\Gamma_1
& \lesssim \!\!\!\!\!\!  \sum_{\lambda_0\ll 0, \  j_1 \ll  j_2 \thickapprox j_3 } \!\!\!\!\!\!   \langle \lambda_0 \rangle ^{1/2}   \|v_{\lambda_0}\|_{ L^{\infty}_{x,t}}\|u_{-I_{j_1}}\|_{ L^{2}_{x,t}}\|u_{ I_{j_2}}\|_{L^{4}_{x,t}}\|\partial_x u_{-\lambda_0+I_{j_3}}\|_{L^{4}_{x, t}} \notag\\
& \lesssim \!\!\!\!\!\!  \sum_{\lambda_0\ll 0, \  j_1 \ll  j_2 \thickapprox j_3 } \!\!\!\!\!\!   (\langle \lambda_0 \rangle + 2^{j_2} )^{-1/2}   \|v_{\lambda_0}\|_{V^2_\Delta}\|u_{-I_{j_1}}\|_{V^2_\Delta}\|u_{ I_{j_2}}\|_{L^{4}_{x,t}}\|\partial_x u_{-\lambda_0+I_{j_3}}\|_{L^{4}_{x, t}},
\label{trilinear59}
\end{align}
for which the right hand side of \eqref{trilinear59} is the same one as the second inequality of \eqref{trilinear58}.

If $u_{I_{j_2}}$  gains the highest dispersion modulation, using \eqref{highestmod57},  the bilinear estimate \eqref{bilinear} and Corollaries \ref{L4toXj}, \ref{L4toXjlambda},  we have
\begin{align}
\Gamma_1
& \lesssim \!\!\!\!\!\!  \sum_{\lambda_0\ll 0, \  j_1 \ll  j_2 \thickapprox j_3 } \!\!\!\!\!\!   \langle \lambda_0 \rangle ^{1/2}   \|v_{\lambda_0}\|_{ L^{\infty}_{x,t}} \|u_{ I_{j_2}}\|_{L^{2}_{x,t}}  \|u_{-I_{j_1}} \partial_x \overline{u}_{-\lambda_0+I_{j_3}}\|_{L^{2}_{x, t}} \notag\\
& \lesssim \!\!\!\!\!\!  \sum_{\lambda_0\ll 0, \  j_1 \ll  j_2 \thickapprox j_3 } \!\!\!\!\!\!   (\langle \lambda_0 \rangle + 2^{j_2} )^{\varepsilon}   \|v_{\lambda_0}\|_{V^2_\Delta}\|u_{-I_{j_1}}\|_{V^2_\Delta} \|u_{ I_{j_2}}\|_{V^2_\Delta}\| u_{-\lambda_0+I_{j_3}}\|_{V^2_\Delta} \notag\\
& \lesssim  T^{ \varepsilon/4}  \sum_{  j_1 \ll  j_3, \ \lambda_0\ll 0} 2^{-j_1/p}  2^{(-2/p+ \varepsilon)j_3}   \|v_{\lambda_0}\|_{V^{2}_{\Delta}} \|u\|_{X^{1/2}_{p,\Delta}(-\lambda_0+I_{j_3})} \|u\|^2_{X^{1/2}_{p,\Delta}} \notag\\
  & \lesssim T^{ \varepsilon/4}\|u\|^{3}_{X^{1/2}_{p, \Delta}}\|v\|_{Y_{p', \Delta}}.  \label{trilinear60}
\end{align}
 If $u_{-\lambda_0+ I_{j_3}}$  has the highest dispersion modulation, one can use similar way as in \eqref{trilinear60} to get the same estimate and we omit the details.

{\it Case $a_+a_+h$.} Finally, we consider the following case as in Table 3.

\begin{table}[h]
\begin{center}

\begin{tabular}{|c|c|c|c|}
\hline
 ${\rm Case}$   &  $\lambda_1\in $ & $\lambda_2\in $ & $\lambda_3\in $    \\
\hline
$a_+a_+h$  & $ [0, \infty)$  & $ [0, \infty)$  &  $[\lambda_0, \infty)$  \\
\hline
\end{tabular}
\end{center}
\caption{ $\lambda_1,\lambda_2,\lambda_3 \geq 0$  .}
\end{table}

In view of FCC \eqref{E:frequencyrelation} we see that $\lambda_3 \in [-\lambda_0-C, \infty)$.
We consider the dyadic collections of $\lambda_k$:
$$
\lambda_k \in [ 0, \infty) = \bigcup_{j_k}  I_{j_k},  \ k=1,2;  \ \ \lambda_3 \in [-\lambda_0-C, \infty)  = \bigcup_{ j_3\geq -1}(-\lambda_0  +I_{j_3}), \ \ j_{-1}=[-C,0).
$$
From FCC \eqref{E:frequencyrelation} it follows that $j_3 \thickapprox (j_1\vee j_2) $.  We need to estimate
\begin{align}
\mathscr{L}^{-}_{l_-a_+h_+} (u,v) &  :=   \sum_{  \lambda_0 \ll 0 } \left\langle\lambda_0\right\rangle^{1/2}\left|\int_{\left[0,T\right] \times \mathbb{R}}\overline{v}_{\lambda_0}  u_{[0, \infty)}   u_{[0, \infty)}   \partial_x \overline{u}_{[ -\lambda_0-C, \infty)} \,dxdt \right| \notag\\
 &  \lesssim   \left(\sum_{\lambda_0\ll 0, \ j_1\ll  j_2 \thickapprox j_3 }  + \sum_{\lambda_0\ll 0, \ j_2 \ll  j_1 \thickapprox j_3 } + \sum_{\lambda_0\ll 0, \ j_1 \thickapprox j_2 \thickapprox j_3  } \right)   \langle\lambda_0 \rangle^{1/2} \nonumber \\
& \ \ \ \ \ \times \int_{\left[0,T\right] \times \mathbb{R}} |\overline{v}_{\lambda_0} u_{I_{j_1}} u_{I_{j_2}}  \partial_x \overline{u}_{-\lambda_0+ I_{j_3}}| \,dxdt:= \Upsilon_1+ \Upsilon_2+\Upsilon_3. \label{trilinear61}
 \end{align}
It suffices to estimate $\Upsilon_1$. The highest dispersion modulation in the right hand side of $\Upsilon_1$ satisfies
\begin{equation}\label{highestmod62}
\max_{0\leq k \leq 3} |\xi^2_k +\tau_k|   \gtrsim   (\langle \lambda_0\rangle  +2^{j_1})  (\langle \lambda_0\rangle  +2^{j_2}) .
\end{equation}
Now we compare the highest dispersion modulation between \eqref{highestmod57} and \eqref{highestmod62},  we see that the the highest modulation in $\Upsilon_1$ is larger than that of $\Gamma_1$. If $v_{\lambda_0}$ has the highest dispersion modulation,  using the same estimates as in \eqref{trilinear58} for $u_{I_{j_1}}, \ u_{I_{j_2}}, \  \partial_x \overline{u}_{-\lambda_0+ I_{j_3}}$ in the spaces $V^2_\Delta, \ L^4_{x,t}$, we can obtain that $\Upsilon_1$ has the same upper bound as that of $\Gamma_1$  in \eqref{trilinear58}. The other cases are also similar to those estimates in \eqref{trilinear59} and \eqref{trilinear60}.

{\bf Step 2.} We consider the case that $\lambda_0$ is the secondly minimal integer in $\lambda_0,...,\lambda_3$. Namely, there is a bijection $\pi: \{1,2,3\} \to \{1,2,3\}$ such that
\begin{align}
\lambda_{\pi(1)} \leq \lambda_0 \leq \lambda_{\pi(2)} \vee \lambda_{\pi(3)}.  \label{trilinear64}
\end{align}
{\bf Step 2.1.} We assume that $\lambda_0 \gg 0$ and $\lambda_1=\lambda_{\pi(1)}$.  By \eqref{trilinear64}, we need to estimate
\begin{align}
\mathscr{L} (u,v) &  :=   \sum_{  \lambda_0 \gg  0 } \left\langle\lambda_0\right\rangle^{1/2} \int_{\left[0,T\right] \times \mathbb{R}}\overline{v}_{\lambda_0}  u_{(-\infty, \lambda_0)}   u_{[\lambda_0, \infty)}   \partial_x \overline{u}_{[\lambda_0, \infty)} \,dxdt  \notag\\
& =   \sum_{  \lambda_0 \gg  0 } \left\langle\lambda_0\right\rangle^{1/2} \int_{\left[0,T\right] \times \mathbb{R}}\overline{v}_{\lambda_0}  u_{(-\infty, -\lambda_0/2)}   u_{[\lambda_0, \infty)}   \partial_x \overline{u}_{[\lambda_0, \infty)} \,dxdt \notag\\
    & \ \ + \sum_{  \lambda_0 \gg  0 } \left\langle\lambda_0\right\rangle^{1/2} \int_{\left[0,T\right] \times \mathbb{R}}\overline{v}_{\lambda_0}  u_{[-\lambda_0/2, \lambda_0/2)}   u_{[\lambda_0, \infty)}   \partial_x \overline{u}_{[\lambda_0, \infty)} \,dxdt \notag\\
 & \ \   + \sum_{  \lambda_0 \gg  0 } \left\langle\lambda_0\right\rangle^{1/2} \int_{\left[0,T\right] \times \mathbb{R}}\overline{v}_{\lambda_0}  u_{[\lambda_0/2, \lambda_0]}   u_{[\lambda_0, \infty)}   \partial_x \overline{u}_{[\lambda_0, \infty)} \,dxdt \notag\\
& := \mathscr{L}^{-} (u,v) + \mathscr{L}^{0} (u,v) + \mathscr{L}^{+} (u,v)
. \label{trilinear65}
 \end{align}
First, we estimate $\mathscr{L}^{0} (u,v)$. From FCC \eqref{E:frequencyrelation} it follows that
\begin{align}
\mathscr{L}^{0} (u,v) & =   \sum_{  \lambda_0 \gg  0 } \left\langle\lambda_0\right\rangle^{1/2} \int_{\left[0,T\right] \times \mathbb{R}}\overline{v}_{\lambda_0}  u_{[-\lambda_0/2, 0)}   u_{[2\lambda_0 -C, \infty)}   \partial_x \overline{u}_{[\lambda_0, \infty)} \,dxdt \notag\\
& \quad  +  \sum_{  \lambda_0 \gg 0 } \left\langle\lambda_0\right\rangle^{1/2} \int_{\left[0,T\right] \times \mathbb{R}}\overline{v}_{\lambda_0}  u_{[0, \lambda_0/2)}   u_{[3\lambda_0/2-C, \infty)}   \partial_x \overline{u}_{[\lambda_0, \infty)} \,dxdt \notag\\
& := \mathscr{L}^{0,-} (u,v) + \mathscr{L}^{0,+} (u,v). \label{trilinear66}
 \end{align}
We further divide $\mathscr{L}^{0,+} (u,v)$ into two parts
\begin{align}
\mathscr{L}^{0,+} (u,v) &  =  \sum_{  \lambda_0 \gg  0 } \left\langle\lambda_0\right\rangle^{1/2} \int_{\left[0,T\right] \times \mathbb{R}}\overline{v}_{\lambda_0}  u_{[0, \lambda_0/2)}   u_{[3\lambda_0/2-C, 2\lambda_0)}   \partial_x \overline{u}_{[\lambda_0, \infty)} \,dxdt \notag\\
 &  \quad \quad +  \sum_{  \lambda_0 \gg  0 } \left\langle\lambda_0\right\rangle^{1/2} \int_{\left[0,T\right] \times \mathbb{R}}\overline{v}_{\lambda_0}  u_{[0, \lambda_0/2)}   u_{[2\lambda_0, \infty)}   \partial_x \overline{u}_{[\lambda_0, \infty)} \,dxdt \notag\\
& := \mathscr{L}^{0,+}_l (u,v) + \mathscr{L}^{0,+}_h (u,v)
. \label{trilinear67}
 \end{align}
Again, in view of FCC \eqref{E:frequencyrelation} we have
\begin{align}
\mathscr{L}^{0,+}_l  (u,v) &  =  \sum_{  \lambda_0 \gg  0 } \left\langle\lambda_0\right\rangle^{1/2} \int_{\left[0,T\right] \times \mathbb{R}}\overline{v}_{\lambda_0}  u_{[0, \lambda_0/2)}   u_{[3\lambda_0/2-C, 2\lambda_0)}   \partial_x \overline{u}_{[\lambda_0, 3\lambda_0/2+C)} \,dxdt. \label{trilinear68}
 \end{align}
In the right hand side of $\mathscr{L}^{0,+}_l  (u,v)$, by DMCC \eqref{E:modulationrelation} we have
$$
\bigvee^3_{k=0} |\xi^2_k + \tau_k| \gtrsim \langle \lambda_0\rangle^2.
$$
By Corollary \ref{L4toXj} we have
\begin{align}
 \|u_{([, \lambda_0/2)}\|_{L^\infty_t L^2_x \cap V^2_\Delta} \leq \sum_{j\leq \log_2^{\lambda_0} +1} \|u_{I_j}\|_{V^2_\Delta}  \lesssim \|u\|_{X^{1/2}_{p,\Delta}} .   \label{trilinear69}
\end{align}
If $v_{\lambda_0}$ gains the highest dispersion modulation, by H\"older's inequality, the dispersion modulation decay estimate \eqref{dispersiondecay}, \eqref{trilinear69} and Lemma \ref{L4toXlambda}, we have for $0<\varepsilon <1/4p$,
\begin{align}
|\mathscr{L}^{0,+}_l  (u,v)| &   \lesssim   \sum_{  \lambda_0 \gg  0 } \langle\lambda_0 \rangle^{-1/2}  \|{v}_{\lambda_0}\|_{V^2_\Delta}  \|u_{[0, \lambda_0/2)}\|_{L^\infty_t L^2_x}   \|u_{[3\lambda_0/2-C, 2\lambda_0)} \|_{L^4_{x,t}}  \|\partial_x \overline{u}_{[\lambda_0, 3\lambda_0/2+C)}\|_{L^4_{x,t}} \notag\\
&   \lesssim  T^{\varepsilon/2} \sum_{  \lambda_0 \gg  0 } \langle\lambda_0 \rangle^{-2/p + 2\varepsilon}  \|{v}_{\lambda_0}\|_{V^2_\Delta}      \|u  \|^3_{X^{1/2}_{p,\Delta} }   \lesssim   T^{\varepsilon/2} \|u  \|^3_{X^{1/2}_{p,\Delta} }. \label{trilinear70}
 \end{align}
If $u_{[0, \lambda_0/2)}$ the  highest dispersion modulation, taking $L^\infty_{x,t}, L^2_{x,t}, L^4_{x,t}, L^4_{x,t}$ norms to ${v}_{\lambda_0}$, $ u_{[0, \lambda_0/2)},$ $ u_{[3\lambda_0/2-C, 2\lambda_0)}, $ $\partial_x \overline{u}_{[\lambda_0, 3\lambda_0/2+C)}$, then applying \eqref{dispersiondecay},  we have
\begin{align}
|\mathscr{L}^{0,+}_l  (u,v)|
 & \lesssim   \sum_{  \lambda_0 \gg  0 } \langle\lambda_0 \rangle^{-1/2}  \|{v}_{\lambda_0}\|_{V^2_\Delta}  \|u_{[0, \lambda_0/2)}\|_{V^2_\Delta}   \|u_{[3\lambda_0/2-C, 2\lambda_0)} \|_{L^4_{x,t}}  \|\partial_x \overline{u}_{[\lambda_0, 3\lambda_0/2+C)}\|_{L^4_{x,t}}, \notag
 \end{align}
which reduces the same estimate as the first inequality in \eqref{trilinear70}.

Let $u_{[3\lambda_0/2-C, 2\lambda_0)}$ have the highest dispersion modulation. Taking $L^\infty_{x,t}, L^2_{x,t}, L^2_{x,t}$ norms to ${v}_{\lambda_0}$,  $ u_{[3\lambda_0/2-C, 2\lambda_0)}, $ $u_{[0, \lambda_0/2)}\partial_x \overline{u}_{[\lambda_0, 3\lambda_0/2+C)}$, then applying the dispersion modulation decay \eqref{dispersiondecay} to $u_{[3\lambda_0/2-C, 2\lambda_0)}$ and the bilinear estimate \eqref{bilinear} to  $u_{(0, \lambda_0/2]} \partial_x \overline{u}_{[\lambda_0, 3\lambda_0/2+C)}$, one obtains that $\mathscr{L}^{0,+}_l  (u,v)$ has the desired estimate. When $\partial_x u_{[\lambda_0, 3\lambda_0/2+C)}$ has the highest dispersion modulation, the argument is similar.

Now we estimate $\mathscr{L}^{0,+}_h  (u,v)$. We adopt the following decompositions:
$$
\lambda_1 \in (0,\lambda_0/2) =\bigcup_{j_1>0} I_{j_1}, \ \ \lambda_2 \in [2\lambda_0, \infty) = \bigcup_{j_2\geq 0} (2\lambda_0 +I_{j_2}), \ \  \lambda_3 \in [\lambda_0, \infty) = \bigcup_{j_3\geq 0} (\lambda_0 +I_{j_3}).
$$
By FCC \eqref{E:frequencyrelation}, we have
$$
j_3 \thickapprox j_1 \vee j_2, \ \ j_1\leq \log_2^{\lambda_0} +1.
$$
\begin{align}
|\mathscr{L}^{0,+}_h (u,v)| &   \leq   \sum_{  \lambda_0 \gg  0, \ j_3 \thickapprox j_1 \vee j_2 } \left\langle\lambda_0\right\rangle^{1/2} \int_{\left[0,T\right] \times \mathbb{R}} |\overline{v}_{\lambda_0}  u_{I_{j_1}}   u_{ 2\lambda_0 +I_{j_2}}   \partial_x \overline{u}_{ \lambda_0 + I_{j_3}}| \,dxdt \notag\\
&   \leq   \left( \sum_{  \lambda_0 \gg  0, \ j_1\ll  j_3 \thickapprox  j_2 } + \sum_{  \lambda_0 \gg  0, \ j_2 \ll  j_3 \thickapprox  j_1 } + \sum_{  \lambda_0 \gg  0, \ j_1 \thickapprox   j_2 \thickapprox  j_3 } \right) \left\langle\lambda_0\right\rangle^{1/2} \notag\\
 & \ \ \ \ \times \int_{\left[0,T\right] \times \mathbb{R}} |\overline{v}_{\lambda_0}  u_{I_{j_1}}   u_{ 2\lambda_0 +I_{j_2}}   \partial_x \overline{u}_{ \lambda_0 + I_{j_3}}| \,dxdt = I+II+III. \label{trilinear71}
 \end{align}
In the right hand side of \eqref{trilinear71}, by DMCC \eqref{E:modulationrelation} we have
\begin{align}
\bigvee^3_{k=0} |\xi^2_k + \tau_k| \gtrsim \langle \lambda_0\rangle (\langle\lambda_0 \rangle + 2^{j_2}). \label{trilinear72}
\end{align}
 Noticing that \eqref{trilinear72} is the same as \eqref{highestmod57}, we can follow the same ideas in the estimates of $\Gamma_1$  to get the bound of $I$, see \eqref{trilinear58}--\eqref{trilinear60}.

Let us observe that $j_1, j_2, j_3  \lesssim \ln \lambda_0$ in $II$ and $III$. It follows that the estimates of $II$ and $III$ are easier than that of $I$. We omit the details.

$\mathscr{L}^{0,-} (u,v)$ can be handled  in a similar way as that of  $\mathscr{L}^{0,+}_h  (u,v)$ and we omit the details.

Now we estimate $\mathscr{L}^{+} (u,v)$. We use the dyadic decomposition
$$
\lambda_1 \in (\lambda_0/2, \lambda_0) =\bigcup_{j_1>0}(\lambda_0 - I_{j_1}), \ \ \lambda_k \in [\lambda_0, \infty) = \bigcup_{j_k\geq 0} (\lambda_0 +I_{j_k}), \ \ k=2,3.
$$
By FCC \eqref{E:frequencyrelation}, we have
$$
j_2 \thickapprox j_1 \vee j_3 \ \ j_1\leq \log_2^{\lambda_0} +1.
$$
\begin{align}
|\mathscr{L}^{+} (u,v)| &   \leq   \sum_{  \lambda_0 \gg  0, \ j_2 \thickapprox j_1 \vee j_3 } \left\langle\lambda_0\right\rangle^{1/2} \int_{\left[0,T\right] \times \mathbb{R}} |\overline{v}_{\lambda_0}  u_{\lambda_0- I_{j_1}}   u_{\lambda_0 +I_{j_2}}   \partial_x \overline{u}_{ \lambda_0 + I_{j_3}}| \,dxdt. \label{trilinear75}
 \end{align}
Comparing \eqref{trilinear75} with \eqref{decompL}, we see that  $\mathscr{L}^{+}(u,v)$ in \eqref{trilinear75} is rather similar to \eqref{decompL} in Step 1.1. We omit the details of the proof.

We estimate $\mathscr{L}^{-}(u,v)$. In view of FCC \eqref{E:frequencyrelation} we see that
\begin{align}
\mathscr{L}^{-}  (u,v) &  =  \sum_{  \lambda_0 \gg  0 } \left\langle\lambda_0\right\rangle^{1/2} \int_{\left[0,T\right] \times \mathbb{R}}\overline{v}_{\lambda_0}  u_{(-\infty, \ -\lambda_0/2]}   u_{[5\lambda_0/2-C, \ \infty)}   \partial_x \overline{u}_{[\lambda_0, \ \infty)} \,dxdt. \label{trilinear76}
 \end{align}
Decompose $\lambda_k$ in the following dyadic way:
\begin{align*}
& \lambda_1 \in (-\infty, -\lambda_0/2] =\bigcup_{j_1\geq 0}(-\lambda_0/2 - I_{j_1}), \ \ \lambda_3 \in [\lambda_0, \infty) = \bigcup_{j_k\geq 0} (\lambda_0 +I_{j_3}) \\
& \lambda_2 \in [5\lambda_0/2 -C, \infty) =\bigcup_{j_2\geq -1}(5\lambda_0/2 + I_{j_2}), \ j_{-1}=[-C,0).
\end{align*}
It follows that $j_2 \thickapprox j_1\vee j_3$.
\begin{align}
|\mathscr{L}^{-}  (u,v)| & \leq  \sum_{  \lambda_0 \gg  0,  \ j_2 \thickapprox j_1\vee j_3} \left\langle\lambda_0\right\rangle^{1/2} \int_{\left[0,T\right] \times \mathbb{R}} |\overline{v}_{\lambda_0}  u_{ -\lambda_0/2-I_{j_1}}    u_{5\lambda_0/2 + I_{j_2}}   \partial_x \overline{u}_{\lambda_0 +I_{j_3}}| \,dxdt. \label{trilinear77}
 \end{align}
By DMCC \eqref{E:modulationrelation} we have
\begin{equation}\label{highestmod78}
\max_{0\leq k \leq 3} |\xi^2_k +\tau_k|   \gtrsim   (\langle \lambda_0\rangle  +2^{j_1})  (\langle \lambda_0\rangle  +2^{j_2}) .
\end{equation}
By \eqref{highestmod78}, we see that the dispersion modulation estimate \eqref{dispersiondecay} gives better decay. So, the estimate of  $\mathscr{L}^{-}  (u,v)$ is easier than that of the above cases and we will not perform the details.

{\bf Step 2.2.} We consider the case $\lambda_0 \ll 0$, $\lambda_1 <\lambda_0 \leq \lambda_2\vee \lambda_3$.  According to the size of $\lambda_2$, we divide the proof into four subcases $\lambda_2\in h_-,$ $\lambda_2\in l_-$, $\lambda_2 \in [0, -3\lambda_0/4)$ and $\lambda_2\in [-3\lambda_0/4, \infty)$. Moreover, in view of FCC \eqref{E:frequencyrelation},  in order to keep the left hand side of \eqref{trilinear4} nonzero,  it suffices to consider the following four subcases in Table 4.

\begin{table}[h]
\begin{center}

\begin{tabular}{|c|l|l|l|}
\hline
 ${\rm Case}$   &  $\lambda_1\in $ & $\lambda_2\in $ & $\lambda_3\in $    \\
\hline
$1$  & $ [5\lambda_0/4-C, \lambda_0)$  & $ [\lambda_0, 3\lambda_0/4)$  &  $[\lambda_0, 3\lambda_0/4 +C)$  \\
\hline
$2$ & $ [2\lambda_0 -C, \lambda_0)$  & $ [3\lambda_0/4, \ 0)$ & $ [\lambda_0, C)$  \\
\hline
$3$ & $ (11\lambda_0/4-C, \lambda_0)$  & $ [0, -3\lambda_0/4)$  & $ [\lambda_0,  -3\lambda_0/4 +C)$\\
\hline
$4$ & $ (-\infty, \lambda_0)$  & $ [-3\lambda_0/4, \infty)$  & $ [\lambda_0,  \infty)$\\
\hline
\end{tabular}
\end{center}
\caption{$\lambda_0 \ll 0$, \ \  $\lambda_1<\lambda_0 \leq \lambda_2\vee \lambda_3$. }
\end{table}

{\it Case 1.} Let us observe that all $\lambda_k$ ($k=1,2,3$) are localized in a neighbourhood of $\lambda_0$, which is essentially the same as in Case $h_-h_-a$ as in Step 1.2, we omit the details of the proof.

{\it Case 2.}  Similar to Case $h_-l_-l_-$ as in Step 1.2 and we omit its proof.

{\it Case 3.} By separating $\lambda_3 \in [\lambda_0,  -3\lambda_0/4 +C) = [\lambda_0, 0) \cup [0, -3\lambda_0/4 +C)$, we see that this case is quite similar to Case 2.

{\it Case 4.} Observing that $\lambda_1$ and $\lambda_2$ are far away from $0$, but $\lambda_3\in [\lambda_0, \infty)$ containing $0$, we need to further split $[\lambda_0,\infty)$ into $ [\lambda_0,0)$ and $ [0,\infty)$. We consider the following two subcases of Case 4, see Table 5.

\begin{table}[h]
\begin{center}

\begin{tabular}{|c|l|l|l|}
\hline
 ${\rm Case}$   &  $\lambda_1\in $ & $\lambda_2\in $ & $\lambda_3\in $    \\
\hline
$4.1$ & $ (-\infty, \lambda_0)$  & $ [-3\lambda_0/4, \infty)$  & $ [\lambda_0,  0)$\\
\hline
$4.2$ & $ (-\infty, \lambda_0)$  & $ [-3\lambda_0/4, \infty)$  & $ [0,  \infty)$\\
\hline
\end{tabular}
\end{center}
\caption{Two subcases of Case 4. }
\end{table}

{\it Case 4.1.} By FCC \eqref{E:frequencyrelation}, we have $\lambda_1\in (-\infty, 7\lambda_0/4 +C)$. So, we need to estimate
\begin{align}
\mathscr{L}^{2,\infty}  (u,v) & : =  \sum_{  \lambda_0 \ll  0 } \left\langle\lambda_0\right\rangle^{1/2} \int_{\left[0,T\right] \times \mathbb{R}}\overline{v}_{\lambda_0}  u_{(-\infty, \ 7\lambda_0/4 +C)}   u_{[-3\lambda_0/4, \ \infty)}   \partial_x \overline{u}_{[\lambda_0, \ 0)} \,dxdt. \label{trilinear79}
 \end{align}
We can decompose $\lambda_k$ in the following way:
\begin{align*}
& \lambda_1\in (-\infty, 7\lambda_0/4 +C) =\bigcup_{j_1\geq -1} (7\lambda_0/4 -I_{j_1}), I_{-1} =[-C,0),  \\
& \lambda_2 \in ( -3\lambda_0/4, \infty) =\bigcup_{j_1\geq 0} (-3\lambda_0/4 +I_{j_1}), \ \ \lambda_ 3 \in [\lambda_0, 0) =\bigcup_{j_3\geq 0} -I_{j_3}.
\end{align*}
It follows that $j_1 \thickapprox j_2 \vee j_3$, $j_3 \leq 1+ \log_2^{|\lambda_0|}$
\begin{align}
|\mathscr{L}^{2,\infty}  (u,v)| & \leq  \sum_{  \lambda_0 \ll  0,  \ j_1 \thickapprox j_2\vee j_3} \left\langle\lambda_0\right\rangle^{1/2} \int_{\left[0,T\right] \times \mathbb{R}} |\overline{v}_{\lambda_0}  u_{ 7\lambda_0/4-I_{j_1}}    u_{-3\lambda_0/4 + I_{j_2}}   \partial_x \overline{u}_{-I_{j_3}}| \,dxdt. \label{trilinear80}
 \end{align}
By DMCC \eqref{E:modulationrelation} we have
\begin{equation}\label{highestmod81}
\max_{0\leq k \leq 3} |\xi^2_k +\tau_k|   \gtrsim   (\langle \lambda_0\rangle  +2^{j_1})  (\langle \lambda_0\rangle  +2^{j_2}) .
\end{equation}
One can imitate the procedure as in Case $l_-a_+h_+$ in Step 1.2 to obtain the result, as desired.

{\it Case 4.2.}  We need to estimate
\begin{align}
\mathscr{L}^{3,\infty}  (u,v) & : =  \sum_{  \lambda_0 \ll  0 } \left\langle\lambda_0\right\rangle^{1/2} \int_{\left[0,T\right] \times \mathbb{R}}\overline{v}_{\lambda_0}  u_{(-\infty, \ \lambda_0)}   u_{[-3\lambda_0/4, \ \infty)}   \partial_x \overline{u}_{[0, \ \infty)} \,dxdt. \label{trilinear82}
 \end{align}
We can decompose $\lambda_k$ in the following way:
\begin{align*}
& \lambda_1\in (-\infty, \lambda_0) =\bigcup_{j_1>0} (\lambda_0  -I_{j_1}),   \\
& \lambda_2 \in ( -3\lambda_0/4, \infty)  = \bigcup_{j_2 }  I_{j_2}, \ \ \lambda_ 3 \in [0, \infty) =\bigcup_{j_3\geq 0} I_{j_3}.
\end{align*}
It follows that $j_2 \thickapprox j_1 \vee j_3$, $j_2 \geq   \log_2^{|\lambda_0|} -C$
\begin{align}
|\mathscr{L}^{3,\infty}  (u,v)| & \leq  \sum_{  \lambda_0 \ll  0,  \ j_1\vee j_3\thickapprox j_2 \geq \log^{|\lambda_0|}_2 -C} \left\langle\lambda_0\right\rangle^{1/2} \int_{\left[0,T\right] \times \mathbb{R}} |\overline{v}_{\lambda_0}  u_{ \lambda_0-I_{j_1}}    u_{ I_{j_2}}   \partial_x \overline{u}_{I_{j_3}}| \,dxdt. \label{trilinear83}
 \end{align}
By DMCC \eqref{E:modulationrelation} we have for $j_1\geq 10$,
\begin{equation}\label{highestmod84}
\max_{0\leq k \leq 3} |\xi^2_k +\tau_k|   \gtrsim       2^{j_1}   (\langle \lambda_0\rangle  +2^{j_2}) .
\end{equation}
Comparing \eqref{highestmod84} with \eqref{highestmod35},  \eqref{trilinear83} and \eqref{trilinear32}, one can imitate the procedures as in the Case $h_-a_+a$ to obtain the result, as desired. The details are omitted.

{\bf Step 2.3.} We consider the remaining cases. If $\pi(1)=2$, namely $\lambda_2$ is the smallest one in all $\lambda_k$, by symmetry we see that it is the same as the case $\pi(1)=1$. If $\pi(1)=3$, it follows from FCC \eqref{E:frequencyrelation} we have $\lambda_0 \thickapprox \lambda_1 \thickapprox \lambda_2 \thickapprox \lambda_3$. So, the summation in the left hand side of \eqref{trilinear4} is essentially on $\lambda_0$, the other summations are finitely many. By H\"older's inequality we have the result, as desired.

{\bf Step 3.} We consider the case $|\lambda_0|  \lesssim 1$.  It follows that $\langle \lambda_0\rangle^{1/2} \sim 1 $  in the left hand side of \eqref{trilinear4}, which means that we have gained half order derivative. So, this case becomes easier to handle and the details of the proof are omitted.

\section{Quintic linear estimates}

For the sake of convenience, we denote for $w=(w^{(0)},...,w^{(5)}) $,
\begin{align}
 \mathscr{L}(w) = \sum_{\lambda_0,...,\lambda_5\in \mathbb{Z}} \langle \lambda_0\rangle^{1/2} \int_{\mathbb{R}\times [0,T]} \prod^5_{k=0} w^{(k)}_{\lambda_k} (x,t) dxdt, \label{6linear}
\end{align}
In the following we always assume that
\begin{align}
w^{(k)}_{\lambda_k}= u_{\lambda_k}, \ \ k=1,3,5; \ \ w^{(0)}_{\lambda_0}= \overline{v}_{\lambda_0}, \ w^{(k)}= \overline{u}_{\lambda_k}, \ k=2,4, \label{6linear1}
\end{align}
which means that for $w=(\bar{v},u,\bar{u},u,\bar{u},u)$,
\begin{align}
 \mathscr{L}(w)&  = \sum_{\lambda_0,...,\lambda_5\in \mathbb{Z}} \langle \lambda_0\rangle^{1/2} \int_{\mathbb{R}\times [0,T]} \overline{v}_{\lambda_0}u_{\lambda_1}  \overline{u}_{\lambda_2} u_{\lambda_3} \overline{u}_{\lambda_4} u_{\lambda_5} (x,t) dxdt  \nonumber \\
 &
 = \sum_{\lambda_0 \in \mathbb{Z}} \langle \lambda_0\rangle^{1/2} \int_{\mathbb{R}\times [0,T]} \overline{v}_{\lambda_0}|u|^4u    (x,t) dxdt. \label{6linear2}
\end{align}

Our goal in this section is to show that

\begin{lem}\label{nonestimate}
Let $p\in [4, \infty)$, $0<T<1$, $0<\varepsilon \ll 1$. Let $\mathscr{L}(w)$ be as in \eqref{6linear} and \eqref{6linear1}.   We have
\begin{align}
|\mathscr{L}(w) | \lesssim T^{\varepsilon}  \|v\|_{Y^{0}_{p',\Delta}} \|u\|^5_{X^{1/2}_{p,\Delta}}    \label{6linear3}
\end{align}
and by duality,
\begin{align}
\||u|^4u \|_{X^{1/2}_{p,\Delta}}  \lesssim T^{\varepsilon}  \|u\|^5_{X^{1/2}_{p,\Delta}}.     \label{6linear4}
\end{align}
\end{lem}

We divide the proof of Lemma \ref{nonestimate} into a few steps according to the size of $\lambda_0$.  We can assume that $\lambda_0>0$, since in the opposite case one can substitute $\lambda_0,...,\lambda_5$ by $-\lambda_0,..,-\lambda_5$.

{\bf Step 1.} Let us assume that $\lambda_0,...,\lambda_5$ satisfy
\begin{align}
\lambda_0 = \max_{0\leq k\leq 5} |\lambda_k|.  \label{6linear5}
\end{align}
For short, considering the higher, lower and all frequency of $\lambda_k$,  we use the following notations:
$$
\left\{
\begin{array}{l}
\lambda_k \in h \Leftrightarrow \lambda_k\in [c\lambda_0, \lambda_0] \\
\lambda_k \in h_- \Leftrightarrow \lambda_k\in [-\lambda_0, -c\lambda_0]\\
\lambda_k \in l \Leftrightarrow \lambda_k\in [0, c \lambda_0]\\
\lambda_k \in l_- \Leftrightarrow \lambda_k\in [-c\lambda_0, 0]\\
\lambda_k \in a \Leftrightarrow \lambda_k\in [-\lambda_0, \lambda_0]\\
\lambda_k \in a_+ \Leftrightarrow \lambda_k\in [0, \lambda_0]\\
\lambda_k \in a_- \Leftrightarrow \lambda_k\in [-\lambda_0, 0]
\end{array}
\right.
$$
for some $c>0$. First, we consider the case that there are two higher frequency in $\lambda_1,...,\lambda_5$, say, $\lambda_1,\lambda_3$ belong to higher frequency intervals.   We denote by $(\lambda_k)\in hhaaa$ that all $\lambda_0,...,\lambda_5$ satisfy conditions \eqref{6linear5} and
\begin{align}
\lambda_1, \lambda_3 \in h, \quad   \lambda_2, \lambda_4, \lambda_5  \in a.  \label{6linear6}
\end{align}
For $w=(\bar{v},u,\bar{u},u,\bar{u},u)$, we write
\begin{align}
 \mathscr{L}_{hhaaa}(w) = \sum_{(\lambda_k)\in hhaaa } \langle \lambda_0\rangle^{1/2} \int_{\mathbb{R}\times [0,T]} \prod^5_{k=0} w^{(k)}_{\lambda_k} (x,t) dxdt, \label{6linear8}
\end{align}
we will always use the notation
\begin{align}
 \mathscr{L}_{bcdef}(w) = \sum_{(\lambda_k)\in  bcdef} \langle \lambda_0\rangle^{1/2} \int_{\mathbb{R}\times [0,T]} \prod^5_{k=0} w^{(k)}_{\lambda_k} (x,t) dxdt, \label{6linear8a}
\end{align}
where $b,c,d,e,f \in \{a, a_+, a_-, h, h_-, l, l_-\}$.

\begin{lem}\label{nonestimate2h}
Let $p\in [4, \infty)$, $0<T<1$, $0<\varepsilon \ll 1$.    We have
\begin{align}
|\mathscr{L}_{bcdef}(w) | \lesssim T^{\varepsilon}  \|v\|_{Y^{0}_{p',\Delta}} \|u\|^5_{X^{1/2}_{p,\Delta}}     \label{6linear31}
\end{align}
if at least two ones in $``b,c,d,e,f"$ belong to $\{h, h_-\}$.
\end{lem}

{\bf Proof.} {\it Case $hha_+a_+a_+.$} We denote by $(\lambda_k) \in hha_+a_+a_+$ the case that  $\lambda_0,..., \lambda_5$ satisfy \eqref{6linear5} and
\begin{align}
 \lambda_k \in h, \ k=1,3, \ \ \lambda_2, \lambda_4, \lambda_5 \in a_+.     \label{6linear35}
\end{align}
We decompose $\lambda_2, \lambda_4, \lambda_5$ in a dyadic way:
\begin{align}
\lambda_k \in [0, \lambda_0]= \bigcup_{j_k\geq 0}  I_{j_k}, \ \  k=2,4,5; \     \ \ 0\leq j_2,  j_4,  j_5 \leq \log^{\lambda_0}_2.    \label{6linear36}
\end{align}
We can assume that $\lambda_0\gg 1$.  Assume that $j_{\max} =j_2\vee j_4\vee j_5$,  $j_{\min} =j_2\wedge j_4\wedge j_5$  and $j_{\rm med}$ is the secondly larger one in $\{j_2, j_4, j_5\}$   in \eqref{6linear36}.  By H\"older's inequality,
\begin{align}
 |\mathscr{L}_{hha_+a_+a_+ }(w)| \lesssim & \sum_{\lambda_0; \  j_{\min}, j_{\rm med}, j_{\max} \leq \log^{\lambda_0}_2 } \langle \lambda_0 \rangle^{1/2}  \|v_{\lambda_0}\|_{L^\infty_{x,t}} \prod_{k=1,3} \|w^{(k)}_{[c \lambda_0, \lambda_0]} \|_{L^4_{x,t\in [0,T]}} \nonumber\\
  & \times  \|u_{I_{j_{\max}}} \|_{L^4_{x,t\in [0,T]}}   \|u_{I_{j_{\rm med}}} \|_{L^4_{x,t\in [0,T]}}   \| u_{I_{j_{\min}}} \|_{L^\infty_{x,t\in [0,T]}}      \label{6linear37}
\end{align}
In view of $V^2_\Delta \subset L^\infty_t L^2_x $, $\|v_{\lambda_0}\|_{L^\infty_x}  \lesssim \|v_{\lambda_0}\|_{L^2_x}$,
  by Corollaries \ref{L4toXj} and \ref{L4toXlambda}, one has that
\begin{align}
 |\mathscr{L}_{hha_+a_+a_+ }(w)| \lesssim & T^{4\varepsilon/3}  \sum_{ \lambda_0; \  j_{\min}, j_{\rm med}, j_{\max} \leq \log^{\lambda_0}_2 }
  \langle \lambda_0\rangle^{1/2 }  \|v_{\lambda_0}\|_{V^2_\Delta}   \langle \lambda_0\rangle^{ 2(-1/4 -1/p +  \varepsilon)} \|u \|^2_{X^{1/2}_{p,\Delta}} \nonumber\\
  & \ \  \times 2^{(j_{\max}+j_{\rm med})(-1/4-1/p+ \varepsilon) } \|u \|^2_{X^{1/2}_{p,\Delta}}
 2^{j_{\min}(1/2-1/p)} \|u\|_{X^{1/2}_{p,\Delta}}.       \label{6linear38}
\end{align}
Choosing $0<\varepsilon <1/2p$, and making the summation on $j_{\min}$, $j_{\rm med}$ and $j_{\max}$ in order, one obtain that
\begin{align}
 |\mathscr{L}_{hha_+a_+a_+ }(w)|
 \lesssim & T^{4\varepsilon/3} \sum_{ \lambda_0} \langle \lambda_0\rangle^{ -2/ p + 2\varepsilon }   \|v_{\lambda_0}\|_{V^2_\Delta}    \|u \|^5_{X^{1/2}_{p,\Delta}}.   \label{6linear39}
\end{align}
Noticing that $\{\langle \lambda_0\rangle^{ -2/p +2\varepsilon}\} \in \ell^p$, by H\"older's inequality we have
\begin{align}
 |\mathscr{L}_{hha_+a_+a_+ }(w)|
\lesssim & T^{ \varepsilon }     \|v \|_{Y^{0}_{p',\Delta}}    \|u \|^5_{X^{1/2}_{p,\Delta}}. \ \ (0<\varepsilon \leq 1/2p)    \label{6linear40}
\end{align}

{\it Case $hha_+a_+a_-$.}   We denote by $(\lambda_k) \in hha_+a_+a_-$ the case that $\lambda_0,...,\lambda_5$ satisfy \eqref{6linear5} and
\begin{align}
 \lambda_k \in [c\lambda_0, \lambda_0], \ k=1, 3, \ \ \lambda_2, \lambda_4 \in [0, \lambda_0],  \ \lambda_5 \in [-\lambda_0, 0].     \label{6linear41}
\end{align}
We decompose $\lambda_2, \lambda_4, \lambda_5$ by:
\begin{align}
\lambda_k \in [0, \lambda_0]= \bigcup_{j_k\geq 0}  I_{j_k}, \ k=2,4;\ \   \lambda_5 \in [-\lambda_0, 0]= \bigcup_{j_5\geq 0}  -I_{j_5}, \     \ \   j_2, j_4,  j_5 \leq \log^{\lambda_0}_2.    \label{6linear42}
\end{align}
Repeating the procedures as in Case $hha_+a_+a_+$,  we can show that
\begin{align}
 |\mathscr{L}_{hha_+a_+a_-}(w)|
\lesssim & T^{ \varepsilon }     \|v \|_{Y^{0}_{p',\Delta}}    \|u \|^5_{X^{1/2}_{p,\Delta}}.   \label{6linear43}
\end{align}
Recall that in the proof above, the condition that two $\lambda_k$ is localized in higher frequency, say $\lambda_1, \lambda_3\in [c\lambda_0,\lambda_0]$ can guarantee that $\{\langle \lambda_0\rangle^{ -2/ p + 2\varepsilon }\}$ is convergent in $\ell^p$. Hence, applying the same way as in the above, we can obtain the result for the other cases which contain two higher frequency and three all frequency.

If some $\lambda_k$ is only in lower frequency $[0, c\lambda_0]$, the proof is almost the same as in the above. If all $\lambda_k$ have three or more higher frequency $[c\lambda_0, \lambda_0]$, the proof is easier than the Case $hha_+a_+a_+$, since the dyadic decomposition starting at $0$ for the third higher frequency has at most finite dyadic intervals.  $\hfill\Box$\\

In \eqref{6linear2}, for the sake of symmetry, we can assume that
$$
\lambda_1\geq \lambda_3\geq \lambda_5,  \ \ \lambda_2\geq \lambda_4 \leqno{\rm (H1)}
$$
One easily sees that $\lambda_0,..., \lambda_5$  satisfy the following frequency constrainted condition:
$$
\lambda_0 + \lambda_2 + \lambda_4 = \lambda_1 + \lambda_3 + \lambda_5 + l,  \ |l|\leq 10.  \  \ \leqno{\rm (H2)}
$$
The non-trivial case is that
$$
\lambda_0 =\max_{0\leq k\leq 5} |\lambda_k| \gg 1.  \  \ \leqno{\rm (H3)}
$$
The case $\lambda_0=\max_{0\leq k\leq 5} |\lambda_k| \lesssim 1$  implies that the summation in \eqref{6linear2} has at most finite terms.  So, in view of (H1), (H2) and (H3), we see that the orders of $\lambda_0,...,\lambda_5$ have the following 10 cases:
\begin{align*}
\lambda_0\geq \lambda_2 \geq \lambda_1 \geq \lambda_3\geq  \lambda_5 \geq \lambda_4, \tag{\rm Ord1}\\
\lambda_0\geq \lambda_1 \geq \lambda_2 \geq \lambda_3\geq  \lambda_5 \geq \lambda_4, \tag{\rm  Ord2}\\
\lambda_0\geq \lambda_1 \geq \lambda_3 \geq \lambda_5\geq  \lambda_2 \geq \lambda_4, \tag{\rm  Ord3}\\
\lambda_0\geq \lambda_1 \geq \lambda_3 \geq \lambda_2\geq  \lambda_4 \geq \lambda_5, \tag{\rm  Ord4}\\
\lambda_0\geq \lambda_1 \geq \lambda_3 \geq \lambda_2\geq  \lambda_5 \geq \lambda_4, \tag{\rm  Ord5}\\
\lambda_0\geq \lambda_1 \geq \lambda_2 \geq \lambda_3\geq  \lambda_4 \geq \lambda_5, \tag{\rm  Ord6}\\
\lambda_0\geq \lambda_1 \geq \lambda_2 \geq \lambda_4\geq  \lambda_3 \geq \lambda_5, \tag{\rm  Ord7}\\
\lambda_0\geq \lambda_2 \geq \lambda_4 \geq \lambda_1\geq  \lambda_3 \geq \lambda_5, \tag{\rm  Ord8}\\
\lambda_0\geq \lambda_2 \geq \lambda_1 \geq \lambda_3\geq  \lambda_4 \geq \lambda_5, \tag{\rm  Ord9}\\
\lambda_0\geq \lambda_2 \geq \lambda_1 \geq \lambda_4\geq  \lambda_3 \geq \lambda_5. \tag{\rm  Ord10}
\end{align*}

If there are at least two higher frequency in $\lambda_1,...,\lambda_5$ which are localized in $ [c\lambda_0, \lambda_0]\cup [-\lambda_0, -c\lambda_0]$, in the proof of Lemma \ref{nonestimate2h} we neither consider the orders  of $\lambda_1,...,\lambda_5$ nor use the constraint condition (H2).  Now we study the case that there is only one higher frequency in $\lambda_1,...,\lambda_5$. For $\lambda_1,...,\lambda_5$, if there is only one higher frequency which is localized in $ [c\lambda_0, \lambda_0]\cup [-\lambda_0, -c\lambda_0]$, we see that it must be the biggest one localized in $ [c\lambda_0, \lambda_0]$ or the smallest one localized in $ [-\lambda_0, -c\lambda_0]$, the other frequency is localized in $ [-c\lambda_0, c\lambda_0]$.

First, let us consider (Ord4) case according to the high-low frequency. The case $\lambda_5\in  h_-$ and $\lambda_1,...,\lambda_4\in l \ {\rm or } \ l_-$ never happens for small $c>0$. So, it suffices to consider the case $\lambda_1\in h$.  we divide the proof into a few cases, see Table 6.

\begin{table}[h]
\begin{center}
Case (Ord4)

\begin{tabular}{|c|c|c|c|c|c|}
\hline
 ${\rm Case}$   &  $\lambda_1\in $ & $\lambda_3\in $ & $\lambda_2\in $ &  $\lambda_4\in $ & $\lambda_5\in$  \\
\hline
$hllll$  & $h$  & $l$  & $l$  & $l$   & $l$ \\
\hline
$hllll_-$ & $h$  & $l$ & $l$  & $l$   &  $l_-$  \\
\hline
$hlll_-l_-$ & $h$  & $l$  & $l$  & $l_- $   & $l_-$ \\
\hline
$hll_-l_-l_-$ & $ h$  & $l$  & $l_-$  & $l_- $   & $l_- $ \\
\hline
$hl_-l_-l_-l_-$ & $ h$  & $l_-$  & $l_-$  & $l_-$   & $l_-$ \\
\hline
\end{tabular}
\end{center}
\caption{$\lambda_0\geq\lambda_1\geq\lambda_3\geq \lambda_2 \geq \lambda_2\geq \lambda_5$,  only one higher frequency in $\lambda_1,...,\lambda_5$.}
\end{table}

{\it Case hllll.}  We denote by $(\lambda_k) \in hllll$ the case that  $\lambda_0,..., \lambda_5$ satisfy \eqref{6linear5}, (Ord4) and
\begin{align}
 \lambda_1 \in h,  \ \ \lambda_k \in l, \ k=2,3,4,5.     \label{6linear44}
\end{align}
We decompose $\lambda_1,..., \lambda_5$ in a dyadic way:
\begin{align}
\lambda_1 \in [c\lambda_0,\lambda_0] = \bigcup_{j_1\geq 0}(\lambda_0-I_{j_1}), \ \  \lambda_k \in [0, \lambda_0]= \bigcup_{j_k\geq 0}  I_{j_k}, \ \  k=2,...,5.         \label{6linear45}
\end{align}
By condition (Ord4) we see that
\begin{align}
  0\leq j_5\leq  j_4 \leq  j_2 \leq j_3 \leq \log^{\lambda_0}_2.    \label{6linear46}
\end{align}
One has that
\begin{align}
 \mathscr{L}_{hllll}(w) \leq \!\!\!\!\!\!\!\! \sum_{\lambda_0; \ j_1; \  j_5\leq  j_4 \leq  j_2 \leq j_3 \leq \log^{\lambda_0}_2 } \!\!\!\!\!\!\!\! \langle \lambda_0\rangle^{1/2} \int_{\mathbb{R}\times [0,T]} |\overline{v}_{\lambda_0} u_{\lambda_0-I_{j_1}}  \prod^5_{k=2} w^{(k)}_{I_{j_k}} (x,t)| dxdt, \label{6linear47}
\end{align}
Using the frequency constraint condition (H2), we see that $j_1> j_3+10$ implies that $\mathscr{L}_{hllll}(w) =0$. So, we have $j_1\leq j_3+10$ in  \eqref{6linear47}.  By H\"older's inequality, we have
\begin{align}
 |\mathscr{L}_{hllll}(w)| \lesssim & \!\!\!\! \sum_{\lambda_0; \ j_1; \  j_5\leq  j_4 \leq  j_2 \leq j_3 \leq \log^{\lambda_0}_2 }  \!\!\!\! \langle \lambda_0\rangle^{1/2} \|\overline{v}_{\lambda_0}\|_{L^\infty_{x,t}} \|u_{\lambda_0-I_{j_1}}\|_{L^4_{x,t\in [0,T]}}  \nonumber\\
  & \ \ \times   \|u_{I_{j_3}}\|_{L^4_{x,t\in [0,T]}}\|u_{I_{j_2}}\|_{L^4_{x,t\in [0,T]}}\|u_{I_{j_4}}\|_{L^4_{x,t\in [0,T]}} \|u_{I_{j_5}}\|_{L^\infty_{x,t}}.   \label{6linear48}
\end{align}
By Corollaries \ref{L4toXj} and \ref{L4toXjlambda}, we have
\begin{align}
 |\mathscr{L}_{hllll}(w)| \lesssim & \  T^{4\varepsilon/3} \!\!\! \!\!\!\!\!\! \sum_{\lambda_0; \ j_1\leq j_3+10; \  j_5\leq  j_4 \leq  j_2 \leq j_3 }  \!\!\! \!\!\!\!\!\! \langle\lambda_0\rangle^{1/2} \|v_{\lambda_0}\|_{V^2_\Delta} \langle\lambda_0\rangle^{-1/2}  2^{j_1(1/4-1/p+\varepsilon)} \|u\|_{ X^{1/2}_{p,\Delta} (\lambda_0-I_{j_1}) }\nonumber\\  & \ \ \times  2^{(j_2+j_3+j_4)(-1/4-1/p+\varepsilon)} 2^{j_5(1/2-1/p)}
  \|u \|^4_{X^{1/2}_{p,\Delta}  } .   \label{6linear49}
\end{align}
Taking $0<\varepsilon \leq 1/4p$ and summarizing over all $j_5, j_4, j_2$ and $j_3$ in order, we obtain that
\begin{align}
 |\mathscr{L}_{hllll}(w)| \lesssim & \  T^{ \varepsilon }   \sum_{\lambda_0; \ j_1\leq j_3+10}   \!\!   \|v_{\lambda_0}\|_{V^2_\Delta} 2^{j_1(1/4-1/p+\varepsilon)} 2^{ j_3(-1/4-1/p+\varepsilon)} \|u\|_{ X^{1/2}_{p,\Delta} (\lambda_0-I_{j_1}) }
  \|u \|^4_{X^{1/2}_{p,\Delta}  } \nonumber\\
 \lesssim & \  T^{ \varepsilon }   \sum_{\lambda_0; \ j_1 }   \!\!   \|v_{\lambda_0}\|_{V^2_\Delta} 2^{j_1( -2/p+2\varepsilon)}   \|u\|_{ X^{1/2}_{p,\Delta} (\lambda_0-I_{j_1}) }
  \|u \|^4_{X^{1/2}_{p,\Delta}}.    \label{6linear50}
\end{align}
Noticing that $0<\varepsilon \leq 1/4p$, by H\"older's inequality, one has that
\begin{align}
 |\mathscr{L}_{hllll}(w)| \lesssim &
 \  T^{ \varepsilon }   \|v \|_{Y^{0}_{p',\Delta}} \sum_{ j_1 }      2^{-3j_1 /2p}  \left( \sum_{\lambda_0, \ \lambda\in \lambda_0-I_{j_1}} \langle \lambda\rangle^{p/2}  \|u_\lambda\|^p_{V^2_\Delta}\right)^{1/p}
  \|u \|^4_{X^{1/2}_{p,\Delta}} \nonumber\\
\lesssim &
 \  T^{ \varepsilon }   \|v \|_{Y^{0}_{p',\Delta}}
  \|u \|^5_{X^{1/2}_{p,\Delta}} \sum_{ j_1 }      2^{- j_1 /2p}  \nonumber\\
\lesssim &
 \  T^{ \varepsilon }   \|v \|_{Y^{0}_{p',\Delta}} \|u \|^5_{X^{1/2}_{p,\Delta}}.            \label{6linear51}
\end{align}

{\it Case $hllll_-$.}  We denote by $(\lambda_k) \in hllll_-$ the case that  $\lambda_0,..., \lambda_5$ satisfy \eqref{6linear5}, (Ord4) and
\begin{align}
 \lambda_1 \in h;  \ \ \lambda_k \in l, \ k=2,3,4;  \ \ \lambda_5\in l_-.    \label{6linear54}
\end{align}
We decompose $\lambda_1,..., \lambda_5$ in a dyadic way:
\begin{align}
& \lambda_1 \in [c\lambda_0,\lambda_0] = \bigcup_{j_1\geq 0}(\lambda_0-I_{j_1}); \ \  \lambda_k \in [0, \lambda_0]= \bigcup_{j_k\geq 0}  I_{j_k}, \ \  k=2,3,4;\nonumber\\
& \ \ \lambda_5 \in [-c\lambda_0,0]= \bigcup_{j_5\geq 0}  -I_{j_5}.         \label{6linear55}
\end{align}
By condition (Ord4) we see that
\begin{align}
  0\leq j_4 \leq  j_2 \leq j_3 \leq \log^{\lambda_0}_2-C.    \label{6linear56}
\end{align}
One has that
\begin{align}
 \mathscr{L}_{hllll_-}(w) \leq   \!\!\!  \!\!\!  \!\!\! \sum_{\lambda_0; \ j_1; \  j_5; \  j_4 \leq  j_2 \leq j_3 \leq \log^{\lambda_0}_2}  \!\!\!  \!\!\!   \langle \lambda_0\rangle^{1/2} \int_{\mathbb{R}\times [0,T]} |\overline{v}_{\lambda_0} u_{\lambda_0-I_{j_1}} u_{-I_{j_5}} \prod^4_{k=2} w^{(k)}_{I_{j_k}} (x,t)| dxdt. \label{6linear57}
\end{align}
Using the frequency constraint condition (H2), we see that
\begin{align}
j_1 \leq  j_3 + C.   \label{6linear58}
\end{align}
Denote
$$j_{\max} = j_2 \vee j_5, \ \ j_{\min} = j_4 \wedge j_5, \ \ j_{\rm med} \in \{j_2,j_4,j_5\} \setminus \{j_{\max}, j_{\min}\}.$$
Following the ideas as in the estimate of Case $hllll$, by H\"older's inequality and Corollaries \ref{L4toXj}, \ref{L4toXjlambda}, we have
\begin{align}
 |\mathscr{L}_{hllll_-}(w)| \lesssim & \  T^{\varepsilon} \!\!\!\! \sum_{\lambda_0; \ j_1\leq j_3+C, \  j_{\max}, j_{\rm med}, j_{\min}}  \!\!\!\! \langle\lambda_0\rangle^{1/2} \|v_{\lambda_0}\|_{V^2_\Delta} \langle\lambda_0\rangle^{-1/2}  2^{j_1(1/4-1/p+\varepsilon)} \|u\|_{ X^{1/2}_{p,\Delta} (\lambda_0-I_{j_1}) }\nonumber\\  & \ \ \times  2^{(j_3+j_{\max}+j_{\rm med})(-1/4-1/p+\varepsilon)} 2^{j_{\min}(1/2-1/p)}
  \|u \|^4_{X^{1/2}_{p,\Delta}  }.     \label{6linear59}
\end{align}
Summarizing over all $j_{\min}, j_{\rm med}, j_{\max}$ in orders, from \eqref{6linear59} one can control $|\mathscr{L}_{hllll_-}(w)|$ by the right hand side of \eqref{6linear50}.  So, we have
\begin{align}
 |\mathscr{L}_{hllll_-}(w)| \lesssim
 \  T^{ \varepsilon }   \|v \|_{Y^{0}_{p',\Delta}} \|u \|^5_{X^{1/2}_{p,\Delta}}.            \label{6linear60}
\end{align}

{\it Case $hlll_-l_-$.}  We denote by $(\lambda_k) \in hlll_-l_-$ the case that  $\lambda_0,..., \lambda_5$ satisfy \eqref{6linear5}, (Ord4) and
\begin{align}
 \lambda_1 \in h;  \ \ \lambda_2, \lambda_3 \in l;  \ \ \lambda_4, \lambda_5\in l_-.    \label{6linear61}
\end{align}
We decompose $\lambda_1,..., \lambda_5$ by
\begin{align}
& \lambda_1 \in [c\lambda_0,\lambda_0] = \bigcup_{j_1\geq 0}(\lambda_0-I_{j_1}); \ \  \lambda_k \in [0, \lambda_0]= \bigcup_{j_k\geq 0}  I_{j_k}, \ \  k=2,3;\nonumber\\
& \ \ \lambda_k \in [-c\lambda_0,0]= \bigcup_{j_k\geq 0}  -I_{j_k}, \ \ k=4,5.         \label{6linear65}
\end{align}
By condition (Ord4) we see that
\begin{align}
  0\leq  j_2 \leq j_3 \leq \log^{\lambda_0}_2-C, \ \  0\leq  j_4 \leq j_5\leq \log^{\lambda_0}_2-C.    \label{6linear66}
\end{align}
One has that 
\begin{align}
 \mathscr{L}_{hlll_-l_-}(w) \leq   \!\!\!  \!\!\!  \!\!\! \sum_{\lambda_0; \ j_1; \  j_4 \leq j_5; \ j_2 \leq j_3 \leq \log^{\lambda_0}_2}  \!\!\!  \!\!\!   \langle \lambda_0\rangle^{1/2} \int_{\mathbb{R}\times [0,T]} |\overline{v}_{\lambda_0} u_{\lambda_0-I_{j_1}} \overline{u}_{I_{j_2}}  u_{I_{j_3}} \overline{u}_{-I_{j_4}} u_{-I_{j_5}} (x,t)| dxdt. \label{6linear67}
\end{align}
In view if the frequency constraint condition (H2), we see that
$$
I_{j_1} \cap (I_{j_3}+ I_{j_4}- I_{j_2}- I_{j_5}) \neq \varnothing.
$$
Otherwise  $\mathscr{L}_{hlll_-l_-}(w)=0$. It follows that

$(1)$ If $j_4<j_5$, then $j_1< j_3+10$.

$(2)$ If $j_4= j_5$, then $j_1< j_3\vee j_4+10$.

In the case $j_4<j_5$, we can use a similar way as in Case $hllll_-$. Indeed,  putting
$$j_{\max} = j_2 \vee j_5, \ \ j_{\min} = j_4 \wedge j_2, \ \ j_{\rm med} \in \{j_2,j_4,j_5\} \setminus \{j_{\max}, j_{\min}\},$$
then we can repeat the procedure as in Case $hllll_-$ to have the estimate
\begin{align}
 |\mathscr{L}_{hlll_-l_-}(w) |  \lesssim T^\varepsilon \|v\|_{Y^{0}_{p',\Delta}} \|u\|^5_{X^{1/2}_{p,\Delta}}.   \label{6linear68}
\end{align}

In the case $j_4=j_5$, we need to separate the proof into the case $j_3\geq j_4$ and $j_3< j_4$, respectively.  We have

$(2a)$ If $j_3\leq j_4 $, then $ j_2\leq j_3\leq j_4=j_5$.

$(2b)$ If $j_3\geq j_4 $, then $ j_3\geq j_4=j_5 \geq j_2, $ or $j_3\geq j_2 \geq  j_4=j_5 $.

In the case $(2b)$, we can use the same way as in Case $hllll_-$ to show that \eqref{6linear68} holds. In the case $(2a)$, we have from Corollaries \ref{L4toXj} and \ref{L4toXjlambda} that
\begin{align}
 |\mathscr{L}_{hlll_-l_-}(w)| \lesssim & \sum_{\lambda_0; \ j_1\leq j_4+C; \ j_2 \leq j_3 \leq  j_5=  j_4} \langle \lambda_0\rangle^{1/2} \|\overline{v}_{\lambda_0}\|_{L^\infty_{x,t}} \|u_{\lambda_0-I_{j_1}}\|_{L^4_{x,t\in [0,T]}}  \nonumber\\
  & \ \ \times   \|u_{I_{j_3}}\|_{L^4_{x,t\in [0,T]}}\|u_{-I_{j_4}}\|_{L^4_{x,t\in [0,T]}} \|u_{-I_{j_5}}\|_{ L^4_{x,t\in [0,T]}}\|u_{I_{j_2}}\|_{L^\infty_{x,t}} .   \nonumber\\
 \lesssim & \  T^{4\varepsilon/3} \!\!\! \!\!\!\!\!\! \sum_{\lambda_0; \ j_1\leq j_4+C; \ j_2 \leq j_3 \leq  j_5=  j_4 }  \!\!\! \!\!\!\!\!\! \langle\lambda_0\rangle^{1/2} \|v_{\lambda_0}\|_{V^2_\Delta} \langle\lambda_0\rangle^{-1/2}  2^{j_1(1/4-1/p+\varepsilon)} \|u\|_{ X^{1/2}_{p,\Delta} (\lambda_0-I_{j_1}) }\nonumber\\  & \ \ \times  2^{(j_5+j_3+j_4)(-1/4-1/p+\varepsilon)} 2^{j_2(1/2-1/p)}
  \|u \|^4_{X^{1/2}_{p,\Delta}  } .   \label{6linear69}
\end{align}
Making the summation in turn to $j_2, j_3$ and $j_4$, we can get \eqref{6linear68}.

{\it Case $hll_-l_-l_-$.}  We denote by $(\lambda_k) \in hll_-l_-l_-$ the case that  $\lambda_0,..., \lambda_5$ satisfy \eqref{6linear5}, (Ord4) and
\begin{align}
 \lambda_1 \in h;  \ \ \lambda_3 \in l;  \ \ \lambda_2, \lambda_4, \lambda_5\in l_-.    \label{6linear71}
\end{align}
We decompose $\lambda_1,..., \lambda_5$ by
\begin{align}
& \lambda_1 \in [c\lambda_0,\lambda_0] = \bigcup_{j_1\geq 0}(\lambda_0-I_{j_1}); \ \  \lambda_3 \in [0, \lambda_0]= \bigcup_{j_3\geq 0}  I_{j_3};\nonumber\\
& \ \ \lambda_k \in [-c\lambda_0,0]= \bigcup_{j_k\geq 0}  -I_{j_k}, \ \ k=2, 4,5.         \label{6linear75}
\end{align}
By condition (Ord4) we see that
\begin{align}
   0 \leq j_2 \leq j_4 \leq j_5\leq \log^{\lambda_0}_2-C.    \label{6linear76}
\end{align}
Using the frequency constraint condition (H2), we see that

\mbox{\rm (i) } If $ j_3\geq j_5$, then $j_1\leq j_3+10$.

(ii) If $j_3 < j_5$, then $j_1< j_4+10$.

In the case (i), we can follow the same way as in Case $hllll$ to obtain the estimate
\begin{align}
 |\mathscr{L}_{hll_-l_-l_-}(w)| \lesssim  T^\varepsilon \|v\|_{Y^{0}_{p',\Delta}}
  \|u \|^4_{X^{1/2}_{p,\Delta}  } .   \label{6linear79}
\end{align}
In the case (ii), we need to further analyze $j_2, j_3, j_4$ and denote
$$
j_{\max} = j_3 \vee j_4, \ \ j_{\min} = j_2 \wedge j_3, \ \ j_{\rm med} \in \{j_2,j_3,j_4\} \setminus \{j_{\max}, j_{\min}\}.
$$
Then we can repeat the procedures as in Case $hllll_-$ to obtain \eqref{6linear79}.

{\it Case $hl_-l_-l_-l_-$.}  We denote by $(\lambda_k) \in hl_-l_-l_-l_-$ the case that  $\lambda_0,..., \lambda_5$ satisfy \eqref{6linear5}, (Ord4) and
\begin{align}
 \lambda_1 \in h;     \ \ \lambda_2, ..., \lambda_5\in l_-.    \label{6linear81}
\end{align}
We decompose $\lambda_1,..., \lambda_5$ by
\begin{align}
& \lambda_1 \in [c\lambda_0,\lambda_0] = \bigcup_{j_1\geq 0}(\lambda_0-I_{j_1});
& \ \ \lambda_k \in [-c\lambda_0,0]= \bigcup_{j_k\geq 0}  -I_{j_k}, \   k=2,...,5.         \label{6linear85}
\end{align}
By condition (Ord4) we see that
\begin{align}
   0\leq j_3 \leq j_2 \leq j_4 \leq j_5\leq \log^{\lambda_0}_2-C.    \label{6linear86}
\end{align}
Using the frequency constraint condition (H2), we see that
$$ j_5 \leq j_4 +2, \ \ j_1\leq j_4+10.$$
Noticing that $j_3 \leq j_2 \leq j_4 \leq j_5 \leq j_4+2$, using a similar way as in Case $hllll$ we can obtain the estimate
 \begin{align}
 |\mathscr{L}_{hl_-l_-l_-l_-}(w)| \lesssim  T^\varepsilon \|v\|_{Y^{0}_{p',\Delta}}
  \|u \|^4_{X^{1/2}_{p,\Delta}  } .   \label{6linear89}
\end{align}

Next, we consider (Ord1) case according to the high-low frequency. The case $\lambda_2\in  h $ and $\lambda_1,\lambda_3,\lambda_4,\lambda_5\in l \ {\rm or } \ l_-$ never happens for small $c>0$. So, it suffices to consider the case $\lambda_4\in h_-$.  We divide the proof into a few cases, see Table 7.

\begin{table}[h]

\begin{center}
Case (Ord1)

\begin{tabular}{|c|c|c|c|c|c| }
\hline
 ${\rm Case}$   &  $\lambda_2\in $ & $\lambda_1\in $ & $\lambda_3\in $ &  $\lambda_5\in $ & $\lambda_4\in$   \\
\hline
$llllh_-$  & $l$  & $l$  & $l$  & $l$   & $h_-$   \\
\hline
$llll_-h_-$ & $l$  & $l$ & $l$  & $l_-$   &  $h_-$    \\
\hline
$hlll_-l_-h_-$ & $l$  & $l$  & $l_-$  & $l_- $   & $h_-$    \\
\hline
$ll_-l_-l_-h_-$ & $ l$  & $l_-$  & $l_-$  & $l_- $   & $h_- $   \\
\hline
$l_-l_-l_-l_-h_-$ & $ l_-$  & $l_-$  & $l_-$  & $l_-$   & $h_-$    \\
\hline
\end{tabular}
\end{center}
\caption{$\lambda_0\geq\lambda_2\geq\lambda_1\geq \lambda_3 \geq \lambda_5\geq \lambda_4$,  only one higher frequency in $\lambda_1,...,\lambda_5$.}
\end{table}

We consider the decomposition of $\lambda_k$:
$$
\lambda_4 \in [-\lambda_0, -c\lambda_0] = \bigcup_{j_4\geq 0}(-\lambda_0 + I_{j_4}), \ \  \lambda_k \in [0,c\lambda_0] = \bigcup_{j_k\geq 0} I_{j_k};  \ \  \lambda_l \in [-c\lambda_0,0] = \bigcup_{j_l\geq 0}- I_{j_l}
$$
for any $k,l\in \{1,2,3,5\}$. For short, we denote by $(\lambda_k) \in llllh_-$ the case that $\lambda_0,...,\lambda_5$ satisfy \eqref{6linear5}, (Ord1) and $\lambda_1,...,\lambda_3,\lambda_5\in l$, $\lambda_4\in h_-$.  Similarly, we will use the notations  $(\lambda_k) \in llll_-h_-,..., l_-l_-l_-l_-h_-$.

 Using the frequency constraint condition (H2), we see that
\begin{align*}
& j_4\leq j_1+10 \ \ {\rm if} \ \ (\lambda_k) \in  llllh_- \cup llll_-h_- \cup    lll_-l_-h_-. \\
& j_1\vee ...\vee j_5 \lesssim 1, \ \ {\rm if} \ \ (\lambda_k) \in ll_-l_-l_-h_- \cup l_-l_-l_-l_-h_- .
\end{align*}
The Cases $ll_-l_-l_-h_-$ and $l_-l_-l_-l_-h_-$ are trivial and we omit the discussions.  For the Cases $llllh_-$, $llll_-h_-$  and $lll_-l_-h_-$, we need to further analyze $j_2, j_3, j_5$ and write
$$
j_{\max} = j_2\vee j_3 \vee j_5, \ \ j_{\min} = j_2 \wedge j_3 \wedge j_5, \ \ j_{\rm med} \in \{j_2,j_3,j_5\} \setminus \{j_{\max}, j_{\min}\}.
$$
Similar to  \eqref{6linear59}, we have

\begin{align}
 |\mathscr{L}_{llllh_-}& (w)| + |\mathscr{L}_{llll_-h_-}(w)| + |\mathscr{L}_{lll_-l_-h_-}(w)| \nonumber\\
 \lesssim & \  T^{\varepsilon}  \sum_{\lambda_0; \ j_4\leq j_1+C, \  j_{\max}, j_{\rm med}, j_{\min}}  \langle\lambda_0\rangle^{1/2} \|v_{\lambda_0}\|_{V^2_\Delta} \langle\lambda_0\rangle^{-1/2}  2^{j_4(1/4-1/p+\varepsilon)} \|u\|_{ X^{1/2}_{p,\Delta} (-\lambda_0+I_{j_4}) }\nonumber\\  & \ \ \times  2^{(j_1+j_{\max}+j_{\rm med})(-1/4-1/p+\varepsilon)} 2^{j_{\min}(1/2-1/p)}
  \|u \|^4_{X^{1/2}_{p,\Delta}  }.     \label{6linear89}
\end{align}
Using the same way as above and making the summation in turn to  $j_{\min}, j_{\rm med},j_{\max}$ and then to $\lambda_0$ and $j_4$, we have
\begin{align}
 |\mathscr{L}_{llllh_-}  (w)| + |\mathscr{L}_{llll_-h_-}(w)| + |\mathscr{L}_{lll_-l_-h_-}(w)|
 \lesssim  \  T^{\varepsilon}   \|v\|_{Y^{0}_{p',\Delta}}  \|u \|^4_{X^{1/2}_{p,\Delta}  }.     \label{6linear90}
\end{align}

If $\lambda_0,..,\lambda_5$ has the orders as in (Ord2) and only one of them is lying in higher frequency, there are two possible cases $\lambda_1\in h$ and $\lambda_4\in h_-$. We discuss these two cases separately.

\begin{table}[h]

\begin{center}
Case (Ord2): $\lambda_1\in h$

\begin{tabular}{|c|c|c|c|c|c|}
\hline
 ${\rm Case}$   &  $\lambda_1\in $ & $\lambda_2\in $ & $\lambda_3\in $ &  $\lambda_5\in $ & $\lambda_4\in$  \\
\hline
$hllll$  & $h$  & $l$  & $l$  & $l$   & $l$ \\
\hline
$hllll_-$ & $h$  & $l$ & $l$  & $l$   &  $l_-$  \\
\hline
$hlll_-l_-$ & $h$  & $l$  & $l$  & $l_- $   & $l_-$ \\
\hline
$hll_-l_-l_-$ & $ h$  & $l$  & $l_-$  & $l_- $   & $l_- $ \\
\hline
$hl_-l_-l_-l_-$ & $ h$  & $l_-$  & $l_-$  & $l_-$   & $l_-$ \\
\hline
\end{tabular}
\end{center}
\caption{$\lambda_0\geq\lambda_1\geq\lambda_2\geq \lambda_3 \geq \lambda_5\geq \lambda_4$,  only one higher frequency in $\lambda_1,...,\lambda_5$.}
\end{table}
As in Table 8, we denote by $(\lambda_k) \in hllll$ that $\lambda_0,..,\lambda_5$ satisfy \eqref{6linear5}, (Ord2) and $\lambda_1\in h$, $\lambda_2,...,\lambda_5\in l$.  Similarly for the notations  $(\lambda_k) \in hllll_-$,...,$hl_-l_-l_-l_-$.  We consider the decomposition of $\lambda_k$:
$$
\lambda_1 \in [c\lambda_0,  \lambda_0] = \bigcup_{j_4\geq 0}( \lambda_0 - I_{j_1}); \ \  \lambda_k \in [0,c\lambda_0] = \bigcup_{j_k\geq 0} I_{j_k},  \ \  \lambda_l \in [-c\lambda_0,0] = \bigcup_{j_l\geq 0}- I_{j_l}
$$
for any $k,l\in \{2,3,4,5\}$. We have the following constraint conditions on $j_1,...,j_5$:
\begin{align*}
& j_1\leq j_3+C, \ \ j_2\geq j_3\geq j_5 \geq j_4  \ \ {\rm if} \ \ (\lambda_k) \in  hllll,   \\
& j_1\leq j_3 \vee j_4 +C, \ \ j_2\geq j_3\geq j_5   \ \ {\rm if} \ \ (\lambda_k) \in  hllll_-,   \\
& j_1\leq  j_3 \vee j_4+C, \ \ j_2\geq j_3, j_4\geq j_5  \ \ {\rm if} \ \ (\lambda_k) \in  hlll_-l_-,   \\
& j_1\leq     j_4, \ \  j_3\leq   j_5\leq j_4  \ \ {\rm if} \ \ (\lambda_k) \in  hll_-l_-l_-,   \\
& j_1\leq     j_4+C, \ \  j_2\leq j_3\leq   j_5\leq j_4  \ \ {\rm if} \ \ (\lambda_k) \in  hl_-l_-l_-l_-.
 \end{align*}
Using the above constraint conditions, we can follow the same ideas as in (Ord1) case to have the estimate
\begin{align}
 |\mathscr{L}_{bcdef}  (w)|
 \lesssim  \  T^{\varepsilon}   \|v\|_{Y^{0}_{p',\Delta}}  \|u \|^4_{X^{1/2}_{p,\Delta}  },     \label{6linear100}
\end{align}
where $bcdef \in \{hllll, hllll_-, hlll_-l_-, hll_-l_-l_-, hl_-l_-l_-l_- \}$.

\begin{table}[h]

\begin{center}
Case (Ord2): $\lambda_4\in h_-$

\begin{tabular}{|c|c|c|c|c|c| }
\hline
 ${\rm Case}$   &  $\lambda_1\in $ & $\lambda_2\in $ & $\lambda_3\in $ &  $\lambda_5\in $ & $\lambda_4\in$   \\
\hline
$llllh_-$  & $l$  & $l$  & $l$  & $l$   & $h_-$   \\
\hline
$llll_-h_-$ & $l$  & $l$ & $l$  & $l_-$   &  $h_-$    \\
\hline
$hlll_-l_-h_-$ & $l$  & $l$  & $l_-$  & $l_- $   & $h_-$    \\
\hline
$ll_-l_-l_-h_-$ & $ l$  & $l_-$  & $l_-$  & $l_- $   & $h_- $   \\
\hline
$l_-l_-l_-l_-h_-$ & $ l_-$  & $l_-$  & $l_-$  & $l_-$   & $h_-$    \\
\hline
\end{tabular}
\end{center}
\caption{$\lambda_0\geq\lambda_1\geq\lambda_2\geq \lambda_3 \geq \lambda_5\geq \lambda_4$,  only one higher frequency in $\lambda_1,...,\lambda_5$.}
\end{table}

As in Table 9, we denote by $(\lambda_k) \in  llllh_-$ that $\lambda_0,..,\lambda_5$ satisfy \eqref{6linear5}, (Ord2) and $\lambda_1\in h$, $\lambda_2,...,\lambda_5\in l$.  Similarly for the notations  $(\lambda_k) \in llll_-h_-$,...,$l_-l_-l_-l_-h_-$.  We consider the decomposition of $\lambda_k$:
$$
\lambda_4 \in [-\lambda_0,  -c\lambda_0] = \bigcup_{j_4\geq 0}( -\lambda_0 + I_{j_4}); \ \  \lambda_k \in [0,c\lambda_0] = \bigcup_{j_k\geq 0} I_{j_k},  \ \  \lambda_l \in [-c\lambda_0,0] = \bigcup_{j_l\geq 0}- I_{j_l}
$$
for any $k,l\in \{1,2,3,5\}$. We have the following constraint conditions on $j_1,...,j_5$:
\begin{align*}
& j_4\leq j_1+C, \ j_1\geq j_2\geq j_3\geq j_5  \ \ {\rm if} \ \ (\lambda_k) \in llll h_-,   \\
& j_4\leq j_1  +C, \ j_1\geq j_2\geq j_3, j_5\leq j_1+C   \ \ {\rm if} \ \ (\lambda_k) \in llll_-h_-,   \\
& j_4\leq  j_1 +C, \ j_1\geq j_2, \  j_3 \leq j_5 \leq j_1+C  \ \ {\rm if} \ \ (\lambda_k) \in lll_-l_-h_-,   \\
& j_4\leq  j_1 +C, \ j_2\leq j_3\leq   j_5 \leq j_1+C   \ \ {\rm if} \ \ (\lambda_k) \in ll_-l_-l_-h_-,   \\
& j_4\leq  C, \ j_1\leq j_2\leq j_3 \leq   j_5 \ \ {\rm if} \ \ (\lambda_k) \in l_-l_-l_-l_-h_-.
 \end{align*}
Using the above constraint conditions, we can follow the same ideas as in (Ord4) case to have the estimate
\begin{align}
 |\mathscr{L}_{bcdef}  (w)|
 \lesssim  \  T^{\varepsilon}   \|v\|_{Y^{0}_{p',\Delta}}  \|u \|^4_{X^{1/2}_{p,\Delta}  },     \label{6linear101}
\end{align}
where $bcdef \in \{llllh_-, llll_-h_-, lll_-l_-h_-, ll_-l_-l_-h_-, l_-l_-l_-l_-h_- \}$ as in Table 9.

If $\lambda_0,..,\lambda_5$ has the orders as in (Ord3) and only one of them is  in higher frequency, then there are two possible cases $\lambda_1\in h$ and $\lambda_4\in h_-$, see Tables 10 and 11.

\begin{table}[h]

\begin{center}
Case (Ord3): $\lambda_1\in h$

\begin{tabular}{|c|c|c|c|c|c|}
\hline
 ${\rm Case}$   &  $\lambda_1\in $ & $\lambda_3\in $ & $\lambda_5\in $ &  $\lambda_2\in $ & $\lambda_4\in$  \\
\hline
$hllll$  & $h$  & $l$  & $l$  & $l$   & $l$ \\
\hline
$hllll_-$ & $h$  & $l$ & $l$  & $l$   &  $l_-$  \\
\hline
$hlll_-l_-$ & $h$  & $l$  & $l$  & $l_- $   & $l_-$ \\
\hline
$hll_-l_-l_-$ & $ h$  & $l$  & $l_-$  & $l_- $   & $l_- $ \\
\hline
$hl_-l_-l_-l_-$ & $ h$  & $l_-$  & $l_-$  & $l_-$   & $l_-$ \\
\hline
\end{tabular}
\end{center}
\caption{$\lambda_0\geq\lambda_1\geq\lambda_3\geq \lambda_5 \geq \lambda_2\geq \lambda_4$, only one higher frequency in $\lambda_1,...,\lambda_5$.}
\end{table}
As in Table 10, we denote by $(\lambda_k) \in hllll$ that $\lambda_0,..,\lambda_5$ satisfy \eqref{6linear5}, (Ord3) and $\lambda_1\in h$, $\lambda_2,...,\lambda_5\in l$.  Similarly for the notations  $(\lambda_k) \in hllll_-$,...,$hl_-l_-l_-l_-$.  We consider the decomposition of $\lambda_k$:
$$
\lambda_1 \in [c\lambda_0,  \lambda_0] = \bigcup_{j_4\geq 0}( \lambda_0 - I_{j_1}); \ \  \lambda_k \in [0,c\lambda_0] = \bigcup_{j_k\geq 0} I_{j_k},  \ \  \lambda_l \in [-c\lambda_0,0] = \bigcup_{j_l\geq 0}- I_{j_l}
$$
for any $k,l\in \{2,3,4,5\}$. We have the following constraint conditions on $j_1,...,j_5$:
\begin{align*}
& j_1\leq j_3+C, \ \ j_3\geq j_5\geq j_2 \geq j_4  \ \ {\rm if} \ \ (\lambda_k) \in  hllll,   \\
& j_1\leq j_3 \vee j_4 +C, \ \ j_3\geq j_5\geq j_2   \ \ {\rm if} \ \ (\lambda_k) \in  hllll_-,   \\
& j_1\leq  j_3 \vee j_4+C, \ \ j_3\geq j_5, j_2\leq j_4  \ \ {\rm if} \ \ (\lambda_k) \in  hlll_-l_-,   \\
& j_1\leq   j_3 \vee  j_4, \ \  j_5\leq   j_2\leq j_4  \ \ {\rm if} \ \ (\lambda_k) \in  hll_-l_-l_-,   \\
& j_1\leq     j_4+C, \ \  j_3\leq j_5\leq   j_2\leq j_4  \ \ {\rm if} \ \ (\lambda_k) \in  hl_-l_-l_-l_-.
 \end{align*}
Using the above constraint conditions, we can follow the same ideas as in (Ord1) case to have the estimate
\begin{align}
 |\mathscr{L}_{bcdef}  (w)|
 \lesssim  \  T^{\varepsilon}   \|v\|_{Y^{0}_{p',\Delta}}  \|u \|^4_{X^{1/2}_{p,\Delta}  },     \label{6linear102}
\end{align}
where $bcdef \in \{hllll, hllll_-, hlll_-l_-, hll_-l_-l_-, hl_-l_-l_-l_- \}$ in Table 10.

\begin{table}[h]

\begin{center}
Case (Ord3): $\lambda_4\in h_-$

\begin{tabular}{|c|c|c|c|c|c| }
\hline
 ${\rm Case}$   &  $\lambda_1\in $ & $\lambda_3\in $ & $\lambda_5\in $ &  $\lambda_2\in $ & $\lambda_4\in$   \\
\hline
$llllh_-$  & $l$  & $l$  & $l$  & $l$   & $h_-$   \\
\hline
$llll_-h_-$ & $l$  & $l$ & $l$  & $l_-$   &  $h_-$    \\
\hline
$hlll_-l_-h_-$ & $l$  & $l$  & $l_-$  & $l_- $   & $h_-$    \\
\hline
$ll_-l_-l_-h_-$ & $ l$  & $l_-$  & $l_-$  & $l_- $   & $h_- $   \\
\hline
$l_-l_-l_-l_-h_-$ & $ l_-$  & $l_-$  & $l_-$  & $l_-$   & $h_-$    \\
\hline
\end{tabular}
\end{center}
\caption{$\lambda_0\geq\lambda_1\geq\lambda_3\geq \lambda_5 \geq \lambda_2\geq \lambda_4$,  only one higher frequency in $\lambda_1,...,\lambda_5$.}
\end{table}

As in Table 11, we denote by $(\lambda_k) \in  llllh_-$ that $\lambda_0,..,\lambda_5$ satisfy \eqref{6linear5}, (Ord3) and $\lambda_1\in h$, $\lambda_2,...,\lambda_5\in l$.  Similarly for the notations  $(\lambda_k) \in llll_-h_-$,...,$l_-l_-l_-l_-h_-$.  We consider the decomposition of $\lambda_k$:
$$
\lambda_4 \in [-\lambda_0,  -c\lambda_0] = \bigcup_{j_4\geq 0}( -\lambda_0 + I_{j_4}); \ \  \lambda_k \in [0,c\lambda_0] = \bigcup_{j_k\geq 0} I_{j_k},  \ \  \lambda_l \in [-c\lambda_0,0] = \bigcup_{j_l\geq 0}- I_{j_l}
$$
for any $k,l\in \{1,2,3,5\}$. We have the following constraint conditions on $j_1,...,j_5$:
\begin{align*}
& j_4\leq j_1+C, \ j_1\geq j_3\geq j_5\geq j_2  \ \ {\rm if} \ \ (\lambda_k) \in llll h_-,   \\
& j_4\leq j_1 \vee j_2  +C, \ j_1\geq j_3\geq j_5,    \ \ {\rm if} \ \ (\lambda_k) \in llll_-h_-,   \\
& j_4\leq  j_1 \vee j_2 +C, \ j_1\geq j_3, \  j_5 \leq j_2,   \ \ {\rm if} \ \ (\lambda_k) \in lll_-l_-h_-,   \\
& j_4\leq  j_1 \vee j_2 +C, \ j_3\leq j_5\leq   j_2,   \ \ {\rm if} \ \ (\lambda_k) \in ll_-l_-l_-h_-,   \\
& j_4\leq j_2+ C, \ j_1\leq j_3\leq j_5 \leq   j_2 \ \ {\rm if} \ \ (\lambda_k) \in l_-l_-l_-l_-h_-.
 \end{align*}
Using the above constraint conditions, we can follow the same ideas as in (Ord4) case to have the estimate
\begin{align}
 |\mathscr{L}_{bcdef}  (w)|
 \lesssim  \  T^{\varepsilon}   \|v\|_{Y^{0}_{p',\Delta}}  \|u \|^4_{X^{1/2}_{p,\Delta}  },     \label{6linear103}
\end{align}
where $bcdef \in \{llllh_-, llll_-h_-, lll_-l_-h_-, ll_-l_-l_-h_-, l_-l_-l_-l_-h_- \}$ as in Table 11.

If $\lambda_0,..,\lambda_5$ has the orders as in (Ord5) and only one of them is in higher frequency, then two possible cases $\lambda_1\in h$ and $\lambda_4\in h_-$ will be happened, see Tables 12 and 13.

\begin{table}[h]

\begin{center}
Case (Ord5): $\lambda_1\in h$

\begin{tabular}{|c|c|c|c|c|c|}
\hline
 ${\rm Case}$   &  $\lambda_1\in $ & $\lambda_3\in $ & $\lambda_2\in $ &  $\lambda_5\in $ & $\lambda_4\in$  \\
\hline
$hllll$  & $h$  & $l$  & $l$  & $l$   & $l$ \\
\hline
$hllll_-$ & $h$  & $l$ & $l$  & $l$   &  $l_-$  \\
\hline
$hlll_-l_-$ & $h$  & $l$  & $l$  & $l_- $   & $l_-$ \\
\hline
$hll_-l_-l_-$ & $ h$  & $l$  & $l_-$  & $l_- $   & $l_- $ \\
\hline
$hl_-l_-l_-l_-$ & $ h$  & $l_-$  & $l_-$  & $l_-$   & $l_-$ \\
\hline
\end{tabular}
\end{center}
\caption{$\lambda_0\geq\lambda_1\geq\lambda_3\geq \lambda_2 \geq \lambda_5 \geq \lambda_4$,  only one higher frequency in $\lambda_1,...,\lambda_5$.}
\end{table}
As in Table 12, we denote by $(\lambda_k) \in hllll$ that $\lambda_0,..,\lambda_5$ satisfy \eqref{6linear5}, (Ord5) and $\lambda_1\in h$, $\lambda_2,...,\lambda_5\in l$.  Similarly for the notations  $(\lambda_k) \in hllll_-$,...,$hl_-l_-l_-l_-$.  We consider the decomposition of $\lambda_k$:
$$
\lambda_1 \in [c\lambda_0,  \lambda_0] = \bigcup_{j_4\geq 0}( \lambda_0 - I_{j_1}); \ \  \lambda_k \in [0,c\lambda_0] = \bigcup_{j_k\geq 0} I_{j_k},  \ \  \lambda_l \in [-c\lambda_0,0] = \bigcup_{j_l\geq 0}- I_{j_l}
$$
for any $k,l\in \{2,3,4,5\}$. We have the following constraint conditions on $j_1,...,j_5$:
\begin{align*}
& j_1\leq j_3+C, \ \ j_3\geq j_2\geq j_5 \geq j_4  \ \ {\rm if} \ \ (\lambda_k) \in  hllll,   \\
& j_1\leq j_3 \vee j_4 +C, \ \ j_3\geq j_2\geq j_5   \ \ {\rm if} \ \ (\lambda_k) \in  hllll_-,   \\
& j_1\leq  j_3 \vee j_4+C, \ \ j_3\geq j_2, j_5\leq j_4  \ \ {\rm if} \ \ (\lambda_k) \in  hlll_-l_-,   \\
& j_1\leq   j_3 \vee  j_4, \ \  j_2\leq   j_5\leq j_4  \ \ {\rm if} \ \ (\lambda_k) \in  hll_-l_-l_-,   \\
& j_1\leq     j_4+C, \ \  j_3\leq j_2\leq   j_5 \leq j_4  \ \ {\rm if} \ \ (\lambda_k) \in  hl_-l_-l_-l_-.
 \end{align*}
Using the above constraint conditions, we can follow the same ideas as in (Ord1) case to have the estimate
\begin{align}
 |\mathscr{L}_{bcdef}  (w)|
 \lesssim  \  T^{\varepsilon}   \|v\|_{Y^{0}_{p',\Delta}}  \|u \|^4_{X^{1/2}_{p,\Delta}  },     \label{6linear112}
\end{align}
where $bcdef \in \{hllll, hllll_-, hlll_-l_-, hll_-l_-l_-, hl_-l_-l_-l_- \}$ in Table 12.

\begin{table}[h]

\begin{center}
Case (Ord5): $\lambda_4\in h_-$

\begin{tabular}{|c|c|c|c|c|c| }
\hline
 ${\rm Case}$   &  $\lambda_1\in $ & $\lambda_3\in $ & $\lambda_2\in $ &  $\lambda_5\in $ & $\lambda_4\in$   \\
\hline
$llllh_-$  & $l$  & $l$  & $l$  & $l$   & $h_-$   \\
\hline
$llll_-h_-$ & $l$  & $l$ & $l$  & $l_-$   &  $h_-$    \\
\hline
$hlll_-l_-h_-$ & $l$  & $l$  & $l_-$  & $l_- $   & $h_-$    \\
\hline
$ll_-l_-l_-h_-$ & $ l$  & $l_-$  & $l_-$  & $l_- $   & $h_- $   \\
\hline
$l_-l_-l_-l_-h_-$ & $ l_-$  & $l_-$  & $l_-$  & $l_-$   & $h_-$    \\
\hline
\end{tabular}
\end{center}
\caption{$\lambda_0\geq\lambda_1\geq\lambda_3\geq \lambda_2 \geq \lambda_5 \geq \lambda_4$,  only one higher frequency in $\lambda_1,...,\lambda_5$.}
\end{table}

As in Table 13, we denote by $(\lambda_k) \in  llllh_-$ that $\lambda_0,..,\lambda_5$ satisfy \eqref{6linear5}, (Ord5) and $\lambda_1\in h$, $\lambda_2,...,\lambda_5\in l$.  Similarly for the notations  $(\lambda_k) \in llll_-h_-$,...,$l_-l_-l_-l_-h_-$.  We consider the decomposition of $\lambda_k$:
$$
\lambda_4 \in [-\lambda_0,  -c\lambda_0] = \bigcup_{j_4\geq 0}( -\lambda_0 + I_{j_4}); \ \  \lambda_k \in [0,c\lambda_0] = \bigcup_{j_k\geq 0} I_{j_k},  \ \  \lambda_l \in [-c\lambda_0,0] = \bigcup_{j_l\geq 0}- I_{j_l}
$$
for any $k,l\in \{1,2,3,5\}$. We have the following constraint conditions on $j_1,...,j_5$:
\begin{align*}
& j_4\leq j_1+C, \ j_1\geq j_3\geq j_2\geq j_5  \ \ {\rm if} \ \ (\lambda_k) \in llll h_-,   \\
& j_4\leq j_1   +C, \ j_1\geq j_3\geq j_2,    \ \ {\rm if} \ \ (\lambda_k) \in llll_-h_-,   \\
& j_4\leq  j_1 \vee j_2 +C, \ j_1\geq j_3, \  j_5 \geq j_2,   \ \ {\rm if} \ \ (\lambda_k) \in lll_-l_-h_-,   \\
& j_4\leq  j_1 \vee j_2 +C, \ j_3\leq j_2\leq   j_5,   \ \ {\rm if} \ \ (\lambda_k) \in ll_-l_-l_-h_-,   \\
& j_4\leq j_2+ C, \ j_1\leq j_3\leq j_2 \leq   j_5, \ \ {\rm if} \ \ (\lambda_k) \in l_-l_-l_-l_-h_-.
 \end{align*}
Using the above constraint conditions, we can follow the same ideas as in (Ord4) case to have the estimate
\begin{align}
 |\mathscr{L}_{bcdef}  (w)|
 \lesssim  \  T^{\varepsilon}   \|v\|_{Y^{0}_{p',\Delta}}  \|u \|^4_{X^{1/2}_{p,\Delta}  },     \label{6linear113}
\end{align}
where $bcdef \in \{llllh_-, llll_-h_-, lll_-l_-h_-, ll_-l_-l_-h_-, l_-l_-l_-l_-h_- \}$ as in Table 8.

If $\lambda_0,...,\lambda_5$ satisfy any case of (Ord6)-(Ord10), we easily see that
$$
\lambda_0\geq \lambda_1, \ \ \lambda_2\geq \lambda_3,  \ \ \lambda_4\geq \lambda_5.
$$
By the frequency constraint condition (H2), one can conclude that
$$
\lambda_0\leq \lambda_1 +20, \ \  \lambda_2\leq \lambda_3 +20, \ \ \lambda_4\leq \lambda_5 +20.
$$
Hence, we have
$$
\lambda_0 \thickapprox  \lambda_1, \ \  \lambda_2 \thickapprox \lambda_3, \ \ \lambda_4 \thickapprox \lambda_5.
$$
We further have for the cases (Ord8) and (Ord10),
$$
\lambda_0 \thickapprox  \lambda_1 \thickapprox \lambda_2 \thickapprox \lambda_3 \thickapprox \lambda_4 \thickapprox \lambda_5.
$$
For the case (Ord7),
$$
 \lambda_2 \thickapprox \lambda_3 \thickapprox \lambda_4 \thickapprox \lambda_5.
$$
For the case (Ord9),
$$
\lambda_0 \thickapprox  \lambda_1 \thickapprox \lambda_2 \thickapprox \lambda_3.
$$
If $\lambda_0 \thickapprox  \lambda_1$, the summations on both $\lambda_0$ and $\lambda_1$ are the summation on $\lambda_0$ together with a finite summation on $\lambda_1$. So,  the proof in the case (Ord6) is easier than that of the case (Ord1). The details of the proof are omitted.
Up to now, we have finished the proof of Step 1.

{\bf Step 2.} We assume that $|\lambda_0|$ is the second largest one in $|\lambda_0|,...,|\lambda_5|$. We can assume, without loss of generality that $\lambda_0 \gg 1$.  There  exists $i\in \{1,...,5\}$  such that $|\lambda_i|= \max_{0\leq k\leq 5} |\lambda_k|$.  First, we point that this case is quite similar to that of $\lambda_0$ to be the largest one as in Step 1. Similar to Lemma \ref{nonestimate2h}, we have

 \begin{lem}\label{nonestimate2hstep2}
Let $0<c_k<C_k$, $k=1,2.$ For $w=(\bar{v},u,\bar{u},u,\bar{u},u)$, we write
\begin{align}
 \mathscr{L}_{2h}(w) = \sum_{ \lambda_0>0 } \langle \lambda_0\rangle^{1/2} \int_{\mathbb{R}\times [0,T]} \prod^2_{k=1} w^{(k)}_{[c_k\lambda_0, C_k\lambda_0]} \prod^5_{k=3} w^{(k)}_{[-C\lambda_0, C\lambda_0]} (x,t) dxdt, \label{6linear8astep2}
\end{align}
Let $p\in [4, \infty)$, $0<T<1$, $0<\varepsilon \ll 1$.    We have
\begin{align}
|\mathscr{L}_{2h}(w) | \lesssim T^{\varepsilon}  \|v\|_{Y^{0}_{p',\Delta}} \|u\|^5_{X^{1/2}_{p,\Delta}}.      \label{6linear31step2}
\end{align}
Substituting $ w^{(k)}_{[c_k\lambda_0, C_k\lambda_0]} $ by $ w^{(k)}_{[-C_k\lambda_0, -c_k\lambda_0]} $ for $k=1$ or $k=2$ in \eqref{6linear8astep2}, \eqref{6linear31step2} also holds.
 \end{lem}
The proof of Lemma \ref{nonestimate2hstep2} is almost the same as that of Lemma \ref{nonestimate2h} and we omit it.    If $|\lambda_i|$ is the largest one in $|\lambda_0|, ..., |\lambda_5|$, we have from the frequency constraint condition
$$
\lambda_0 \leq  |\lambda_i | \leq 5\lambda_0 + 10 \leq  20 \lambda_0.
$$
By Lemma \ref{nonestimate2hstep2} we can assume that $\lambda_k \in [-c\lambda_0, c\lambda_0]$ if $k\in \{1,...,5\}\setminus \{i\}$.  Namely, it suffice to consider the case that $\lambda_i$ is higher frequency and the other $\lambda_k$  ($k\neq i$) are lower frequency.
By symmetry we can further assume that
$$
\lambda_2\geq \lambda_4, \ \ \lambda_1\geq \lambda_3 \geq \lambda_5.
$$

{\it Case 1.} We consider the case $|\lambda_1|\vee |\lambda_3|\vee |\lambda_5| =\max_{0\leq k\leq 5} |\lambda_k|$. It follows that $|\lambda_1|  \vee |\lambda_5| =\max_{0\leq k\leq 5} |\lambda_k|$.  We further claim that
$|\lambda_1|    =\max_{0\leq k\leq 5} |\lambda_k|$. If not, then $-\lambda_5  =\max_{0\leq k\leq 5} |\lambda_k|$, which contradicts the frequency constraint condition.  It is easy to see that  $\lambda_1>0$.  Hence, we have
$$
\lambda_1\geq \lambda_0 \geq \lambda_k, \ \ k=2,...,5.
$$
The possible cases are the following
\begin{align*}
\lambda_1\geq \lambda_0 \geq \lambda_2 \geq \lambda_3\geq  \lambda_4 \geq \lambda_5, \tag{\rm 2Ord1}\\
 \lambda_1\geq \lambda_0 \geq \lambda_2 \geq \lambda_3\geq  \lambda_5 \geq \lambda_4, \tag{\rm 2Ord2}\\
 \lambda_1\geq \lambda_0 \geq \lambda_2 \geq \lambda_4\geq  \lambda_3 \geq \lambda_5, \tag{\rm 2Ord3}\\
 \lambda_1\geq \lambda_0 \geq \lambda_3 \geq \lambda_2\geq  \lambda_4 \geq \lambda_5, \tag{\rm 2Ord4}\\
\lambda_1\geq \lambda_0 \geq \lambda_3 \geq \lambda_2\geq  \lambda_5 \geq \lambda_4, \tag{\rm 2Ord5}\\
 \lambda_1\geq \lambda_0 \geq \lambda_3 \geq \lambda_5\geq  \lambda_2 \geq \lambda_4. \tag{\rm 2Ord6}
\end{align*}
We consider Case (2Ord1). Similar to the Step 1,  we have the 5 subcases as in Table 14.
\begin{table}[h]
\begin{center}
Case (2Ord1) : $\lambda_1 $ is maximal

\begin{tabular}{|c|c|c|c|c|c|}
\hline
 ${\rm Case}$   &  $\lambda_1\in $ & $\lambda_2\in $ & $\lambda_3\in $ &  $\lambda_4\in $ & $\lambda_5\in$  \\
\hline
$2hllll$  & $[\lambda_0, 20\lambda_0]$  & $l$  & $l$  & $l$   & $l$ \\
\hline
$2hllll_-$ & $[\lambda_0, 20\lambda_0]$  & $l$ & $l$  & $l$   &  $l_-$  \\
\hline
$2hlll_-l_-$ & $[\lambda_0, 20\lambda_0]$  & $l$  & $l$  & $l_- $   & $l_-$ \\
\hline
$2hll_-l_-l_-$ & $[\lambda_0, 20\lambda_0]$  & $l$  & $l_-$  & $l_- $   & $l_- $ \\
\hline
$2hl_-l_-l_-l_-$ & $[\lambda_0, 20\lambda_0]$  & $l_-$  & $l_-$  & $l_-$   & $l_-$ \\
\hline
\end{tabular}
\end{center}
\caption{$\lambda_1\geq\lambda_0\geq\lambda_2\geq \lambda_3 \geq \lambda_4 \geq \lambda_5$,  only one higher frequency in $\lambda_1,...,\lambda_5$.}
\end{table}
As in Table 14, we denote by $(\lambda_k) \in 2hllll$ that $\lambda_0,..,\lambda_5$ satisfy  (2Ord1) and $\lambda_1\in [\lambda_0, 20 \lambda_0]$, $\lambda_2,...,\lambda_5\in l$.   We consider the decomposition of $\lambda_k$:
$$
\lambda_1 \in [ \lambda_0,  20 \lambda_0] = \bigcup_{j_1\geq 0}( \lambda_0 + I_{j_1}); \ \  \lambda_k \in [0,c\lambda_0] = \bigcup_{j_k\geq 0} I_{j_k},  \ \  k=2,...,5.
$$
 We have the following constraint conditions on $j_1,...,j_5$:
\begin{align*}
& j_1\leq j_3+C, \ \ j_2\geq j_3\geq j_4 \geq j_5.
 \end{align*}
Let us write
\begin{align}
 \mathscr{L}_{2hllll}(w) = \sum_{(\lambda_k) \in 2hllll } \langle \lambda_0\rangle^{1/2} \int_{\mathbb{R}\times [0,T]} \prod^5_{k=0} w^{(k)}_{\lambda_k} (x,t) dxdt, \label{6linearstep2}
\end{align}
Applying the dyadic decomposition above and H\"older's inequality,
\begin{align}
 |\mathscr{L}_{2hllll}(w)| \lesssim & \sum_{\lambda_0; \ j_1\leq j_2+C; \  j_5\leq  j_4 \leq  j_3 \leq j_2  } \langle \lambda_0\rangle^{1/2} \|\overline{v}_{\lambda_0}\|_{L^\infty_{x,t}} \|u_{\lambda_0+I_{j_1}}\|_{L^4_{x,t\in [0,T]}}  \nonumber\\
  & \ \ \times   \|u_{I_{j_3}}\|_{L^4_{x,t\in [0,T]}}\|u_{I_{j_2}}\|_{L^4_{x,t\in [0,T]}}\|u_{I_{j_4}}\|_{L^4_{x,t\in [0,T]}} \|u_{I_{j_5}}\|_{L^\infty_{x,t}}.   \label{6linear48step2}
\end{align}
This estimate reduces to \eqref{6linear48}. Using the same way as in the estimates of \eqref{6linear48}, we have
\begin{align}
 |\mathscr{L}_{2hllll}(w)|  \lesssim &
 \  T^{ \varepsilon }   \|v \|_{Y^{0}_{p',\Delta}} \|u \|^5_{X^{1/2}_{p,\Delta}}.            \label{6linear51step2}
\end{align}
We can use similar way as in Step 1 and case (2Ord1) to handle the cases (2Ord2)-(2Ord6) and we omit the details.

{\it Case 2. } $|\lambda_2| \vee |\lambda_4|    =\max_{0\leq k\leq 5} |\lambda_k|$. By Lemma \ref{nonestimate2hstep2} we can assume that only one frequency $\lambda_2$ or $\lambda_4$ is localized in higher frequency intervals $[-20\lambda_0,-\lambda_0]\cup [\lambda_0, 20\lambda_0]$. In view of the frequency constraint condition we see that $-\lambda_4  =\max_{0\leq k\leq 5} |\lambda_k|$ and $\lambda_4 \in [-20\lambda_0, -\lambda_0]$.
Hence, we have
$$
\lambda_0 \geq \lambda_k \geq \lambda_4  , \ \ k=1,2,3,5.
$$
The possible cases are the following
\begin{align*}
\lambda_0 \geq\lambda_1 \geq \lambda_2 \geq \lambda_3\geq  \lambda_5 \geq \lambda_4, \tag{\rm 3Ord1}\\
\lambda_0 \geq \lambda_1\geq \lambda_3 \geq \lambda_2\geq  \lambda_5 \geq \lambda_4, \tag{\rm 3Ord2}\\
 \lambda_0 \geq \lambda_1\geq \lambda_3 \geq \lambda_5\geq  \lambda_2 \geq \lambda_4, \tag{\rm 3Ord3}\\
\lambda_0 \geq\lambda_2 \geq \lambda_1 \geq \lambda_3\geq  \lambda_5 \geq \lambda_4. \tag{\rm 3Ord4}
\end{align*}
For example, we consider the case (3Ord1).  As in Table 15, we have 5 subcases
\begin{table}[h]
\begin{center}
Case (3Ord1) : $-\lambda_4 $ is maximal

\begin{tabular}{|l|c|c|c|c|c|}
\hline
 ${\rm Case}$   &  $\lambda_1\in $ & $\lambda_2\in $ & $\lambda_3\in $ &  $\lambda_5\in $ & $\lambda_4\in$  \\
\hline
$3llllh_-$  & $l$  & $l$  & $l$   & $l$ & $[-20\lambda_0,-\lambda_0]$\\
\hline
$3llll_-h_-$ & $l$ & $l$  & $l$   &  $l_-$  & $[-20\lambda_0,-\lambda_0]$\\
\hline
$3lll_-l_-h_-$ & $l$  & $l$  & $l_- $   & $l_-$ & $[-20\lambda_0,-\lambda_0]$ \\
\hline
$3ll_-l_-l_-h_-$ & $l$  & $l_-$  & $l_- $   & $l_- $ & $[-20\lambda_0,-\lambda_0]$\\
\hline
$2l_-l_-l_-l_-h_-$ & $l_-$  & $l_-$  & $l_-$   & $l_-$ & $[-20\lambda_0,-\lambda_0]$ \\
\hline
\end{tabular}
\end{center}
\caption{$\lambda_0\geq\lambda_1\geq\lambda_2\geq \lambda_3 \geq \lambda_5 \geq \lambda_4$, only one higher frequency in $\lambda_1,...,\lambda_5$.}
\end{table}
These case are almost the same as those in Case (Ord2): $\lambda_4\in h_-$ and we do not perform the details of the proof.

{\bf Step 3.} We assume that $|\lambda_0|$ is the third largest one in $|\lambda_0|,...,|\lambda_5|$. We can assume, without loss of generality that $\lambda_0\geq 0$.  There  exists a rearrangement of $1,...,5$, which is denoted by $\pi(1),..., \pi(5) $  such that
$$
|\lambda_{\pi(1)}| \vee  |\lambda_{\pi(2)}| \geq \lambda_0 \geq  |\lambda_{\pi(3)}| \vee   |\lambda_{\pi(4)}| \vee |\lambda_{\pi(5)}|
 $$
For $\lambda_{\pi(1)}, \lambda_{\pi(2)}$, we apply the dyadic decomposition starting at $\pm \lambda_0$,
\begin{align}\label{decomp1}
 \lambda_{\pi(k)}  \in [\lambda_0 , \infty) = \bigcup_{j_{\pi(k)}\geq 0} (\lambda_0 + I_{j_{\pi(k)}}), \ \ or \   \lambda_{\pi(k)}  \in (-\infty, \lambda_0] = \bigcup_{j_{\pi(k)}\geq 0} (-\lambda_0 - I_{j_{\pi(k)}}).
\end{align}
For $\lambda_{\pi(3)}, \lambda_{\pi(4)}, \lambda_{\pi(5)}$, we adopt the dyadic decomposition starting at $ 0$:
\begin{align}\label{decomp2}
 \lambda_{\pi(k)}  \in [ 0 , \lambda_0] = \bigcup_{j_{\pi(k)}\geq 0}   I_{j_{\pi(k)}}, \ \ or \   \lambda_{\pi(k)}  \in  [-\lambda_0, 0] = \bigcup_{j_{\pi(k)}\geq 0} - I_{j_{\pi(k)}}.
\end{align}
For example, we consider the following case which is written as $(\lambda_k) \in  \mathbb{Z}^{6}_{+, 3}$:
$$
\lambda_0 > 0, \ \  \lambda_{\pi(1)},  \lambda_{\pi(2)} \in [\lambda_0 , \infty) , \ \  \lambda_{\pi(3)}, \lambda_{\pi(4)}, \lambda_{\pi(5)}  \in [ 0 , \lambda_0].
 $$
Let us write
\begin{align}
 \mathscr{L}^3 (w) = \sum_{(\lambda_k)\in  \mathbb{Z}^{6}_{+, 3} } \langle \lambda_0\rangle^{1/2} \int_{\mathbb{R}\times [0,T]} \prod^5_{k=0} w^{(k)}_{\lambda_k} (x,t) dxdt, \label{36linear8a}
\end{align}
Denote $j_{\max}= j_{\pi(3)} \vee j_{\pi(4)} \vee j_{\pi(5)}$,   $j_{\min}= j_{\pi(3)} \wedge j_{\pi(4)} \wedge j_{\pi(5)}$,  $j_{\rm med} \in \{ j_{\pi(3)} , j_{\pi(4)}, j_{\pi(5)} \} \setminus \{j_{\max}, j_{\min}\}$.  Applying H\"older's inequality, we have
\begin{align}
 |\mathscr{L}^3(w)| \lesssim & \sum_{\lambda_0; \ j_{\pi(1)}, \ j_{\pi(2)}, \ j_{\min}, j_{\rm med}, j_{\max} } \langle \lambda_0 \rangle^{1/2}  \|v_{\lambda_0}\|_{L^\infty_{x,t}}   \| u_{\lambda_0+I_{j_{\pi(1)}}} \|_{L^4_{x,t\in [0,T]}} \| u_{\lambda_0+I_{j_{\pi(2)}}} \|_{L^4_{x,t\in [0,T]}}\nonumber\\
  & \times  \|u_{I_{j_{\max}}} \|_{L^4_{x,t\in [0,T]}}   \|u_{I_{j_{\rm med}}} \|_{L^4_{x,t\in [0,T]}}   \| u_{I_{j_{\min}}} \|_{L^\infty_{x,t\in [0,T]}}      \label{36linear37}
\end{align}
In view of $V^2_\Delta \subset L^\infty_t L^2_x $, $\|v_{\lambda_0}\|_{L^\infty_x}  \lesssim \|v_{\lambda_0}\|_{L^2_x}$,
  by Corollary \ref{L4toXj} and Lemma \ref{L4toX}, one has that
\begin{align}
 |\mathscr{L}^3(w)| \lesssim & T^{4\varepsilon/3}  \sum_{ \lambda_0;  \ j_{\pi(1)}, \ j_{\pi(2)}, \  j_{\min}, j_{\rm med}, j_{\max} }
  \langle \lambda_0\rangle^{1/2 }  \|v_{\lambda_0}\|_{V^2_\Delta}   \langle \lambda_0 + 2^{j_{\pi(1)}}\rangle^{-1/2} \langle \lambda_0 +  2^{j_{\pi(2)}}\rangle^{-1/2} \nonumber\\
  & \ \  \times 2^{(j_{\pi(1)} +  j_{\pi(2)}) (1/4 -1/p +  \varepsilon)} 2^{(j_{\max}+j_{\rm med})(-1/4-1/p+ \varepsilon) }
 2^{j_{\min}(1/2-1/p)} \|u \|^5_{X^{1/2}_{p,\Delta}}  .       \label{36linear38}
\end{align}
Choosing $0<\varepsilon \leq 1/4p$, and making the summation on $j_{\min}$, $j_{\rm med}$ and $j_{\max}$ in order, one obtain that
\begin{align}
 |\mathscr{L}^3(w)| \lesssim & T^{4\varepsilon/3}  \sum_{ \lambda_0;  \ j_{\pi(1)}, \ j_{\pi(2)} }
  \langle \lambda_0\rangle^{1/2 }  \|v_{\lambda_0}\|_{V^2_\Delta}   \langle \lambda_0 + 2^{j_{\pi(1)}}\rangle^{-1/2} \langle \lambda_0 +  2^{j_{\pi(2)}}\rangle^{-1/2}    \nonumber\\
& \ \  \times 2^{(j_{\pi(1)} +  j_{\pi(2)}) (1/4 -1/p +  \varepsilon)}  \|u \|^5_{X^{1/2}_{p,\Delta}}.       \label{36linear39}
\end{align}
By splitting the summation $\sum_{  j_{\pi(1)} }$ into two parts $\sum_{  j_{\pi(1)} \leq \log^{\lambda_0}_2 }$ and  $\sum_{  j_{\pi(1)} > \log^{\lambda_0}_2 }$, one sees that
\begin{align}
 |\mathscr{L}^3(w)| \lesssim & T^{4\varepsilon/3}  \sum_{ \lambda_0  }
  \langle \lambda_0\rangle^{-2/p+2\varepsilon }  \|v_{\lambda_0}\|_{V^2_\Delta}      \|u \|^5_{X^{1/2}_{p,\Delta}}.       \label{36linear40}
\end{align}
Noticing that $\{\langle \lambda_0\rangle^{ -2/p +2\varepsilon}\} \in \ell^p$, by H\"older's inequality we have
\begin{align}
 |\mathscr{L}^3(w)|
\lesssim & T^{ \varepsilon }     \|v \|_{Y^{0}_{p',\Delta}}    \|u \|^5_{X^{1/2}_{p,\Delta}}. \ \ (0<\varepsilon \leq 1/4p)    \label{36linear41}
\end{align}
The proof in the above is also adapted to the other cases, say $(\lambda_k) \in  \mathbb{Z}^{6}_{\pm, 3}$:
$$
\lambda_0 > 0, \ \  \lambda_{\pi(1)}  \in [\lambda_0 , \infty) , \  \lambda_{\pi(2)} \in (-\infty, -\lambda_0],  \ \  \lambda_{\pi(3)} \in [-\lambda_0, 0], \ \lambda_{\pi(4)}, \lambda_{\pi(5)}  \in [ 0 , \lambda_0].
 $$
Since we did not use the frequency constraint condition, there is no essential difference between the cases $(\lambda_k) \in  \mathbb{Z}^{6}_{+, 3}$ and $(\lambda_k) \in  \mathbb{Z}^{6}_{\pm, 3}$. We omit the details of the proof for the other cases.  When $\lambda_0=0$, the proof is much easier than the case $\lambda_0>0$.

{\bf Step 4.} We assume that $|\lambda_0|$ is the fourth, or five largest one, or the minimal one in $|\lambda_0|,...,|\lambda_5|$. We can assume, without loss of generality that $\lambda_0\geq 0$. Let us consider the case $\lambda_0\gg 1$ be the fourth  largest one in $|\lambda_0|,...,|\lambda_5|$.  There  exists a rearrangement of $1,...,5$, which is denoted by $\pi(1),..., \pi(5) $  such that
$$
|\lambda_{\pi(1)}| \vee  |\lambda_{\pi(2)}| \vee |\lambda_{\pi(3)}| \geq \lambda_0 \geq  |\lambda_{\pi(4)}| \vee |\lambda_{\pi(5)}|
 $$
For $\lambda_{\pi(1)}, ..., \lambda_{\pi(3)}$, we apply the dyadic decomposition \eqref{decomp1} starting at $\pm \lambda_0$.
For $\lambda_{\pi(4)}, \lambda_{\pi(5)}$, we adopt the dyadic decomposition \eqref{decomp2} starting at $ 0$. Then we can use a similar way as in Step 3 to obtain the result, as desired.  The other cases can be shown along this line.  Up to now we have finished the proof of Lemma \ref{nonestimate}.


\end{document}